
%
%
%
\input amstex
\documentstyle{amsppt}
\NoBlackBoxes
\TagsOnRight
\voffset= 1cm
\hoffset=.5cm
\hsize=6.1in
\vsize=8.2in
\strut
\vskip5truemm
\font\small=cmr8

\def\det{\operatorname{det}}
\def\herm{\operatorname{Herm}}
\define\C{{\Bbb C}}
\define\R{{\Bbb R}}

\define\M{{\Cal E}}

\define\G{{\Cal G}}

\define\fs{{\Cal S}}

\define\ch{\operatorname{ch}}
\define\Tr{\operatorname{Tr}}

\define\Det{\operatorname{Det}}

\define\End{\operatorname{End}}

\define\GL{\operatorname{GL}}

\define\A{{\Cal A}}
\redefine\B{{\Cal B}}

\define\F{{\Cal F}}

\define\im{\operatorname{im}}

\define\cl{\operatorname{cl}}
\redefine\H{\Cal H}

\define\E{\Cal E}

\def\<{\langle}
\def\>{\rangle}

\define\pd#1#2{\dfrac{\partial#1}{\partial#2}}
\documentstyle{amsppt}

\topmatter

\title{Correspondences, von Neumann algebras}\\ {and holomorphic L$^2$
torsion} \endtitle

\author ALAN L. CAREY, MICHAEL S. FARBER and VARGHESE MATHAI\endauthor

\address Department of Pure Mathematics, University of Adelaide, Adelaide 5005,
Australia. \endaddress

\email acarey$\@$maths.adelaide.au.edu\endemail

\address School of Mathematical Sciences, Tel-Aviv University, Tel-Aviv 69978,
Israel. \endaddress

\email farber$\@$math.tau.ac.il\endemail

\address Department of Pure Mathematics, University of Adelaide, Adelaide 5005,
Australia. \endaddress

\email vmathai$\@$maths.adelaide.au.edu\endemail



\abstract
Given a holomorphic Hilbertian bundle on a compact complex manifold,
we introduce the notion of holomorphic $L^2$ torsion, which lies in the
determinant line of the twisted $L^2$ Dolbeault cohomology and represents a
volume element there. Here we utilise the theory of determinant lines of
Hilbertian modules over finite von Neumann algebras as developed in \cite{CFM}.
This specialises to the Ray-Singer-Quillen holomorphic torsion in the
finite dimensional case.
We compute a metric variation formula for the holomorphic $L^2$ torsion,
which shows that it is {\it not} in general independent of the choice of
Hermitian metrics on the complex manifold and on the holomorphic Hilbertian
bundle, which are needed to define it. We therefore initiate the theory of
correspondences of determinant lines, that enables us to define a relative
holomorphic $L^2$ torsion for a pair of flat Hilbertian bundles,
which we prove is independent of the choice of
Hermitian metrics on the complex manifold and on the flat Hilbertian
bundles.
\endabstract




\endtopmatter

\document
\heading{\S 0. Introduction}\endheading

Ray and Singer (cf. \cite{RS}) introduced the notion of holomorphic torsion of a
holomorphic bundle over a compact complex manifold. In  \cite{Q}, Quillen
viewed the
holomorphic torsion as an element in the real determinant line of the twisted
Dolbeault cohomology, or equivalently, as a metric in the dual of the
determinant line of the twisted Dolbeault cohomology. Since then there have
been many generalisations in the finite dimensional case, particularly by
Bismut, Freed, Gillet and Soule, \cite{BF, BGS}.

In this paper, we investigate generalisations
of aspects of this previous work to the case of infinite dimensional
representations of
the fundamental group.
Our approach is to introduce the concepts of holomorphic Hilbertian bundles
and of connections compatible with the holomorphic structure.
These bundles have fibres which are von Neumann algebra modules.
 We are able to define the
{\it determinant line bundle} of a holomorphic Hilbertian
bundle over a compact complex manifold, generalising the
construction of the determinant line of a finitely
generated Hilbertian module
that was developed in our earlier paper \cite{CFM}.
A nonzero element of the
determinant line bundle can be naturally
viewed as a volume form on the Hilbertian bundle.
This enables us to make sense of the notions of volume form and determinant
line bundle in this infinite dimensional and non-commutative situation.
Given an isomorphism of the determinant line bundles of holomorphic
Hilbertian bundles, we introduce the concept of a {\it correspondence}
between the  determinant lines of the twisted $L^2$ Dolbeault cohomologies.
This was previously studied in the finite dimensional situation in \cite{F}.

Restricting our attention to the class of manifolds studied
in [BFKM] (the so-called determinant or $D$-class
examples)
we then define the holomorphic $L^2$ torsion of a holomorphic Hilbertian
bundle;
it reduces to the classical constructions in the finite dimensional situation.
This new torsion invariant lives in the determinant line
of the twisted $L^2$ Dolbeault cohomology. Some key results in our
paper are a metric variation formula for the holomorphic $L^2$ torsion,
and the definition of a correspondence
between the determinant lines of the twisted $L^2$ Dolbeault cohomologies
for a pair of flat holomorphic Hilbertian bundles, and finally the definition
of a metric independent relative holomorphic $L^2$ torsion associated
to a correspondence between determinant line bundles of flat 
Hilbertian bundles.
To prove that a correspondence between determinant
line bundles of flat Hilbertian bundles is well defined, we need to prove
a generalised local index theorem for almost K\"ahler manifolds, and as a
consequence,
we give an alternate proof of Bismut's local index theorem for almost
K\"ahler manifolds
\cite{B}, where we use instead the methods of Donnelly \cite{D} and
Getzler \cite{Ge}.

The paper  is organized as follows. In the first section, we recall
some preliminary material
on Hilbertian modules over finite von Neumann algebras,
the canonical trace on the commutant of a finitely generated Hilbertian
module, the Fuglede-Kadison determinant on Hilbertian modules
and the construction of determinant lines for finitely generated
Hilbertian modules. Details of the material in this section can be
found in \cite{CFM}. In section \S 2, we define Hilbertian bundles
and connections on these. The definition of a connection is tricky in the
infinite dimensional context, and we use some fundamental theorems in von
Neumann
algebras to make sense of our definition. Then we define holomorphic
Hilbertian bundles
and connections compatible with the holomorphic structure
as well as Cauchy-Riemann operators on these.
In section \S 3, we study the properties of the zeta function associated to
holomorphic
Hilbertian bundles of $D$-class. In section \S 4, we define
the holomorphic $L^2$ torsion as an element in the determinant line
of reduced $L^2$ Dolbeaut cohomology. Here we also prove metric variation
formulae and we deduce that holomorphic $L^2$ torsion {\it does} depend
on the choices of Hermitian metrics on the compact complex manifold and
on the holomorphic Hilbertian bundle. However, in sections \S 5 and \S 6,
we give situations when a relative version of the holomorphic $L^2$ torsion
is indeed independent of the choice of metric.
In section \S 5, we are able to deduce
the following theorem (Theorem 5.5 in the text) from the variation formula: let
${\E}$ and ${\F}$ be two flat Hilbert bundles of $D$-class over
a compact Hermitian manifold $X$. Then one can define a relative holomorphic
$L^2$ torsion
$$
\rho_{\E,\F}^p\in \det(H^{p,\ast}(X,\E))\otimes \det(H^{p,\ast}(X,\F))^{-1}
$$
which is independent of the choice of Hermitian metric on $X$.
In section \S 6, we
define the notion of the determinant line bundle of a Hilbertian bundle and
also of
correspondences between determinant lines. The proof that a correspondence
is well defined, uses techniques of Bismut \cite{B}, Donnelly \cite{D} and
Getzler \cite{Ge} in their proof of the local index theorem in
different situations. Using the notion of a
correspondence of determinant line bundles,
we prove one of the main theorems in our paper (Theorem 6.12 in the text),
which can be briefly  stated as follows: let
${\E}$ and ${\F}$ be two flat Hilbertian bundles of $D$-class over
a compact almost K\"ahler manifold $X$
and $\varphi:\det(\E)\to\det(\F)$ be an isomorphism of the corresponding
determinant line bundles. Then one can define a relative holomorphic
$L^2$ torsion
$$
\rho_\varphi^p\in \det(H^{p,\ast}(X,\E))\otimes \det(H^{p,\ast}(X,\F))^{-1}.
$$
Using the correspondence defined by the isomorphism $\varphi$, we show
that the relative holomorphic
$L^2$ torsion $\rho_\varphi^p$ is independent of the choices of
Hermitian metrics on $\E$ and $\F$
and the choice of almost K\"ahler metric on $X$ which are needed to define it.
Recall that an almost K\"ahler manifold is a Hermitian manifold whose
"K\"ahler" 2-form $\omega$ is not necessarily closed, but satisfies
the weaker condition $\overline\partial\partial\omega = 0$. A result
of Gauduchon (cf. \cite{Gau}) asserts that every compact complex surface is
almost K\"ahler,
whereas there are many examples of complex surfaces which are not K\"ahler.
In  section \S 7, we give some examples of calculation of the
holomorphic $L^2$ torsion for locally symmetric spaces and Riemann surfaces.

\heading{\bf \S 1. Preliminaries}\endheading

This section contains some preliminary material from \cite{CFM}.

\subheading{1.0. Hilbertian modules over von Neumann algebras}

 Throughout the paper ${\A}$ will denote
 a finite  von Neumann algebra
with a fixed finite,
normal, and faithful trace $\tau:{\Cal A}\rightarrow \C$.  The involution
in ${\A}$ will be denoted $*$ while $\ell^{2}({\Cal A})$ denotes the completion
of  ${\A}$ in the norm derived from the inner product
$\tau(a^*b),\  a,b\in \A$.
A {\it Hilbert module} over ${\A}$
is a Hilbert space $M$
together with a continuous left ${\Cal A}$-module structure such that there
exists an isometric ${\Cal A}$-linear embedding of $M$ into
$\ell^{2}({\Cal A})\otimes H$, for some Hilbert space $H$.
(Note that this embedding is not
part of the structure.)  A Hilbert module $M$ is {\it finitely generated} if it
admits an imbedding as above with
finite dimensional $H$.
To introduce the notion of determinant line requires us
to  forget
the scalar product on $H$ but keep the topology and the ${\Cal A}$-action.

\subheading{1.1. Definition} A {\it  Hilbertian module}
is a topological vector space $M$
with continuous left ${\Cal A}$-action such that there exists a scalar product
$\langle\;,\;\rangle$ on $M$ which generates the topology of $M$ and such
that $M$ together
with $\langle\;,\;\rangle$ and with the ${\Cal A}$-action is a Hilbert
module. Any scalar product $\langle\;,\;\rangle$ on
$M$ with the above properties will be called {\it admissible}.

\subheading{1.2. Remarks and further definitions} The
choice of any other
 admissible scalar product $\langle\;,\;\rangle_1$ gives an isomorphic
 Hilbert module.  In fact there exists an operator
$
A:M\rightarrow M
$
such that
$$
\langle v,w\rangle_1=\langle Av,w\rangle\tag1
$$
for any $v,w\in M$.  The operator $A$ must be a self-adjoint, positive
linear homeomorphism (since the scalar products $\<\ ,\ \>$ and
$\<\ ,\ \>_1$ define the same topology),
which commutes with the $\A$-action. A
{\it finitely generated Hilbertian
module} is one for which the corresponding Hilbert module
 is finitely
generated. Finally, a
 {\it morphism} of Hilbertian modules is a continuous linear map
$f:M\to N$, commuting with the $\A$-action.
Note that the kernel of any morphism $f$ is again a Hilbertian module as is
the closure of the image $\cl(\im(f))$.

\subheading{1.3. The canonical trace on the commutant}

Any choice of an admissible scalar product $\langle\;,\;\rangle$ on $M$,
defines obviously a $*$-operator on $\B$ (by assigning to an operator
its adjoint) and turns $\B$ into a von Neumann algebra.
If we choose another admissible scalar product $\<\ ,\ \>_1$
on $M$ then the new involution will be given by
$$
f\mapsto \ A^{-1}f^*A\qquad \text{for}\quad f\in \B,\tag2
$$
where $A\in\B$ satisfies
$\<v,w\>_1\ =\ \<Av,w\>$ for $v,w\in M$.
The trace on the commutant may now be defined as in \cite{Di}
 and here will be denoted
$\Tr_{\tau}$. It
is finite, normal, and faithful. If $M$ and $N$ are two finitely
generated modules over $\A$, then the canonical traces $\Tr_\tau$ on
$\B(M)$, $\B(N)$ and on $\B(M\oplus N)$ are compatible
in the following sense:
$$\Tr_\tau \left(\matrix
A&B\\
C&D
\endmatrix \right)
= \Tr_\tau (A) \ +\ \Tr_\tau(D),\tag3$$
for all $A\in \B(M),\quad D\in\B(N)$ and any
morphisms $B:M\to N$, and $C:N\to M$.
Note that the {\it von Neumann dimension}
of a Hilbertian submodule $N$ of $M$ is defined
as $\dim_{\tau}(M)=\Tr_{\tau}(P_N)$
where $P_N$ is the orthogonal projection onto $N$.

\subheading{1.4. Fuglede-Kadison determinant for
Hilbertian modules}

Let $\GL(M)$ denote the group of all invertible elements
of the algebra $\B(M)$ equipped with
the norm topology. With this topology it is a Banach Lie group whose
 Lie algebra may be identified with the commutant $\B(M)$. The canonical
trace $\Tr_{\tau}$ on the commutant $\B(M)$
 is a homomorphism of the Lie algebra $\B(M)$ into
 $\C$ and by standard theorems, it defines a group homomorphism
of the universal covering group of $\GL(M)$ into $\C$. This approach
leads to following construction of the Fuglede-Kadison determinant,
compare \cite{HS}.

\proclaim{1.5. Theorem} There exists a function
$\Det_{\tau}: \GL(M)\to \R^{>0}$
(called the Fuglede-Kadison determinant) whose key properties are:
\roster
\item $\Det_{\tau}$ is a group homomorphism and
is continuous if
$\GL(M)$ is supplied with the norm topology;
\item If $A_t$ for $t\in [0,1]$
is a continuous piecewise smooth path in $\GL(M)$ then
$$
\log \lbrack\frac{\Det_{\tau}(A_1)}{\Det_{\tau}(A_0)}\rbrack\ =\
\int_0^1\Re\Tr_{\tau}\lbrack A_t^{-1}A_t^\prime\rbrack dt.\tag4
$$
Here $\Re$ denotes the real part and $A_t^\prime$ denotes the
derivative of $A_t$ with respect to $t$.
\item Let $M$ and $N$ be two finitely generated modules over $\A$, and
$A\in\GL(M)$ and $B\in\GL(N)$ two automorphisms, and
$\gamma:N\to M$ be a homomorphism. Then the map given by the matrix
$$\left(\matrix
A&\gamma\\
0&B
\endmatrix\right)$$
belongs to $\GL(M\oplus N)$ and
$$
\Det_\tau\left(\matrix
A&\gamma\\
0&B
\endmatrix\right)\ =\
\Det_\tau(A)\cdot\Det_\tau(B)\tag5
$$
\endroster
\endproclaim

Given an operator $A\in\GL(M)$, there is a continuous piecewise
smooth path
$A_t\in\GL(M)$ with $t\in [0,1]$ such that $A_0=I$ and $A_1=A$ (it is well
known that the group $\GL(M)$ is pathwise connected, cf. \cite{Di}).
Then from (4) we have the formula:
$$\log \Det_{\tau}(A)\ =\
\int_0^1\Re\Tr_{\tau}\lbrack A_t^{-1}A_t^\prime\rbrack dt.\tag6$$
This integral does not depend on the choice of the path.
 As an example  consider the following situation.
Suppose that a self-adjoint operator $A\in \GL(M)$
has spectral resolution
$$A\ =\ \int_0^\infty \lambda dE_\lambda\tag7$$
where $dE_\lambda$ is the spectral measure.
Then we can choose the path
$$A_t\ =\ t(A-I)\ +\ I,\quad t\in [0,1]$$
joining $A$ with $I$ inside $\GL(M)$. Applying (6) we obtain
$$
\log \Det_\tau(A)\ = \int_0^\infty \ln\lambda d\phi_\lambda\tag8
$$
where $\phi_\lambda = \Tr_\tau E_\lambda$ is the spectral density function.

\subheading{1.6. Operators of determinant class}.

Following \cite{BFKM} and \cite{CFM} we extend the
previous ideas to a wider class of operators. An operator $A$ as in (7)
is said to be $D-class$ ($D$ for determinant) if
$$\int_0^\infty \ln\lambda d\phi_\lambda > -\infty\tag9$$
A scalar product
$\langle v,w\rangle=\langle Av,w\rangle_1$
is said to be $D-admissible$ if $A$ is $D$-class and $\langle\ ,\ \rangle_1$
is any admissible scalar product.
The Fuglede-Kadison determinant extends to such operators
via the formula:
$$\Det_\tau(A)=\exp[\int_0^\infty \ln\lambda d\phi_\lambda].\tag10$$

\subheading{1.7. Determinant line of a Hilbertian module}

 For a Hilbertian module $M$ we defined in \cite{CFM} the
determinant line $\det(M)$
as a real vector
space generated by symbols
$\<\ ,\ \>$, one for any admissible scalar product on $M$, subject to the
following relations: for any pair  $\<\ ,\ \>_1$ and $\<\ ,\ \>_2$
of admissible scalar products on $M$ we require
$$\<\ ,\ \>_2\ =\ \sqrt {\Det_\tau(A)}^{\ -1} \cdot \<\ ,\ \>_1,\tag11$$
where $A\in \GL(M)\cap\B(M)$ is such that
$\<v,w\>_2\ =\ \<Av,w\>_1$
for all $v,w\in M$.
It is not difficult to see  that $\det(M)$ {\it is one-dimensional
generated by the
symbol $\<\ ,\ \>$ of any admissible scalar product on $M$}.
Note also, that the real line has {\it the canonical orientation}, since the
transition coefficient $\sqrt {\Det_\tau(A)}$ is always positive. Thus we
may speak of {\it positive and negative} elements of $\det(M)$.
We think of elements of $\det(M)$ as ``volume forms" on $M$.
If $M$ is trivial module, $M=0$, then we set $\det(M)=\R$, by definition.

 Given two finitely generated Hilbertian modules
$M$ and $N$ over $\A$,
with admissible scalar products $\<\ ,\ \>_M$ and $\<\ ,\ \>_N$
respectively, we may obviously define the
scalar product $\<\ ,\ \>_M \oplus \<\ ,\ \>_N$ on the direct sum. This
defines the isomorphism
$$\det(M)\otimes\det(N)\to\det(M\oplus N).\tag12$$
By property (5) of the Fuglede-Kadison determinant
it is easy to show that
 this homomorphism does not depend on the
choice of the metrics
$\<\ ,\ \>_M$ and $\<\ ,\ \>_N$ and
 preserves the orientations.
Note that, any isomorphism $f:M\to N$ between
finitely generated
Hilbertian modules induces canonically an orientation
preserving isomorphism of the determinant
lines
$f^\ast:\det(M)\to\det(N).$
Indeed, if $\<\ ,\ \>_M$ is an admissible scalar product on
$M$ then set
$$f^\ast(\<\ ,\ \>_M)= \<\ ,\ \>_N,\tag13$$
where $\<\ ,\ \>_N$ is the scalar product on $N$ given by
$\<v,w\>_N=\<f^{-1}(v),f^{-1}(w)\>_M$ for $v,w\in N$.
 This definition does not depend on the choice
of the scalar product $\<\ ,\ \>_M$ on $M$: if we have a different admissible
scalar product $\<\ ,\ \>_M^\prime$ on $M$, where $\<v,w\>_M^\prime =
\<A(v),w\>_M$ with $A\in\GL(M)$ then the induced scalar product on $N$
will be
$$\<v,w\>_N^\prime\ =\ \<(f^{-1}Af)v,w\>_N$$
and our statement follows from property (5) of the Fuglede-Kadison
determinant. Finally we note the
{\it functorial} property:
if $f:M\to N$ and $g:N\to L$ are two isomorphisms between finitely generated
Hilbertian modules then
$(g\circ f)^\ast\ =\ g^\ast\circ f^\ast.$

\proclaim{1.8. Proposition} If $f:M\to M$ is an automorphism of a
finitely generated Hilbertian module $M$, $f\in\GL(M)$, then the induced
homomorphism $f^\ast:\det(M)\to \det(M)$ coincides with the multiplication
by $\Det_\tau(f)\in \R^{>0}$. Furthermore any exact sequence
$$0\to M^\prime@>{\alpha}>>M@>{\beta}>>M^{\prime\prime}\to 0$$
of finitely generated Hilbertian modules determines canonically an isomorphism
$$\det(M^\prime)\otimes\det(M^{\prime\prime})\ \to \det(M),$$
which preserves the orientation of the determinant lines.
\endproclaim

\subheading{1.9. Extension to $D$-admissible scalar products}

Any $D$-admissible scalar product
 determines a
 non-zero element of the determinant
line det($M$) namely
$\Det_\tau(A)^{-1/2} \langle\ ,\ \rangle_1.$
A $D-admissible$ isomorphism $f:M\to N$ is one
for which the inner product $\<v,w\>_M=\<f(v),f(w)\>_N$
on $M$ is $D$-admissible for some and hence any admissible
inner product on $N$. Proposition 1.9 extends to
$D$-admissible isomorphisms and to the obvious notion of
$D$-admissible exact sequence.

\heading{\bf \S 2. Holomorphic Hilbertian $\A$-bundles Bundles and
$\A$-linear Connections}\endheading

In this section, we define Hilbertian $\A$-bundles and
$\A$-linear connections on these. The definition of ($\A$-linear) connection
is tricky in the infinite dimensional case, if one wants to be able to
horizontally lift curves. We use some fundamental theorems in von Neumann
algebras to make sense of our definition.
We also define holomorphic Hilbertian $\A$-bundles bundles and
holomorphic $\A$-linear connections on these.

\subheading{2.1. Hilbertian $\A$-bundles}
A Hilbertian $\A$-bundle with fibre $M$ over $X$ is given by the
following data.
\roster
\item $p:\E\to X$ a smooth bundle of topological vector spaces,
possibly infinite dimensional, such that each fibre $p^{-1}(x),\ x\in X$ is
a separable Hilbertian space (cf.\cite{Lang}).
\item There is a smooth fibrewise action
$\A\times \E\to\E$ which endows each fibre
$p^{-1}(x),\ x\in X$ with a Hilbertian $\A$-module structure, such
that for all $x\in X, p^{-1}(x)$ is isomorphic to $M$ as Hilbertian
$\A$-modules.
\item There is a local trivializing cover of $p : \E\to X$ which
intertwines the $\A$-actions.  More precisely, there is an open
cover $\{U_\alpha\}$ of $X$ such that for each $\alpha$, there is a smooth
isomorphism
$$
   \tau_\alpha : p^{-1}(U_\alpha) \to U_\alpha\times M
$$
which intertwines the $\A$-actions on $p^{-1}(U_\alpha) \subset
\E$ and on $U_\alpha\times M$, and such that ${\operatorname{pr}_1}\circ
\tau_\alpha = p$,
where $\operatorname{pr}_1 : U_\alpha\times M\to U_\alpha$ denotes the
projection onto the first factor.  The restriction of $\tau_\alpha$
$$
   \tau_\alpha : p^{-1}(x) \to \{x\}\times M
$$
is the isomorphism of Hilbertian $\A$-modules $\forall x\in
U_\alpha$, as given in (2).
\endroster

\subheading{2.2. Remarks}
If $\{U_\alpha\}$ is a trivializing open cover of $p : \E\to X$,
then the isomorphisms
$$
   \tau_\beta \circ \tau_\alpha^{-1} : (U_\alpha \cap U_\beta) \times M
      \to (U_\alpha\cap U_\beta)\times M
$$
are of the form $\tau_\beta\circ\tau_\alpha^{-1} = (\operatorname{id},
g_{\alpha\beta})$ where $g_{\alpha\beta} : U_\alpha\cap U_\beta \to
\operatorname{GL}(M)$ are smooth maps and are called the
transition functions of $p : \E\to X$, and they satisfy the
cocycle identity
$$
   g_{\alpha\beta}g_{\beta\gamma}g_{\gamma\alpha}=1\quad
      \forall \alpha,\beta,\gamma.
$$
Now suppose that $\{U_\alpha\}_\alpha$ is an open cover
of $X$, and on each
intersection $U_\alpha\cap U_\beta$, we are given smooth maps
$$
   g_{\alpha\beta} : U_\alpha \cap U_\beta \to
      \operatorname{GL}(M)
$$
satisfying $g_{\alpha\beta} g_{\beta\gamma} g_{\gamma\alpha}=1$ on
$U_\alpha\cap U_\beta \cap U_\gamma$ and $g_{\alpha\alpha}=1$ on
$U_\alpha$, then one can construct a Hilbertian $\A$-bundle $p :
\E\to X$ via the clutching construction viz, consider the disjoint
union $\tilde{\E}=\bigcup_\alpha(U_\alpha\times M)$ with the
product topology, and define the equivalence relation $\sim$ on
$\tilde{\E}$ by $(x,v)\sim (y,w)$ for $(x,v)\in U_\alpha\times M$
and $(y,w)\in U_\beta\times M$ if and only if $x=y$ and
$w=g_{\alpha\beta}(x)v$.  Then the quotient $\tilde{\E}/\sim =
\E\to X$ is easily checked to be a Hilbertian $\A$-bundle
over $X$

\subheading{2.3. Remarks}
This definition generalizes and is compatible with Breuer's definition of
Hilbert $\A$-bundles (cf.\cite{B}, \cite{BFKM}) and also with Lang's
definition \cite{Lang}, where the action of the von Neumann algebra is not
considered.
Actually Breuer \cite{B} considers von Neumann algebras $\A$ which are not
necessarily finite.

\subheading{2.4. Examples}
(a). It follows from Breuer's work (\cite{B}) that there are many examples of
Hilbertian $\A$-bundles, even in the case of simply connected
manifolds.  For example, on the 2-sphere $S^2$, the isomorphism classes of
Hilbertian $\A$-bundles with fibre $\ell^2(A)$, are in 1-1
correspondence with homotopy classes of maps from $S^1$ to
$\operatorname{GL}(\ell^2(A))$.  If $\A$ is a type $II_1$ factor, then by a
result of Araki, Smith and Smith \cite{ASS}, it follows that the
isomorphism classes
of Hilbertian $\A$-bundle
 over $S^2$ is isomorphic to
$\R$ (considered as a discrete group).

(b) Let $\E\to X$ be a Hilbertian $\A$-bundle over $X$. Then
$\Lambda^jT^*_\C X\otimes\E$ is also a Hilbertian $\A$-bundle over $X$, where
$\Lambda^jT^*_\C X$ denotes the jth exterior power of the complexified
cotangent bundle of $X$.
This can be seen as follows.
Let
$$
   g_{\alpha\beta} : U_\alpha \cap U_\beta \to GL(M)
$$
denote the transition functions of the Hilbertian $\A$-bundle $\E$ with fibre
$M$, and
$$
   g'_{\alpha\beta} : U_\alpha \cap U_\beta \to GL(r, \C)
$$
denote the transition functions of the $\C$ bundle $\Lambda^j T^*_\C X\to X$.
Then
$$
   g''_{\alpha\beta} : U_\alpha \cap U_\beta \to GL(\C^r\otimes M)
$$
denotes the transition functions of the Hilbertian $\A$-bundle
$\Lambda^jT^*_\C X\otimes\E$ with fibre $\C^r\otimes M$.

\subheading{2.5. Sections of Hilbertian $\A$-bundles}
A section of a Hilbertian $\A$-bundle $p : \E\to X$ is a
smooth map ${s} : X\to\E$ such that $p\circ {s}$ is
the identity map on $X$.  Let $\{ U_\alpha \}_\alpha$ be a local trivialization
of $p : \E \to X$.  Then a smooth section ${s}$ is given
on $U_\alpha$ by a smooth map $s_\alpha : U_\alpha\to M$.  On $U_\alpha\cap
U_\beta$ one has the relation ${s}_\alpha = g_{\alpha\beta}
{s}_\beta$.

\subheading{2.6. $\A$-linear connections on Hilbertian $\A$-bundles}
An $\A$-linear connection on a Hilbertian $\A$-bundle $p : \E\to X$ is
an $\A$-morphism
$$
   \nabla : \Omega^j (X,\E) \to \Omega^{j+1}(X,\E)
$$
such that for any $A\in\Omega^0
(X,\End_{\A}(\E))$ and
$w\in\Omega^j(X,\E)$, there is
$\nabla A\in\Omega^1(X,\End_{\A}(\E))$
such that
$$
   \nabla(Aw) - A(\nabla w) = (\nabla A)w.
$$
Here $\Omega^j (X,\E)$ denotes the space of smooth sections of the
Hilbertian $\A$-bundle $\Lambda^jT^*_\C X\otimes\E$, and
$\Omega^1(X,\End_{\A}(\E))$
denotes the space of smooth sections of the
Hilbertian $\A$-bundle $T^*_\C X\otimes\End_{\A}(\E)$

\subheading{2.7. Remarks}
Let $V$ be a vector field on $X$.  Then
$$
   \nabla_V A \in \Omega^0(X, \End_{\A}(\E))
$$

\proclaim{2.8. Proposition}
Let $\nabla, \nabla'$ be two connections on the Hilbertian
$\A$-bundle $p : \E\to X$ with fibre $M$.  Then
$$
   \nabla - \nabla' \in \Omega^1 (X,\End_{\A}\E)
$$
\endproclaim

\demo{Proof}
Let $V$ be a vector field on $X$.  Then $\delta_V = \nabla_V - \nabla'_V$
in $C^\infty(X)$ linear, and hence by (\cite{Lang}) is defined pointwise.
$(\delta_V)_x$ is a derivation on the von Neumann algebra
$\End_{\A}(\E_x)$.  Since $(\delta_V)_x$ is everywhere defined, by Lemma 3,
part III, chapter 9 of  (\cite{Dix}), $(\delta_V)_x$ is
bounded.  By Theorem 1, part III, chapter 9 of (\cite{Dix}), there is an element
$B_x(V)\in \End_{\A}(\E_x)$ such that $(\delta_V)_x =
ad B_x(V)$.  That is, $x\to ad B_x(V)$ is
smooth.  The remainder of the proof establishes that there is a
smooth choice $x\to \Tilde{\Tilde{B}}_x(V)$ such that $ad
\Tilde{\Tilde{B}}_x(V) = (\delta_V)_x$.  We first discuss the local
problem.

Let $U$ be an open subset of $X$ and $M$ be a Hilbertian $\A$-module.  Consider
the trivial bundle $U\times M\to U$ over $U$.  By Dixmier's result cited
 above, there
is a map
$$
   x\to ad B_x(V)\quad x\in U
$$
where $B_x(V)\in\End_{\A}(M)$ for all $x\in U$, such that
$$
   ad B_x(V) = (\nabla_V - \nabla'_V)_x
$$
since $\nabla, \nabla'$ are connections and $V$ is smooth, we deduce that
$x\to ad B_x(V)$ is smooth.  However, it isn't {\it a priori} clear
that one can choose $x\to B_x(V)$ to be smooth, as $B_x(V)$ is only defined
modulo the centre of the von Neumann algebra $\End_{\A}(M) = \B(M)$.
To complete the proof we need the next result.

\proclaim{2.9. Lemma}
Let $\A$ be a von Neumann algebra with centre $Z$.  Then there is a smooth
section $s : \A/Z \to \A$ to the natural projection $p : \A\to \A/Z$.
\endproclaim

\demo{Proof}
Let $Z\subset \A\subset B(\ell^2(\A))$, then since $Z$ is a type I von
Neumann algebra
and hence injective, there exists a
projection of norm 1, $ P:B(H) \to Z$ (\cite{HT}). Then $\A \cap ker P$ is a
complementary subspace
to $Z$ and one defines a section to the projection  $p : \A\to \A/Z$.
$$
   s : \A/Z \to \A\quad\text{as } s([v]) = (1-P) v.
$$
Then $s$ is smooth since it is linear.

More explicitly, given a subgroup $G$ of the unitaries in the commutant of
$Z$, $U(Z')$,
which is amenable and whose span is ultra weakly
dense in $Z'$, one can use the invariant mean on $G$ to average over the closure
of the orbit $\{uxu^*: u \in G\}$
and thus obtain a map $P$ so that $P(x)$ is this average for each $x$ and hence
commutes with every $u \in Z'$. That is, $P(x)$ is in $Z''=Z$. Such projections
are called Schwartz projections, according to Kadison. (cf. \cite{Ph}).
\enddemo\hfill$\square$

Returning now to the proof of proposition
2.8, we define the smooth map $\tilde{B}(V)$ by
$$
   \tilde{B}(V) = s \circ ad B(V),
$$
where $s: End_\A(M)/Z \to End_\A(M)$ is the section as in Lemma 2.9
(with $End_\A(M)$ replacing $\A$) . Then
clearly
$$
   \nabla_V = \nabla'_V = ad\tilde{B}(V),
$$
where $x\to\tilde{B}_x(V)$ is smooth.  This solves the problem locally.

Let $\E\to X$ be a Hilbertian bundle with fibre $M$, and
$\{U_\alpha\}$ be a trivialization of $\E\to X$.  We have seen
that on $U_\alpha$, there is a smooth section
$$
   x\to\tilde{B}_{\alpha, x}(V)\quad\text{for } x\in U_{\alpha}
$$
on $\E\big|_{U_\alpha\cap U_\beta}$, we can compare the 2 sections
obtained, $x\to\tilde{B}_{\alpha,x}(V) - \tilde{B}_{\beta, x}(V)\in Z$,
since $ad\tilde{B}_{\alpha, x}(V) = ad\tilde{B}_{\beta, x}(V)$.
Therefore we can define $\lambda_{\alpha\beta}(x) = \tilde{B}_{\alpha,x}(V)
- \tilde{B}_{\beta, x}(V)$ i.e.\ $\lambda_{\alpha\beta} : U_\alpha \cap
U_\beta \to Z$ is a Cech 1-cocycle with values in the sheaf of smooth $Z$
valued functions.  As $Z$ is contractable,
lemme 22 of \cite{DD} applies and so the 1st cohomology
with values in the sheaf of smooth $Z$ valued functions is trivial.
Therefore $\lambda_{\alpha\beta}$ is a coboundary i.e.\ there are
smooth maps
$$
   \varphi_\alpha : U_\alpha\to Z
$$
such that $\lambda_{\alpha\beta} = \varphi_\beta - \varphi_\alpha$.
Then $\{ x\to \tilde{B}_{\alpha, x}(V) + \varphi_{\alpha,x}(V)\}_\alpha$ is
a global section, since on $U_\alpha\cap U_\beta$, one has
$$
   \tilde{B}_{\alpha, x}(V) + \varphi_{\alpha, x}(V) =
      \tilde{B}_{\beta,x}(V) + \varphi_{\beta, x}(V),
$$
i.e.\ one gets a smooth section
$$
   X \to \End_{\A}(\E), \qquad x\to\Tilde{\Tilde{B}}_x(V)
$$
where $\Tilde{\Tilde{B}}_x(V) = \tilde{B}_{\alpha,x}(V) + \varphi_{\alpha,x}(V)$
for $x\in U_\alpha$.  It follows that $\Tilde{\Tilde{B}}\in \Omega^1
(X,\End_{\A}(\E)$.
\enddemo\hfill$\square$

Let $\nabla$ be a connection on $p : \E\to X$ and let $\{ U_\alpha
\}_\alpha$ be a trivialization of $p : \E\to X$.  Since
$\E\big|_{U_\alpha} \cong U_{\alpha}\times M$, one sees that the
 differential $d$ is a connection on $p : \E
\big|_{U_\alpha}\to U_\alpha$.  By proposition 2.8, $\nabla - d\in
\Omega^1(U_\alpha, \End_{\A}M)$ i.e.\ $\nabla=d+B_\alpha$ where
$B_\alpha\in\Omega^1(U_\alpha, \End_{\A}M)$.

On $U_\alpha \cap U_\beta$, one easily derives the relation
$$
B_\beta = g_{\alpha\beta}^{-1} B_\alpha g_{\alpha\beta} +
      g_{\alpha\beta}^{-1} dg_{\alpha\beta}.\tag14
$$
So a connection can also be thought of as a collection $\{d +
B_\alpha\}_\alpha$ where $B_\alpha\in\Omega^1(U_\alpha,
\End_{\A} M)$ and satisfying the relation (14) on the intersection.

\subheading{2.10. Parallel sections and horizontal lifts of curves}

Let $\nabla$ be a connection on $p : \E\to X$.  Let $p :
\E\to X$ be a Hilbertian $\A$-bundle and $I=[0,1]$ be the unit interval.
Let $\gamma : I\to X$ be a curve.  Let $\xi : I\to\E$ be a curve such that
$p_0\xi =
\gamma$.  Then $\xi$ is called a lift of $\gamma$.  $\xi$ is said to be a
{\it horizontal lift} of $\gamma$ if it is parallel along $\gamma$, that
is, if it satisfies the following equation,
$$
   \nabla_{\dot{\gamma}(t)} \xi(t)=0\quad \forall t\in I
$$
where dot denotes the derivative with respect to $t$.
In a local trivialization $U_\alpha$, the equation looks as,
$$
   \dot{\xi}(t) + B_\alpha(\dot{\gamma}(t))\xi(t) = 0\quad
      \forall t\in I\tag15
$$
where $\nabla = d+B_\alpha$ on $U_\alpha$ as before.  Since
$B_\alpha(\dot{\gamma}(t))$ is {\it bounded}, we use a theorem of
ordinary differential equations for Banach space valued functions (see prop
1.1, chapter IV in \cite{Lang}) to see that there is a unique solution to
equation (15) with initial
condition $\xi(0) = v\in M$.  It follows that a connection enables one to
lift curves horizontally.  This enables one to define a ``horizontal''
subbundle $\H$ of $T\E$, which is a complement to the
``vertical'' subbundle $p^*\E\subset T\E$.  This is how
[Lang] discusses connections on infinite dimensional vector bundles.
Conversely, given a choice of ``horizontal'' subbundle $\H$ of
$T\E$, one can define a ``covariant derivative'' (that is,
a connection) as follows.  By
hypothesis $T\E=\H\oplus p^*\E$.  Let
$\operatorname{pr}_2 : T\E\to p^*\E$ denote projection to
the 2nd factor and $\kappa : T\E \to \E$ be the
composition $p\circ \operatorname{pr}_2$ where $p :
p^*\E\to\E$.  Let $V$ be a vector field on $X$.  Define
$\nabla_V s = \kappa(D s(V))$ where $s :
X\to\E$ is a smooth section, and $D{s}$ is its
differential.  Then $\nabla$ locally has the form $\{ d + B_\alpha \}$ on a
trivialization $\{ U_\alpha \}$ of $p : \E \to X$, where $B_\alpha
\in \Omega^1(U_\alpha, \End_{\A}M)$ (see [Lang, Chapter IV, Section
3]) and it satisfies relation (14).  Therefore $\nabla$ defines
a connection on $p : \E \to X$ in the sense of 2.6.

\subheading{2.11. Holomorphic Hilbertian $\A$-bundles}

A Hilbertian $\A$-bundle $p : \E\to X$ with fibre $M$, is
said to be a holomorphic Hilbertian $\A$-bundle if the transition
functions of $p : \E\to X$,
$$
   g_{\alpha\beta} : U_\alpha \cap U_\beta \to
      \operatorname{GL}(M)
$$
are holomorphic maps.  We call $\{ U_\alpha \}_\alpha$ a holomorphic
trivialization of $p : \E\to X$.

\subheading{2.12 Remarks}
$\operatorname{GL}(M)$ is an open subset of a Banach space,
and so it is a complex manifold (of infinite dimension).

\subheading{2.13 Examples of holomorphic Hilbertian $\A$-bundles}
(a) By using the clutching construction again, we see that holomorphic
Hilbertian $\A$-bundles over $S^2$ correspond to holomorphic maps
$$
   g : A_\epsilon \to \operatorname{GL}(M)
$$
where $A_\epsilon = \{ z\in\C : 1 - \epsilon < |z| < 1+\epsilon \}$ is
an annulus, for some small $\epsilon>0$.  Therefore by 2.4, there
are many examples of holomorphic Hilbertian $\A$-bundles over
$S^2$.

(b) Let $p : \E\to X$ be a flat Hilbertian $\A$-bundle
over $X$, i.e.\ $M$ is a finitely generated $(\pi-\A)$ bimodule,
where $\varphi : \pi \to \operatorname{GL}(M)$ is the left
action of $\pi$ on $M$.  Then
$$
   \E = (M\times \tilde{X})/\sim \to X
$$
where $(v,x)\sim (\varphi(g)v, g.x)$ for $g\in\pi,\ v\in M$ and
$x\in\tilde{X}$. Let
$$
   g_{\alpha\beta} : U_\alpha \cap U_\beta \to \pi
$$
denote the transition functions of the universal cover $\tilde{X}$, which
is a principal $\pi$ bundle over $X$.  Here $\{ U_\alpha\}_\alpha$ forms an
open cover of $X$.  Since $\pi$ is a discrete group and $g_{\alpha\beta}$
is smooth, it follows that $g_{\alpha\beta}$ is locally constant, and therefore
holomorphic.  The transition functions of $\E$ are
$\varphi(g_{\alpha\beta})$, which again are locally constant, and therefore
holomorphic.

(c) Let $E\to X$ be a holomorphic $\C$-vector bundle over $X$ and
$\E\to X$ a flat Hilbertian $\A$-bundle over $X$.  Let
$$
   g_{\alpha\beta} : U_\alpha \cap U_\beta \to \operatorname{GL}(r,\C)
$$
denote the holomorphic transition functions where $\{ U_\alpha \}_\alpha$
form an open cover of $X$. Let
$$
   g'_{\alpha\beta} : U_\alpha \cap U_\beta \to
      \operatorname{GL}(M)
$$
denote the transition functions of the flat Hilbertian $\A$-bundle
$\E\to X$.  Since $\E\to X$ is flat, $g'_{\alpha\beta}$
are locally constant and this holomorphic (by the previous example).
Consider the new bundle whose transition functions are given by
$$
   g''_{\alpha\beta} \equiv g_{\alpha\beta} \otimes g'_{\alpha\beta}
      : U_\alpha \cap U_\beta \to \operatorname{GL}(\C^r\otimes M).
$$
Since the $g''_{\alpha\beta}$ are holomorphic, so is the new bundle which is
the tensor product bundle, and which is denoted by
$$
   E\otimes_{\C} \E \to X
$$
We have shown that it is a holomorphic Hilbertian $\A$-bundle over $X$,
with fibre
$\C^r \otimes M$.

\subheading{2.14. Holomorphic sections of holomorphic Hilbertian $\A$-bundles}

Let $p : \E\to X$ be a holomorphic Hilbertian
$\A$-bundle.  A section ${a} : X\to\E$ is said to
be a {\it holomorphic section} if in a holomorphic
local trivialization, $\{ U_\alpha \}_\alpha$,
the expression for ${s}$ in $U_\alpha$,
$$
   {s}_\alpha : U_\alpha \to M
$$
is a holomorphic map.  Note that $M$ is a Banach space, and therefore
a complex manifold.  On $U_\alpha\cap U_\beta$, one has the relation
$$
   {s}_\alpha = g_{\alpha\beta}{s}_\beta
$$
which is holomorphic, since $g_{\alpha\beta}$ is holomorphic.

\subheading{2.15. $\A$-linear Cauchy-Riemann operators}

Let $p : \E\to X$ be a holomorphic Hilbertian $\A$-bundle
over $X$.
With respect to the decomposition
$$
T^*_\C X = T^*X\otimes_\R \C = \big(T^{1,0}X\big)^*\oplus
\big(T^{0,1}X\big)^*, \tag16
$$
the space of smooth differential $j$-forms on $X$
with values in $\E$ decomposes as a direct sum of
spaces of smooth differential $(p,q)$-forms on $X$
with values in $\E$, where $p+q =j$.
This space, which is an $\A$ module, will be denoted by $\Omega^{p,q}(X,\E)$.

Then there is a unique operator
$$
   \bar{\partial} : \Omega^{p,q}(X,\E) \to \Omega^{p,q+1}
      (X,\E)
$$
which in any holomorphic trivialization of $p : \E\to X$, is equal
to
$$
   \bar{\partial} = \sum_{i=1}^n e(d\bar{z}^i)
      \frac{\partial}{\partial \bar{z}^i}
$$
where $e(d\bar{z}^i)$ denotes exterior multiplication by the $1$-form
$d\bar{z}^i$ and $n=\dim_{\C} X$. Note that $\bar{\partial}^2 = 0$.

\subheading{2.16. Holomorphic $\A$-linear connections}

Let $\nabla : \Omega^p(X,\E) \to \Omega^{p+1} (X,\E)$ be
an $\A$-linear connection on a holomorphic Hilbertian
$\A$-bundle $p : \E \to X$.  Then with respect to (16),
there is a decomposition
$$
   \nabla = \nabla' + \nabla''.
$$
Here
$$
   \nabla' : \Omega^{p,q} (X,\E) \to \Omega^{p+1,q}
       (X,\E)
$$
is an $\A$-morphism such that for $A \in
\Omega^0(X,\End_{\A} \E)$ and $w \in
\Omega^j (X,\E)$,
$$
   \nabla'(Aw) - A(\nabla' w) = (\nabla' A) w
$$
where $\nabla' A \in \Omega^{1,0} (X, \End_{\A}
\E)$ is the $(1,0)$ component of $\nabla A$, while
$$
   \nabla'' : \Omega^{p,q} (X, \E) \to \Omega^{p,q+1}(X,\E)
$$
is an $\A$-morphism such that
$$
   \nabla'' (Aw) - A(\nabla'' w) = (\nabla'' A)\, w
$$
where $\nabla'' A \in \Omega^{0,1} (X,\End_{\A}
\E)$ is the $(0,1)$ component of $\nabla A$.

An $\A$-linear connection $\nabla$ on a  holomorphic
Hilbertian $\A$-bundle $p : \E \to X$ is said to be a
{\it holomorphic $\A$-linear connection} if $\nabla'' =
\bar{\partial}$.  In this case, $(\nabla'')^2 = 0$.

Since every holomorphic Hilbertian $\A$-bundle has a
$\A$-linear Cauchy-Riemann operator, it follows that it also
has a holomorphic
$\A$-linear connection.

\subheading{2.17. Examples of holomorphic $\A$-linear connections}

(a) Let $\E\to X$ be a flat Hilbertian $\A$-bundle. Then $\E$
has a {\it canonical flat
$\A$-linear connection} $\nabla$ given by the de Rham exterior derivative,
where we identify the space of smooth differential $j$-forms on $X$
with values in $\E$,  denoted $\Omega^{j}(X,\E)$,
as $\pi$-invariant differential forms in
$M\otimes_\C \Omega^{j}(\widetilde X)$. Here
$M\otimes_\C \Omega^{j}(\widetilde X)$ has the diagonal action.
 (See \cite{CFM} for more details).
Since the de Rham differential $d = \bar{\partial} + {\partial}$,
it is a {\it canonical flat holomorphic $\A$-linear connection}.

(b) Let $E \to X$ be a holomorphic $\C$-vector bundle over $X$, and
$\E \to X$ a flat Hilbertian $\A$-bundle over $X$.  Then
we have seen that $E \otimes_{\C} \E \to X$ is a
holomorphic Hilbertian $\A$-bundle over $X$, with fibre
$\C^r \otimes M$.  Let $\tilde{\nabla}$ be a holomorphic connection
on $E \to X$, and let $\Tilde{\Tilde{\nabla}}$ be the canonical flat
$\A$-linear connection on $\E \to X$.  Then $\nabla =
\tilde{\nabla} \otimes 1  + 1 \otimes \Tilde{\Tilde{\nabla}}$ is easily
checked to yield a holomorphic $\A$-linear connection on the
holomorphic Hilbertian $\A$-bundle $E \otimes_{\C}
\E \to X$.

\heading{\bf \S 3. Zeta functions and $D$-class bundles}\endheading

We now have most of the notation and preliminary results we need
to generalize the classical construction of the holomorphic torsion
of D.B.Ray and I.M.Singer \cite{RS}
to the infinite dimensional case.
This section generalizes \cite{BFKM} and \cite{CFM} for the notion of a
$D$-class
holomorphic Hilbertian bundle and the definition of zeta-functions
for complexes of such bundles.

\subheading{3.1 Hermitian metrics, Hilbert $\A$ bundles, L$^2$ scalar
products and the canonical
holomorphic (Hermitian) $\A$-linear connection}
{\it A Hermitian metric $h$} on a Hilbertian $\A$-bundle $p:\E\to X$
is a smooth family of admissible scalar products on the fibers. Any
Hermitian metric on $p:\E\to X$ defines a wedge product
$$
\wedge : \Omega^{p,q}(X,\E)\otimes \Omega^{r,s}(X,\E)\rightarrow
\Omega^{p+r,q+s}(X)
$$
similar to the finite dimensional case.

Let $p:\E\to X$ be a {\it holomorphic} Hilbertian $\A$-bundle
and $h$ be a Hermitian metric on $\E$. The Hermitian metric on $p:\E\to X$
determines a {\it canonical holomorphic $\A$-linear connection} on $\E$ as
follows.
Let $\nabla$ be a holomorphic $\A$-linear connection on $\E$ which
preserves the Hermitian metric $\E$, that is,
$$
d h(\xi,\eta) = h(\nabla\xi,\eta) + h(\xi,\nabla\eta)
$$
where $\xi$ and $\eta$ are smooth sections of $\E$. Equating forms of the
same type, one has
$$
\partial h(\xi,\eta) = h(\nabla'\xi,\eta) + h(\xi,\nabla''\eta)
$$
and
$$
\bar{\partial} h(\xi,\eta) = h(\nabla''\xi,\eta) + h(\xi,\nabla'\eta).
$$
Since $\nabla'' = \bar{\partial}$, we see that a choice of Hermitian
metric determines a holomorphic $\A$-linear connection, which is called the
{\it canonical holomorphic $\A$-linear connection}.

The Hermitian metric on $p:\E\to X$ together with a Hermitian metric
on $X$ determines a scalar product on $\Omega^{p,q}(X,\E)$ in the
standard way; namely, using the Hodge star operator
$$\ast:\Omega^{p,q}(X,\E)\to \Omega^{n-q, n-p}(X,\E)$$
one sets
$$(\omega,\omega^\prime)\ =\ \int_X \omega\wedge\ast\overline\omega^\prime $$

With this scalar product $\Omega^i(X,\E)$
becomes a pre-Hilbert space. Define the space of
$L^2$ differential ${p,q}$-forms on $X$ with coefficients in
$\E$, denoted $\Omega_{(2)}^{p,q}(X,\E)$, to be the Hilbert
space completion of $\Omega^{p,q}(X,\E)$. We will tend to ignore the
scalar product on $\Omega_{(2)}^{p,q}(X,\E)$ and view it as an
infinite Hilbertian $\A$ module.

\subheading{3.2 Reduced L$^2$ Dolbeault cohomology}
Given a holomorphic Hilbertian $\A$ bundle $p:\E\to X$ together with a
Hermitian metric on $\E$, one
defines the {\it reduced $L^2$ Dolbeault cohomology with coefficients
in $\E$} as the quotient
$$
H^{{p,q}}(X,\E)=\frac{\ker
{\nabla''}/\Omega^{{p,q}}_{(2)}(X,\E)}{\cl(\im\;{\nabla''}/
\Omega^{p,q-1}(X,\E))},
$$
where the Cauchy-Riemann operator ${\nabla''}$ is associated to the canonical
$\A$-linear connection $\nabla$ on $\E$. ${\nabla''}$ on $\E$ extends
to an unbounded, densely defined operator $\Omega_{(2)}^{p,q}(X,\M)\to
\Omega_{(2)}^{p,q+1}(X,\M)$.
Then $H^{p,q}(X,\E)$ is naturally defined as a Hilbertian module over $\A$.
It can also be considered as the cohomology of $X$ with coefficients in a
locally constant sheaf, determined by $\E$.

\subheading{3.3 Hodge decomposition}
The Laplacian $\square_{p,q}$ acting
on $L^2$ $\M$-valued $(p,q)$-forms on $X$ is defined to be
$$
\square_{p,q}={\nabla''} {\nabla''}^{*} + {\nabla''}^*
{\nabla''}:\Omega_{(2)}^{p,q}(X,\M)
\rightarrow\Omega_{(2)}^{p,q}(X,\M)
$$
where ${\nabla''}^{*}$ denotes the formal adjoint of ${\nabla''}$ with
respect to
the
$L^{2}$ scalar product on $\Omega_{(2)}^{p,q}(X,\M)$.
Note that by definition, the
Laplacian is a formally self-adjoint operator which is densely defined. We
also denote by $\square_{p,q}$ the self adjoint extension of the Laplacian.

Let ${\H}^{{p,q}}(X,\M)$ denote the closed subspace of $L^{2}$ harmonic
${p,q}$-forms with coefficients in $\M$, that is, the kernel of
$\square_{p,q}$. Note that ${\H}^{p,q}(X,\M)$ is a Hilbertian ${\A}$-module.
By elliptic regularity (cf. section 2, \cite{BFKM}), one sees that
${\H}^{p,q}(X,\M) \subset \Omega^{p,q}(X,\E)$, that is, every $L^2$ harmonic
$(p,q)$-form with coefficients in $\M$ is smooth.
Standard arguments then show that one has the following Hodge decomposition
(cf. \cite{D}; section 4, \cite{BFKM} and also section 3, \cite{GS})
$$
\Omega_{(2)}^{p,q}(X,\M) = {\H}^{p,q}(X,\M) \oplus \cl(\im\;{\nabla''}/
\Omega^{p,q-1}(X,\E)) \oplus \cl(\im\;{\nabla''}^*/
\Omega^{p,q+1}(X,\E)).
$$
Therefore it follows that the natural map
$
{\H}^{{p,q}}(X,\M)\rightarrow H^{{p,q}}(X,\M)
$
is an isomorphism Hilbertian ${\A}$-modules.
The corresponding $L^2$ Betti numbers are denoted by
$$
b^{{p,q}}(X,\M) = \dim_\tau \left( H^{{p,q}}(X,\M)\right).
$$

\subheading{3.4 Definition }
Let $\square_{p,q} = \int_0^\infty \lambda dE_{p,q}(\lambda)$ denote the
spectral
decomposition of the Laplacian. The {\it spectral density function} is
defined to be
$N_{p,q}(\lambda) = \Tr_\tau(E_{p,q}(\lambda))$ and the {\it theta function}
is defined to be
$\theta_{p,q}(t) = \int_0^\infty e^{-t\lambda} dN_{p,q}(\lambda)=
\Tr_\tau(e^{-t\square_{p,q}}) -
b^{{p,q}}(X,\M)$. Here we use the well known fact that the projection
$E_{p,q}(\lambda)$ and the heat operator $e^{-t\square_{p,q}}$ have smooth
Schwartz kernels which are smooth sections of a bundle over $X\times X$ with
fiber the commutant of $M$, cf. \cite{BFKM}, \cite{GS}, \cite{Luk}.
The symbol $\Tr_\tau$ denotes
application of the canonical trace
on the commutant to the restriction of the kernels to the
diagonal followed by integration over the manifold $X$. This is a trace;
it vanishes
on commutators of smoothing operators.
See also \cite{M},
\cite{L}  and \cite{GS} for the case of the flat bundle defined by the
regular representation of the fundamental group.

\subheading{3.5 Definition} A holomorphic Hilbertian $\A$-bundle
$\M\rightarrow X$
together with a choice of Hermitian metric $h$ on $\E$, is said to be  {\it
D-class} if
$$
   \int_0^1 \log (\lambda) dN_{p,q}(\lambda) >\,-\infty
$$
or equivalently
$$
\int_1^\infty t^{-1} \theta_{p,q}(t) dt < \infty
$$
for all ${p,q}=0,....,n$.
Note that the
$D$-class property of a holomorphic Hilbertian $\A$ bundle does
not depend on the choice of metrics $g$ on $X$ and $h$ on $\E$.

For the most of the paper, we make the assumption that the holomorphic
Hilbertian $\A$-bundle
$\M\rightarrow X$ is D-class. Under this assumption, we
will next define and study the zeta function of the Laplacian $\square_{p,q}$
acting on $\E$ valued $L^2$
differential forms on $X$.

\subheading{3.6 Definition} For $\lambda>0$
the {\it zeta function of the Laplacian} $\square_{p,q}$ is defined on the
half-plane
$\Re(s)>n$ as
$$
\zeta_{{p,q}}(s,\lambda, \E)  = \frac{1}{\Gamma(s)}
\int_{0}^{\infty}t^{s-1}e^{-\lambda t} \theta_{{p,q}}(t)dt.
\tag17
$$

\proclaim{3.7 Lemma }
$\zeta_{{p,q}}(s,\lambda,\E)$ is a holomorphic function in the half-plane $
\Re(s)>n$ (where $n=\dim_\C X$) and has a meromorphic continuation to
$\C$ with no pole at $s=0$. If we assume that the holomorphic
Hilbertian $\A$-bundle $\M\rightarrow X$ is D-class then
$\lim_{\lambda\rightarrow 0} \zeta'_{{p,q}}(0,\lambda, \E)$
exists (where the prime denotes differentiation with respect to $s$)
\endproclaim

\demo{Proof} There is an asymptotic expansion as
$t\rightarrow 0^{+}$ of the trace of the heat kernel
$\Tr_\tau(e^{-t\square_{p,q}})$ (cf. \cite{BFKM} and chapter 13 \cite{R}),
$$
\Tr_\tau(e^{-t\square_{p,q}})\sim t^{-n}\sum_{i=0}^{\infty}t^{i}c_{i,{p,q}}
\tag18
$$
In particular, $\Tr_\tau(e^{-t\square_{p,q}})\leq C  t^{-n}$ for $
0<t\leq 1$. From this we deduce that $\zeta_{{p,q}}(s,\lambda,\E)$ is well
defined on the
half-plane $\Re(s)>n$ and it is holomorphic there.
The meromorphic continuation of $\zeta_{{p,q}}(s,\lambda,\E)$
to the half-plane $
\Re(s)>n-N$ is obtained by considering the first $N$ terms of the
small time asymptotic expansion (18) of $\Tr_\tau(e^{-t\square_{p,q}})$,
$$
\align
\zeta_{{p,q}}(s,\lambda,\E)\ & =
-\;\sum_j\frac{b^{{p,q}}(X,\E)(-\lambda)^j}{(s+j)j!}+\frac{1}{\Gamma(s)}
\left[
\sum_{0\leq i+j\leq N}  \frac{(-\lambda)^jc_{i,{p,q}}}{(s+i+j-{n})j!}
+ R_{N}(s,\lambda) \right]\\
& +\frac{1}{\Gamma(s)}\int_1^\infty t^{s-1}
\theta_{p,q}(t)e^{-t\lambda} dt
\tag19
\endalign
$$
where $R_{N}(s,\lambda)$ is holomorphic in the half plane $\Re(s)>n-N$
with a meromorphic extension to a neighbourhood of $s=0$.
Since the Gamma function has a simple pole at $s=0$, we observe
that the meromorphic continuation of $\zeta_{{p,q}}
(s,\lambda,\E)$ has no pole at $s=0$. The last part of the lemma
now follows
cf  \cite{BFKM}.

\enddemo\hfill$\square$

Let $ \zeta'_{{p,q}}(0,0, \E)=
\lim_{\lambda\rightarrow 0} \zeta'_{{p,q}}(0,\lambda, \E)$.
The following corollary is clear from (19).

\proclaim{3.8 Corollary }
One has
$$
\zeta_{{p,q}}(0,0,\E) = - b^{{p,q}}(X,\M)+ c_{n,{p,q}}
$$
where $c_{n,{p,q}}$ is the ${n}$-th coefficient in the small time
asymptotic expansion of the theta function, cf. (18).
\endproclaim

\heading{\bf \S 4. Holomorphic $L^2$-torsion}\endheading

In this section, we define and study the generalization of Ray-Singer
holomorphic torsion
to the case of holomorphic Hilbertian $\A$-bundles.
{\it For the rest of the section, we make the assumption that the holomorphic
Hilbertian $\A$-bundle $\M\rightarrow X$ is D-class}.
Given a Hermitian manifold $X$ and a metric on a
holomorphic Hilbertian $\A$-bundle $\M$ over $X$ with fibre a Hilbertian ${\A}$
module $M$, the holomorphic $L^2$ torsion $\rho_\E^p$ defined in this section is
a {\it positive element of the determinant line}
$$
\det(H^{p,*}(X,\M)).
$$
We also prove a variational formula for the holomorphic $L^2$ torsion.

\subheading{4.1. The construction of holomorphic $L^2$ torsion} Let $(X,g)$ be
a compact, connected Hermitian manifold of complex dimension $n$ with
$\pi=\pi_{1}(X)$. Let $\E\to X$ be a holomorphic Hilbertian $\A$-bundle over $X$
with fibre $M$ and let $h$ be a Hermitian metric on $\M$. We assume that
$\E$ is of $D$-class.

As before,
let $H^{{p,q}}(X,\M)$
denote the $L^{2}$ cohomology groups of $X$ with coefficients in $\M$.
Then we know that $H^{{p,q}}(X,\M)$ is a Hilbertian ${\Cal A}$-module.
If ${\Cal H}^{{p,q}}(X,\M)$ denotes the space of $L^{2}$ harmonic
${p,q}$-forms with coefficients in $\M$, then it is a Hilbert
${\Cal A}$-module with the admissible scalar product induced from
$\Omega_{(2)}^{p,q}(X,\M)$.  By the Hodge theorem, the natural map
$$
{\Cal H}^{{p,q}}(X,\M)\rightarrow H^{{p,q}}(X,\M)
$$
is an isomorphism of Hilbertian ${\A}$-modules. Thus, we may identify
these modules via this isomorphism, or equivalently,
we may say that this isomorphism defines an admissible scalar product
on the reduced $L^2$ cohomology  $H^{{p,q}}(X,\M)$.
  These admissible scalar products on $H^{{p,q}}(X,\M)$ for all ${p,q}$,
determine elements of the determinant lines $\det(H^{{p,q}}(X,\M))$
 and thus, their product in
$$\det(H^{p,\ast}(X,\M))
=\prod_{q=0}^n \det(H^{{p,q}}(X,\M))^{(-1)^q}$$
is defined. This last element we will denote $\rho^{\prime p}(g, h)$; the
notation emphasizing the dependence on the metrics $g$ and $h$.

Using the results of the previous section, we introduce the graded  zeta
function
$$
{\zeta^p}(s,\lambda,\M)=\sum_{q=0}^n (-1)^{q}q\zeta_{{p,q}}(s,\lambda,\M).
$$
It is a meromorphic function with no pole at $s=0$.
Note also that this
zeta-function depends on the choice of the trace $\tau$ and on the metrics
$g$ and $h$.

\subheading{4.2. Definition} Define the {\it holomorphic} $L^{2}$ {\it torsion}
to be the element of the determinant line
$$
\rho_{\E}^p(g,h)\in\det(H^{p,\ast}(X,\M)),\qquad
\rho_{\E}^p(g, h)=e^{\frac{1}{2} {\zeta^p}'(0,0,\M)}
\cdot
\rho^{\prime p}(g, h).
$$
where ${\zeta^p}'$ denotes the derivative with respect to $s$.
Thus, the holomorphic $L^2$ torsion is a volume form on the reduced
$L^2$ Dolbeault cohomology.

\subheading{4.3. Remarks}  1.  In the case when ${\Cal A}=\C$, we
arrive at the classical definition of the Ray-Singer-Quillen
metric on the determinant of the Dolbeault cohomology.

2. We will prove later in this section a metric
variation formula for the holomorphic $L^{2}$  torsion
as defined in 4.2. Using this, we prove that a relative version of the
holomorphic $L^2$ torsion is independent of the choice of Hermitian metric.

3. Assuming that the reduced $L^{2}$ Dolbeault cohomology $H^{p,*}(X,\M)$
vanishes, we can identify canonically the determinant line $\det({H}^{p,*} (X,
\M))$ with $\R$, and so the torsion $\rho_{\E}^p$ in this case is just
a number.

\subheading{4.4 Metric Variation Formulae}
Suppose that a holomorphic Hilbertian $\A$-bundle $\M\to X$ of $D$-class
is given. This property
does not depend on the choice of the metrics.
Consider a smooth 1-parameter family
of metrics $g_{u}$ on $X$ and $h_{u}$ on $\M$, where $u$ varies in an interval
$(-\epsilon,\epsilon)$.
Let $(,)_u$ denote the $L^2$
scalar product on $\Omega_{(2)}^{p,*}(X, \M)$ determined by $g_u$ and $h_u$.
This family determines an invertible, positive, self-adjoint bundle map
$A_u:\M\to \M$ which is uniquely determined by the relation
$$
(\omega,\omega')_u = (A_u \omega,\omega')_0
$$
for $\omega,\omega'\in \Omega_{(2)}^{p,*}(X, \M)$; it depends smoothly on $u$.

Let $\nabla$ be the canonical $\A$-linear connection on $\E$.
Define the operator
$$
D_{u}={\nabla''}+{\nabla''}^{*}_{u}:
\Omega_{(2)}^{p,*}(X,\M)\rightarrow\Omega_{(2)}^{p,*}(X,\M)
$$
where ${\nabla''}^{*}_{u}$ denotes the formal adjoint of ${\nabla''}$ with
respect to the
$L^{2}$ scalar product $(,)_u$ on $\Omega_{(2)}^{p,*}(X,\M)$.  Then
${\nabla''}^{*}_{u}=
A_u^{-1} {\nabla''}^*_0 A_u$ acting on $\Omega^{p,*}_{(2)}(X,\M)$.
Denote $Z_u = A_u^{-1}{\dot A_u}$, where the dot means the derivative with
respect to $u$.

As in 4.1, let $\zeta_u^p(s,\lambda,\M)$ denote the graded zeta function with
respect to the
metrics $g_u,\,h_u$. The scalar product $(,)_u$ induces a scalar
product on the space of harmonic
forms ${\H}^{p,*}_u(X,\M)$, and via the canonical isomorphism
${\H}^{p,*}_u(X,\M)\to
H^{p,*}(X,\M)$, it induces an admissible scalar product on the
reduced $L^2$ cohomology
$H^{p,*}(X,\M)$. Let $\rho'(u)$ denote the class
in $\det(H^{p,*}(X,\M))$ of this scalar product.  Then the
holomorphic $L^2$ torsion
with respect to the metrics $g_u,\,h_u$ is given, as in 4.2, by
$$
\rho_{\E}^p(u)=e^{{1\over 2}{\zeta^p}'_{u}(0,0,\M)}\rho^{\prime p}(u)\in
\det(H^{p,*}(X,\M)), $$
where ${\zeta^p}'$ means the derivative with respect to $s$.

\proclaim{4.5. Theorem} Let $\E\to X$ be a holomorphic Hilbert bundle of
$D$-class.
Then in the notation above, $u\mapsto \rho_{\E}^p(u)$ is a
smooth map and one has
$$
\frac{\partial}{\partial u}\rho_{\E}^p(u) = c_\E^p(u)  \rho_{\E}^p(u),
$$
where $c_\E^p(u)\in \R$ (cf. (24)) is a local
term.
\endproclaim

The proof of this theorem
will follow from two propositions which we will prove in this section.

Let $P_p(u)$ denote the orthogonal
projection from $\Omega_{(2)}^{p,*}(X, \M)$ onto $\ker D^{2}_{u}$ and
$\Tr^s_\tau(.)$ denote the graded trace,
that is the alternating sum of the von Neumann traces $\Tr_\tau$
on operators on $\Omega_{(2)}^{p,*}(X, \M)$ having smooth Schwartz kernels.

\proclaim{4.6. Proposition} Let $\E\to X$ be a holomorphic Hilbert bundle
of $D$-class.
Then in the notation above, one has
$$
\frac{\partial}{\partial u}\;{\zeta^p}'_{u}(0,0,\M)=\;\Tr^s_\tau
(Z_u P_p(u))-2 c_\E^p(u)
$$
where $c_\E^p(u)\in \R$ (cf. (24)) is a local term.
\endproclaim

\demo{Proof}  We consider the function
$$
F(u,\lambda,s) =\sum_{q=0}^{n}(-1)^{q}q\int^{\infty}_{0}t^{s-1}e^{-t\lambda}
\Tr_\tau(e^{-t
\square_{p,q}(u)}-P_{p,q}(u))dt
$$
which is defined on the half-plane $\Re(s)>n$ and is holomorphic there.
As in (18), one has for each $u$, the small time asymptotic expansion
of the heat kernel,
$$
\Tr_\tau(e^{-t\square_{{p,q}}(u)})\sim\sum_{k=0}^{\infty}c_{k,{p,q}}(u)t^{-n
+k}\quad.
\tag20
$$
Using (20), we see that $F(u,\lambda,s)$ has a meromorphic continuation to
$\C$ with
no pole at $s=0$.
This assertion is analogous to that in Lemma 2.8, and is proved by
an easy modification of that proof.

If we know that $u\to F(u,\lambda,s)$ is a smooth function then
$$
\frac{\partial}{\partial
u}\;{\zeta^p}'_{u}(0,0,\M)
= \left.\lim_{\lambda\rightarrow 0}\frac{\partial}{\partial
s}\Big(\frac{1}{\Gamma(s)}\frac{\partial}{\partial
u} F(u,\lambda,s)\Big)\right|_{s=0}
$$
by the $D$-class assumption.
Hence:
$$\frac{\partial}{\partial
u}\;{\zeta^p}'_{u}(0,0,\M)
= \left.\lim_{\lambda\rightarrow 0}\frac{\partial}{\partial u}\;F(u,\lambda,s)
\right|_{s=0}.
$$
Observing that $\Tr_\tau(P_{{p,q}}(u))=b^{{p,q}}(X,\M)$ is independent of $u$
we see that $u\to F(u,s)$ is smooth provided we can show
that $u\to \Tr_\tau(e^{-t\square_{{p,q}}(u)})$ is a smooth function. By
an application of Duhamel's principle, one has
$$
\frac{1}{u'-u}\Big(\Tr_\tau((e^{-{t\over 2}\square_{{p,q}}(u')}
-  e^{-{t\over 2}\square_{{p,q}}(u)})e^{-{t\over 2}\square_{{p,q}}(u)}) \Big)
$$
$$
= -\int_0^{t\over 2}
\Tr_\tau(e^{-s\square_{p,q}(u')} \frac{1}{u'-u}(\square_{p,q}(u') -
\square_{p,q}(u))
e^{-{t\over 2}\square_{{p,q}}(u)}e^{-({t\over 2}-s)\square_{{p,q}}(u)})ds.
\tag21
$$
Since $||\frac{1}{u'-u}(\square_{p,q}(u') - \square_{p,q}(u)) - {\dot
\square_{p,q}(u)})
e^{-{t\over 2}\square_{{p,q}}(u)}||$ is $O(u'-u)$ as $u'\to u$, one sees
that the limit as $u'\to u$ of (21) exists and
$$
\align
\Tr_\tau\Big(\Big(\frac{\partial}{\partial u}e^{-{t\over
2}\square_{{p,q}}(u)}\Big)
e^{-{t\over 2}\square_{{p,q}}(u)}\Big)
& =-\int_0^{t\over 2}
\Tr_\tau(e^{-s\square_{p,q}(u)} {\dot \square_{p,q}(u)}
e^{-{t\over 2}\square_{{p,q}}(u)}e^{-({t\over 2}-s)\square_{{p,q}}(u)})ds\\
& = -{t\over 2} \Tr_\tau( {\dot \square_{p,q}(u)} e^{-t\square_{{p,q}}(u)}).
\endalign
$$
Therefore $u\to \Tr_\tau(e^{-t\square_{{p,q}}(u)})$ is a smooth function
(and hence so is $u\to F(u,s)$) and by the semigroup property of the
heat kernel, one has
$$
\align
\frac{\partial}{\partial u}\;\Tr_\tau(e^{-t\square_{{p,q}}(u)}-P_{{p,q}}(u))& =
\frac{\partial}{\partial u}\;\Tr_\tau(e^{-t\square_{{p,q}}(u)}) \\
&= 2 \Tr_\tau\Big(\Big(\frac{\partial}{\partial u}e^{-{t\over
2}\square_{{p,q}}(u)}\Big)
e^{-{t\over 2}\square_{{p,q}}(u)}\Big)\\
&= -t {\Tr_\tau}(\dot{\square}_{{p,q}}(u)e^{-t\square_{{p,q}}(u)}).
\endalign
$$
A calculation similar to \cite{RS}, page 152 yields
$$
\dot{\square}_{{p,q}}(u)=-Z_u {\nabla''}^{*}_{u}{\nabla''}+{\nabla''}^{*}_{u}Z_u
{\nabla''}
-\;{\nabla''} Z_u {\nabla''}^{*}_{u}+{\nabla''}\;{\nabla''}^{*}_{u} Z_u.
$$
Since ${\nabla''}\;\square_{{p,q}}(u)=\square_{p,q+1}(u){\nabla''}$ and
${\nabla''}^{*}_{u}\square_{{p,q}}(u)=
\square_{p,q-1}(u){\nabla''}^{*}_{u}$ and using the fact that $\Tr_\tau$ is a
trace,
one has
$$
\align
\Tr_\tau(\dot{\square}_{{p,q}}(u)e^{-t\square_{{p,q}}(u)}) \ & =
\Tr_\tau(Z_u {\nabla''}\;{\nabla''}_{u}^ {*}e^{-t\square_{{p,q}}(u)}) -\;
\Tr_\tau(Z_u {\nabla''}^{*}_{u}{\nabla''}\;e^{-t\square_{p,q-1}(u)}) \\
 & + \;\Tr_\tau(Z_u {\nabla''}\;{\nabla''}^{*}_{u}e^{-t\square_{p,q+1}(u)})
-\;\Tr_\tau(Z_u {\nabla''}^{*}_{u}{\nabla''}_{u}e^{-t\square_{{p,q}}(u)}).
\endalign
$$
So one sees that
$$
\align
\frac{\partial}{\partial
u}\;\sum_{q=0}^{n}(-1)^{q}q\Tr_\tau(e^{-t\square_{{p,q}}(u)}-P_{{p,q}}(u))
&= -t
\sum_{q=0}^{n}(-1)^{q}q\Tr_\tau(\dot{\square}_{{p,q}}(u)e^{-t\square_{{p,q}}
(u)})
\\
& = -t \sum_{q=0}^{n}(-1)^{q}q\Tr_\tau(Z_u
{\square}_{{p,q}}(u)e^{-t\square_{{p,q}}(u)})\\
& = t \frac{\partial}{\partial t}\sum_{q=0}^{n}(-1)^{q}\Tr_\tau(Z_u e^{-t
\square_{{p,q}}(u)}).
\endalign
$$
Using this, one sees that for $\Re(s)>n$,
$$
\frac{\partial}{\partial u}\;F(u,\lambda,s) =\sum_{q=0}^{n}(-1)^{q}
\int^{\infty}_{0}t^{s}e^{-t\lambda}
\;\frac{\partial}{\partial t}\;\Tr_\tau(Z_u
(e^{-t\square_{{p,q}}(u)}
-P_{{p,q}}(u))dt
\tag22$$
Since $Z_u$ is a bounded endomorphism, by a straightforward generalization
of lemma 1.7.7 in \cite{Gi},
there is a small time asymptotic expansion
$$
\Tr_\tau(Z_u e^{-t\square_{{p,q}}(u)})\sim\sum_{k=0}^{\infty}m_{k,{p,q}}
(u)t^{-n+k}.\tag23
$$
In particular, one has
$$
|\Tr_\tau(Z_u e^{-t\square_{{p,q}}(u)})|\leq ct^{-n}
$$
for all $0< t\leq 1$.  If $\Re(s)>n$, we can integrate
 the right-hand side of
(22) by parts to get
$$
\;\sum_{q=0}^{n}(-1)^{q+1}
\int^{\infty}_{0}(st^{s-1}-\lambda t^s)e^{-t\lambda}\Tr_\tau(Z_u
(e^{-t\square_{{p,q}}(u)}-P_{{p,q}}(u)))
dt
$$
By splitting the integral into two parts, one from 0 to 1
and the other from 1 to $\infty$ and using
(23) on the first integral together with  the observations above, one gets
the following
explicit expression
for the meromorphic continuation of
$\displaystyle\frac{\partial}{\partial u}\;F(u,s)$ to the half-plane
$\Re(s)>n-N$
$$
\align
\frac{\partial}{\partial u}\;F(u,s)\ & =
\sum_{q=0}^{n}(-1)^{q}\Tr_\tau(Z_u P_{{p,q}}(u))\frac{1}{\lambda^s}
\int_0^\lambda(st^{s-1}-t^s)e^{-t}dt\\
&+
\sum_{q=0}^{n}(-1)^{q+1}\sum_{0\leq k+r\leq N}
\frac{(-\lambda)^rm_{k,{p,q}}(u)}{r!}(\frac{s} {s-n+k+r}-
\frac{\lambda}{s-n+k+r+1}) + R_{N}(u,\lambda,s)
\endalign
$$
where $R_N(u,\lambda,s)$ is holomorphic
in a neighbourhood of zero. At $s=0$
we have
$$
R_N(u,\lambda,0)= \int^{\infty}_{1}\Tr_\tau(Z_u
(e^{-t\square_{{p,q}}(u)}-P_{{p,q}}(u)))
e^{-t\lambda}dt.
$$
Thus we have
$$
\align
\frac{\partial}{\partial u}\;{\zeta'}_{u}^p(0,\lambda,\M)\ & =
\frac{\partial}{\partial u}\;F(u,0) \\
&= \sum_{q=0}^{n}(-1)^{q+1}\Big(\sum_{k+r=n}\frac{(-\lambda)^r}{r!}(1-\lambda)
m_{k,{p,q}}(u) -\Tr_\tau(Z_u P_{{p,q}}(u))\Big)
+R_N(u,\lambda,0)
\endalign
$$
Hence
$${\zeta^p}'_{u}(0,0,\M)
=  \sum_{q=0}^{n}(-1)^{q}\Tr_\tau(Z_u P_{{p,q}}(u))- 2c_\E^p(u)
$$
where
$$
c_\E^p(u)   =
\;\frac{1}{2}\;\sum_{q=0}^{n}(-1)^{q}m_{n,{p,q}}(u) \tag24
$$
This completes the proof of the proposition.
\enddemo\hfill$\square$

The 1-parameter family of scalar products on $\Omega^{p,*}_{(2)}({X}, \M)$
which are induced
by the 1-parameter family of metrics on $X$ and $\M\to X$, defines an
inclusion isomorphism of Hilbertian modules
$$
I_u : {\Cal H}^{p,\ast}_{u}(X,\M)\rightarrow H^{p,\ast}(X,\M).
$$
Here ${\Cal H}^{p,q}_{u}(X,\M)$ denotes the kernel of $\square_{p,q}(u)$.
There is an induced isomorphism of determinant lines cf. (13) and the discussion
in the paragraph above it.
$$
I_u^* : \det(H^{p,\ast}(X,\M))\rightarrow \det({\Cal
H}^{p,\ast}_{u}(X,\M)).
$$
We first identify $H^{p,\ast}(X,\M)$ with ${\Cal H}^{p,\ast}_{0}(X,\M)$.
Then $I_u$ defines a 1-parameter family of admissible scalar products on
$H^{p,\ast}(X,\M)$, which we can write explicitly as follows:
$$
\langle \eta, \eta'\rangle_u = ( P(u)\eta,
P(u)\eta')_u = (A_u  P(u)\eta, P(u)\eta')_0
$$
where $\eta,\eta'$ are harmonic forms in
${\Cal H}^{p,\ast}_{0}(X,\M)$. The relation between these scalar products
in the determinant line $\det(H^{p,\ast}(X,\M))$ is given as in 1.8  and
(11), by
$$
\langle\ ,\  \rangle_u = \prod_{q=0}^n {{\Det}_{\tau'}
(P_{p,q}(u)^\dagger A_u P_{p,q}(u))}^{{(-1)^{q+1}\over 2}}
\langle\ ,\  \rangle_0.\tag25
$$
where $P_{p,q}(u)^\dagger$ denotes the adjoint of $P_{p,q}(u)$ with respect to
the
fixed admissible scalar product $\langle\ \ ,\  \rangle_0$
and $\Tr_{\tau'}(\cdot)$ is the trace on
$H^{{p,q}}({X}, \M)$.
Using the fact that $H^{p,q}( X, \M)$ is isomorphic
to a submodule of a free Hilbertian module as is
${\Cal H}^{{p,q}}_{u}(X, \M)$, it follows that
$\Tr_{\tau'}(.)$
is equal to $\Tr_\tau(P_{p,q}(u)\cdot P_{p,q}(u))$.  We begin with the following

\proclaim{4.7. Proposition} Let $\E\to X$ be a holomorphic Hilbert bundle
of $D$-class.
Then the function $u\to P_{p,q}(u)$ is smooth and
in the notation of 4.4 and 4.5, one has
$$
\frac{\partial}{\partial u} \rho^{\prime p}(u) =
-\frac{1}{2}\Tr^s_\tau( Z_u P_p(u)) \rho^{\prime p}(u).
$$
\endproclaim

\demo{Proof} We will first prove that $u\to P_{p,q}(u)$ is a smooth function.
First consider the Hodge decomposition in the $u$-metric in the context,
$$
\Omega^{{p,q}}_{(2)}({X}, \M) = {\H}^{p,q}_u(X,\E) \oplus \cl(\im {\nabla''})
\oplus
\cl(\im {\nabla''}^*_u)
$$
and let $\pi$ denote the projection onto $\cl(\im {\nabla''})$, which
does not depend on the $u$-metric.
Let $h \in {\H}^{p,q}_0(X,\E)$ be harmonic in the $u=0$ metric. We will
arrive at a formula for $h_u \equiv P_{p,q}(u) h$, from which which
the differentiability
of $u\to P_{p,q}(u)$ will be clear. Now define $r_u$ by the equation
$$
h_u = h + r_u.
$$
Since $h_u$ is harmonic in the $u$-metric, one has
${\nabla''}^*_u(h_u) = 0$. By the formula for ${\nabla''}^*_u$ in 4.4,
one sees that ${\nabla''}^*_0 A_u(h+r_u) = 0$.
Since ${\nabla''}^*_0$ is injective on $\cl(\im {\nabla''})$, one has that
$\pi(A_u(h+r_u)) = 0$. Since $B_u\equiv \pi\,A_u \pi : \cl(\im {\nabla''})\to
\cl(\im {\nabla''})$ is an isomorphism, one sees that
$r_u = - B_u^{-1}\pi A_u(h)$ and therefore
$$
h_u = h -  B_u^{-1}\pi A_u(h).
$$
Since $u\to A_u$ is smooth, it follows that $u\to B_u$ is smooth
and by the formula above, one concludes that $u\to P_{p,q}(u)$ is
also smooth.

Observe that
$$
P_{p,q}(u)^2 = P_{p,q}(u).
$$
Differentiating with respect to $u$, one has
$$
\dot{P_{p,q}(u)} = P_{p,q}(u)\dot{P_{p,q}(u)} +
\dot{P_{p,q}(u)} P_{p,q}(u).
$$
Therefore
$$
P_{p,q}(u)\dot{P_{p,q}(u)}P_{p,q}(u) = 0.
$$
Therefore
$$
\align
{\Tr_\tau}(\dot{P_{p,q}(u)})
&= 2{\Tr_\tau}(P_{p,q}(u)\dot{P_{p,q}(u)})\\
&=2{\Tr_\tau}(P_{p,q}(u)\dot{P_{p,q}(u)}P_{p,q}(u)) \\
& =0.
\endalign
$$
A similar argument shows that the projection $P_{p,q}^\dagger(u)$
also satisfies
$$
{\Tr_\tau}(\dot{P_{p,q}^\dagger (u)})  = 0
$$

By definition,
$\rho^{\prime p}(u) = \langle\  ,\  \rangle_u \in \det\Big( H^{p,*}(X,\M)\Big)$,
and therefore by differentiating the relation (25), one has
$$
\frac{\partial}{\partial u} \rho^{\prime p}(u) =
-{1\over 2} \Tr^s_{\tau'}( C_u^{-1}\frac{\partial}{\partial u} C_u)
\rho^{\prime p}(u)
$$
where $C_u \equiv P_p(u)^\dagger A_u  P_p(u)$ and
$\Tr^s_{\tau'}(.)$ denotes the graded von Neumann trace on
$H^{p,\bullet}( X, \M)$.
Therefore one sees that
$$
\align
\frac{\partial}{\partial u} \rho^{\prime p}(u)
& =-{1\over 2} \Tr^s_{\tau'}( Z_u P_p(u) + P_p(u) \dot{P_p(u)} +  P_p^\dagger(u)
\dot{P_p^\dagger(u)})\rho^{\prime p}(u)\\
& = -{1\over 2} \Tr^s_{\tau'}( Z_u P_p(u))\rho^{\prime p}(u) =
-{1\over 2} \Tr^s_{\tau}( Z_u P_p(u))\rho^{\prime p}(u).
\endalign
$$

\enddemo\hfill$\square$

\noindent{\bf Proof of Theorem 4.5}.
By Proposition 4.7, one calculates
$$
\align
\frac{\partial}{\partial u} \rho_{\E}^p (u) \ & =
 \frac{1}{2}e^{\frac{1}{2} {\zeta^p}'_{u}(0,0,\M)}  \frac{\partial}{\partial u}
{\zeta^p}'_{u}(0,0,\M)
\rho^{\prime p}(u)
+ e^{\frac{1}{2} {\zeta^p}'_{u}(0,0,\M)} \frac{\partial}
{\partial u} \rho^{\prime p}(u) \\
& = \frac{1}{2} \Big[ \frac{\partial}{\partial u} {\zeta^p}'_{u}(0,0,\M)
- \Tr^s_\tau( Z_u P_p(u)) \Big]
 e^{\frac{1}{2} {\zeta^p}'_{u}(0,0,\M)}
\rho^{\prime p}(u)\\
& = \frac{1}{2} \Big[ \frac{\partial}{\partial u} {\zeta^p}'_{u}(0,0,\M)
-  \Tr^s_\tau( Z_u P_p(u)) \Big] \rho_{\E}^p (u).
\endalign
$$
Therefore by Proposition 4.6, one has
$$
\frac{\partial}{\partial u} \rho_{\E}^p (u) =  c_\E^p(u)  \rho_{\E}^p (u)
$$
where $c_\E^p(u)\in \R$ is as in (24).
This completes the proof of the theorem.
\hfill$\square$

\heading{\bf \S 5. Flat Hilbert $\A$-bundles and Relative Holomorphic $L^2$
Torsion}
\endheading

In this section, we define the relative holomorphic
$L^2$ torsion with respect to a pair of {\it flat Hilbert} (unitary)
$\A$-bundles
$\E$ and $\F$, and we
prove that it is independent of the choice of Hermitian metric on the
complex manifold. Thus it can be viewed as an invariant
volume form on the reduced $L^2$ cohomology  $H^{p,*}(X,\E)\oplus
H^{p,*}(X,\F)^\prime$.
In section \S 6, we will prove the relative holomorphic
$L^2$ torsion with respect to a pair of {\it flat Hilbertian}  $\A$-bundles
$\E$ and $\F$, is independent of the choice of almost K\"ahler metric on an
almost K\"ahler manifold and on the choice of Hermitian metrics on $\E$ and
$\F$.

\subheading{5.1 Relative holomorphic $L^2$ torsion}
It follows from Theorem 4.5 that the holomorphic $L^2$ torsion is {\it not}
independent of the choice of metrics on the complex manifold and on the
flat Hilbertian
bundle. Therefore in order to obtain an invariant, we now consider the
{\it relative} holomorphic $L^2$ torsion for a pair of {\it unitary} flat
Hilbertian bundles over a complex manifold. In the next section, we will
study the
{\it relative} holomorphic $L^2$ torsion for an arbitrary pair of flat
Hilbertian bundles over a complex manifold.

A {\it distance function} $r$ on a manifold $X$ is a map
$r:X\times X\to\Bbb{R}$ such that

\vskip 0.1in

(1) Its square $r^2(x,y)$ is smooth on $X\times X$.

(2) $r(x,x)=0$ and $r(x,y)>0$ if $x\neq y$.

(3) $\displaystyle{\partial^2 \over \partial x_i\partial x_j} r^2(x,y)
      \bigg|_{x=y} = g_{ij}(x)$

\vskip 0.1in

Condition (3) says essentially that $r(x,y)$ coincides with the geodesic
distance from $x$ to $y$, whenever $x$ and $y$ are close.  One can easily
construct such a function using local coordinates and a partition of unity.
Let
$$
   k(t,x,y) = c_1 t^{-n} e^{-c_2{r^2(x,y)\over t}},\quad t>0
$$
and $c_1,c_2$ are some positive constants.  Then one has the following basic
theorem about the fundamental solution of the heat equation,

\proclaim{5.2. Proposition}
The heat kernel $e^{-t\square^{\Cal E}_{p,q}}(x,y)$ is a smooth, symmetric
double form on $X$ and has the property
$$
   {\nabla''}_x e^{-t\square^{\Cal E}_{p,q}}(x,y) =
      {\nabla''}_y^* e^{-t\square^{\Cal E}_{p,q}}(x,y). \tag26
$$
It satisfies the bounds
$$
   \big| D e^{-t\square^{\Cal E}_{p,q}}(x,y) \big| \leq c_3 t^{-{1\over 2}}
      k(t,x,y)\tag27
$$
for $D={\nabla''}$ or ${\nabla''}^*$,\ $x,y$ close to each other
and $0<t\leq 1$.  Finally, there is a small time asymptotic expansion
$$
   e^{-t\square^{\Cal E}_{p,q}}(x,x) \sim \sum_{j=0}^\infty t^{-n+j}
      C_{j,p,q}(x)\tag28
$$
as $t\to 0$, where $C_{j,p,q}$ is a smooth double form on $X$, for all
$j$.
\endproclaim

\demo{Proof} The result is local, and in a local normal coordinate
neighborhood of a
point $x \in X$, where the bundle ${\Cal E}$ is also trivialized, one can
proceed
exactly as in \cite{RS1}, Proposition 5.3.
(cf. \cite{R}, \cite{BFKM})
\enddemo\hfill$\square$

\subheading{5.3} By Theorem 4.5, we see that the holomorphic $L^2$
torsion is not necessarily independent
of the choice of Hermitian metrics on $X$ and $\E\to X$. We will now study the
case when the flat Hilbertian $\A$-bundle $\E\to X$ with fiber $M$
is defined by a
{\it unitary representation} $\pi\rightarrow
{\Cal B}_{\Cal A}(M)$, that is, $M$ is a {\it unitary}
Hilbertian $(\A-\pi)$ bimodule. That is,
$$
\E\equiv (M\times\widetilde X)/\sim \to X
$$
where $(v,x)\sim(vg^{-1},gx)$ for all $g\in \pi$, $x\in \widetilde X$ and $v\in
M$.
The unitary representation defines a flat Hermitian
metric $h$ on $\E\to X$. We call such a bundle a {\it flat Hilbert bundle},
or sometimes a {\it unitary} flat Hilbertian bundle.
Then by definition (cf. 4.2), one has
$$
\rho_{\E}^p(g,h)\in\det(H^{p,*}(X,\M)).
$$

Let ${\Cal F}\to X$ be another flat Hilbert $\A$ bundle with fibre $N$,
such that
$\dim_\tau(M) = \dim_\tau(N)$. Let $\square^{\Cal E}_{p,q}(u)$ and
$\square^{\Cal F}_{p,q}(u)$
denote the Laplacians in the metric $g_u$, acting on
$\Omega^{p,q}_{(2)}(X,{\Cal E})$ and $\Omega^{p,q}_{(2)}(X,{\Cal F})$
respectively.   We first prove the following Proposition.

\proclaim{5.4. Proposition} Let $\E$ and $\F$ be a pair of flat Hilbert
bundles over $X$,
as above. Then there are positive constants $C_1, C$ such that
$$
\Big|\Tr_\tau(Z_u\exp(-t\square^\E_{p,q}(u)) ) -
\Tr_\tau(Z_u\exp(-t\square^\F_{p,q}(u)) )\Big| \leq C_1 e^{-{C\over t}}
$$
for all $0<t\le 1$.
\endproclaim

\demo{Proof}
Let $x\in X$ and assume that the ball $U_\delta = \{y\in X : r^2(x,y)<
\delta\}$ is simply connected, where $r$ is a distance function on $X$
which coincides with the geodesic distance on $U_\delta$.  Since the
Laplacian is a local operator, it follows that $\square_{p,q}^{\Cal E}$
acting on $\Omega_{(2)}^{p,q}(X,{\Cal E})$ over $U_\delta$ coincides
with $\square^{\Cal F}_{p,q}$ acting on $\Omega_{(2)}^{p,q}(X, {\Cal
F})$ over $U_\delta$.  By Duhamel's Principle and by
applying Green's theorem, one has for $x,y\in U_\delta$, one has
$$
\align
    e^{-t\square^{\Cal E}_{p,q}(u)}(x,y) -
      e^{-t\square^{\Cal F}_{p,q}(u)}(x,y)
   & = \int_0^t \int_{r^2(x,z)=\delta} \Bigl[
      e^{-(t-s)\square^{\Cal F}_{p,q}(u)}(z,y) \wedge*
      {\nabla''}^* e^{-s\square^{\Cal E}_{p,q}(u)}(x,z)  \\
   &  \phantom{= \int_0^t \int_{r^2(x,z)=\delta} \Bigl[}\
      -{\nabla''}^* e^{-s\square^{\Cal E}_{p,q}(u)}(x,z) \wedge*
      e^{-(t-s)\square^{\Cal F}_{p,q}(u)}(z,y) \\
   & \phantom{= \int_0^t \int_{r^2(x,z)=\delta} \Bigl[}\
      -e^{-s\square^{\Cal E}_{p,q}(u)}(x,z) \wedge*
      {\nabla''}^* e^{-(t-s)\square^{\Cal F}_{p,q}(u)}(z,y) \\
   &  \phantom{= \int_0^t \int_{r^2(x,z)=\delta} \Bigl[}\
      +{\nabla''}^* e^{-(t-s)\square^{\Cal F}_{p,q}(u)}(z,y) \wedge*
      e^{-s\square^{\Cal E}_{p,q}(u)}(x,z) \Bigr]
\endalign
$$
Using the basic estimate $(27)$ for heat kernels, one has
$$
   \Bigl| \operatorname{Tr}_\tau\bigl(Z_u e^{-t\square^{\Cal E}_{p,q}(u)}
      \bigr) - \operatorname{Tr}_\tau\bigl(Z_u e^{-t\square^{\Cal F}_{p,q}(u)}
      \bigr) \Bigr| \le c_1 t^{-{1\over 2}} e^{-{c_2\delta\over t}}\le C_1
e^{-{C\over t}}
$$
for all $0<t\le 1$.
\enddemo\hfill$\square$

\proclaim{5.5 Theorem} In the notation of 4.3, if
$\E$ and $\F$ are a pair of flat Hilbert bundles over $X$ which are of
$D$-class,
then the relative holomorphic $L^2$ torsion
$$
\rho_{\E,\F}^p = \rho_{\E}^p\otimes ({\rho_{\F}^p})^{-1} \in \det\Big(
H^{p,*}(X,
\E)\Big)\otimes \det\Big( H^{p,*}(X, \F)\Big)^{-1}
$$
is independent of the choice of Hermitian metric on $X$ which is needed to
define it.

\endproclaim

\demo{Proof}
Let $u\to g_u$ be a smooth family of Hermitian metrics on $X$ and
$\square^\E_{p,q}(u)$ and $\square^\F_{p,q}(u)$ denote the Laplacians
on $\E$ and $\F$ respectively, as before.

By Proposition 5.4, one has
$$
\Big|\Tr_\tau(Z_u\exp(-t\square^\E_{p,q}(u)) ) -
\Tr_\tau(Z_u\exp(-t\square^\F_{p,q}(u)) )\Big| \leq C_1 e^{-{C\over t}}
$$
as $t\to 0$. That is, $\Tr_\tau^s(Z_u\exp(-t\square^\E_{p,q}(u)) )$ and
$\Tr_\tau^s(Z_u\exp(-t\square^\F_{p,q}(u)) )$ have the same asymptotic expansion
as $t\to 0$. In particular, one has in the notation of Theorem 4.5,
$$
c_\E(u) = c_\F(u).
$$
Then the relative holomorphic $L^2$ torsion
$$
\gather
\rho_{\E,\F}^p\in\det H^{p,*}(X,\Cal{E})\otimes(\det H^{p,*}(X,\F))^{-1}\\
\rho_{\E,\F}^p(u)=\rho_\Cal{E}^p(u)\otimes(\rho_{\F}^p(u))^{-1}
\endgather
$$
satisfies
$$
\align
{\partial \over\partial u}\rho_{\E,\F}^p(u) =& \left({\partial \over\partial u}
   \rho_\Cal{E}^p(u)\right)\otimes\left(\rho_{\F}^p
   (u)\right)^{-1}-\rho_\Cal{E}^p(u)\otimes{\partial \over\partial u}
   \rho_{\F}^p(u)\cdot \rho_{\F}^p(u)^{-2}\\
=& (c_\Cal{E}(u)-c_{\F}(u)) \rho_{\E,\F}^p(u)\\
=& 0
\endalign
$$
using Theorem 4.5 and the discussion above.  This proves the theorem.
\enddemo\hfill$\square$

\heading{\bf \S 6. Determinant Line Bundles, Correspondences and Relative
Holomorphic $L^2$ Torsion} \endheading

In this section, we introduce the notion of determinant line bundles of
Hilbertian $\A$-bundles over compact manifolds. A main result in this
section is Theorem 6.8, which says that the holomorphic $L^2$ torsion
associated to a
flat Hilbertian bundle over a compact almost K\"ahler manifold, depends only on
the class of the Hermitian metric in the determinant line bundle of the flat
Hilbertian bundle. This enables us to show that a correspondence of determinant
line bundles is well defined on almost K\"ahler manifolds.
Finally, using such a correspondence of
determinant line bundles, we prove in Theorem 6.12 that the relative holomorphic
$L^2$ torsion is independent of the choices of almost K\"ahler metrics on the
complex manifold and Hermitian metrics on the pair of flat Hilbertian bundles
over the complex manifold.

\proclaim{6.1. Lemma}
The subgroup $SL(M) = \Det_\tau^{-1}(1)$ of $GL(M)$ is connected.
\endproclaim

\demo{Proof}
Let $U(M)$ denote the subgroup of all unitary elements in $GL(M)$. Recall
the standard retraction of $GL(M)$ onto $U(M)$, which is given by
$$
T_s : GL(M)\to GL(M)
$$
$$
A \to |A|^s \frac{A}{|A|}
$$
where $T_0 :GL(M)\to U(M)$ is onto and $T_1 = identity$. Clearly $U(M)
\subset SL(M)$
and the retraction $T_s$ above restricts to be a retraction of $SL(M)$ onto
$U(M)$. By the results of \cite{ASS}, it follows that $SL(M)$ is connected.
\enddemo\hfill$\square$

Let $\E\to X$ be a Hilbertian $\A$-bundle over $X$ and
$GL(\Cal{E})$ denote the space of complex $A$-linear automorphisms
of $\Cal{E}$ which induce the identity map on $X$, that is,
$GL(\Cal{E})$ is the gauge group of $\Cal{E}$.
The Fuglede-Kadison determinant, cf Theorem 1.32.
$$
\Det_\tau : GL(M)\rightarrow \Bbb{R}^+
$$
extends to a homomorphism
$$
\Det_\tau : GL(\Cal{E})\rightarrow C^\infty(X,\Bbb{R}^+)
$$
where $C^\infty(X,\Bbb{R}^+)$ denotes the space of smooth positive
functions on $X$.  This extension has all the properties listed in
theorem 1.32. Using the long exact sequence in homotopy and the Lemma
above, one has

\proclaim{6.2. Corollary} Let $\E\to X$ be a Hilbertian $\A$-bundle over $X$
(recall that $X$ is assumed to be connected). Then
the subgroup $SL(\E) = \Det_\tau^{-1}(1)$ of $GL(\E)$ is connected.
\endproclaim

\subheading{6.3. Determinant Line Bundles}
Let $\E\to X$ be a Hilbertian $\A$-bundle over $X$. Then we can define a
natural {\it determinant line bundle} of $\E$ as follows:

Let $\herm(\Cal{E})$ denote the space of all Hermitian metrics on
$\Cal{E}$.  Clearly $\herm(\Cal{E})$ is a convex set and $GL(\Cal{E})$
acts on $\herm(\Cal{E})$ by
$$
\gather
GL(\Cal{E})\times \herm(\Cal{E})\rightarrow \herm(\Cal{E})\\
(a,h)\rightarrow {\bar{a}}^t ha
\endgather
$$
That is, $(a.h)_x(v,w)=h_x(av,aw)$ for all $v, w\in\Cal{E}_x$.

The action of $GL(\Cal{E})$ on $\herm(\Cal{E})$ is transitive, that
is, one can identify $\herm(\Cal{E})$ with the quotient
$$
GL(\Cal{E})\big/ U(\Cal{E}, h_0)\qquad \text{where} \quad U(\Cal{E},h_0)
$$
is the subgroup of $GL(\Cal{E})$ which leaves $h_0\in
\herm(\Cal{E})$ invariant, that is, $U(\Cal{E},h_0)$ is the unitary
transformations with respect to $h_0$.

For a Hilbertian bundle $\Cal{E}$ over $X$, we define
$\det(\Cal{E})$ to be the real vector space generated by the symbols
$h$, one for each Hermitian metric on $\Cal{E}$, subject to the
following relations : for any pair $h_1$, $h_2$ of Hermitian metrics
on $\Cal{E}$, we write the following relation
$$
h_2=\sqrt{\Det_\tau(A)}^{\ -1} h_1
$$
where $A\in GL(\Cal{E})$ is positive, self-adjoint and satisfies
$$
h_2(v,w)=h_1(Av,w)
$$
for all $v, w\in \Cal{E}_x$.

Assume that we have three different Hermitian metrics $h_1$, $h_2$ and
$h_3$ on $\Cal{E}$.

Suppose that
$$
h_2(v,w)=h_1(Av,w) \text{ and } h_3(v,w)=h_2(Bv,w)
$$
for all $v,w\in \Cal{E}_x$ and $A,B\in GL(\Cal{E})$.  Then
$h_3(v,w)=h_1(ABv,w)$a and we have the following relations in
$\det(\Cal{E})$,
$$
\align
h_2 &= \sqrt{\Det_\tau(A)}^{\ -1} h_1\\
h_3 &= \sqrt{\Det_\tau(B)}^{\ -1} h_2\\
h_3 &= \sqrt{\Det_\tau(AB)}^{\ -1} h_1
\endalign
$$

The third relation follows from the first two, from which is follows
that $\det(\Cal{E})$ is a line bundle over $X$.

To summarize, $\det(\Cal{E})$ is a real line bundle over $X$, which
has nowhere zero sections $h$, where $h$ is any Hermitian metric on
$\Cal{E}$.  It has a canonical orientation, since the transition
functions $\sqrt{\Det_\tau(A)}^{\ -1}$ are always positive.

Non zero elements of $\det(\Cal{E})$ should be viewed as volume forms
on $\Cal{E}$.

For {\it flat} Hilbertian $\A$ bundles, the determinant line
bundle can be described in the following alternate way.

Then
$\Cal{E}=M\times_\rho\widetilde{X}$, where $\rho:\pi\rightarrow GL(M)$ is a
representation.  The associated {\it determinant line bundle} is
defined as
$$
\det{\Cal{E}}=\det(M)\times_{\Det_\tau(\rho)}\widetilde{X}.
$$
Here $\Det_\tau(\rho):\pi\rightarrow \Bbb{R}^+$ is a representation
which is defined as
$$
\Det_\tau(\rho) (\gamma) = \Det_\tau(\rho(\gamma))
$$
for $\gamma \in \pi$. Then $\det(\Cal{E})$ has the property that
$$
\det(\Cal{E})_x=\det(\Cal{E}_x) \qquad \forall x\in X.
$$
Clearly
$\det(\Cal{E})$ coincides with the construcion given in the beginning
of 6.3, and $\det(\Cal{E})$ is a {\it flat} real line bundle over $X$.

\subheading{6.4 Almost K\"ahler manifolds} A Hermitian manifold $(X,g)$
is said to be {\it almost K\"ahler}
if the K\"ahler 2-form $\omega$ is not necessarily closed, but instead
satisfies the weaker condition $\overline\partial\partial \omega = 0$.
Gauduchon (cf. \cite{Gau}) proved that every complex manifold of real
dimension less than or equal to 4, is almost K\"ahler.

Let $\nabla^B$ denote the holomorphic Hermitian connection on $TX$ with the
torsion tensor $T^B$ and curvature tensor $R^B$.  Define the smooth 3-form
$B$ by
$$
   B(U,V,W) = (T^B(U,V),W)
$$
for all $U,V,W\in TX$.  Let $\omega$ denote the K\"ahler 2-form on $X$.  Then
one has
$$
   B = i(\partial - \overline\partial) \omega.
$$
Since $X$ is almost K\"ahler, it follows that $B$ is closed and therefore the
following curvature identity holds
$$
   (R^B(U,V)W,Z) = (R^{-B}(Z,W)V,U)
$$
for all $U,V,W,Z\in TX$.  The Dolbeault operator $\sqrt{2}({\nabla''} +
\nabla^{\prime\prime*})$ is a Dirac type operator.  More precisely, Let
$\Lambda = (\det T^{\prime\prime0}X)^{1\over 2}$ and $\fs$ denote the
bundle of spinors on $X$, then as $\Bbb{Z}_2$ graded bundles on $X$, one
has
$$
   \Lambda^{p,*}T^*X\otimes {\Cal E} =
      \fs\otimes\Lambda\otimes\Lambda^{p,0} T^*X\otimes {\Cal E}.
$$
Let $\nabla^L$ denote the Levi-Civita connection on $X$ and ${\Cal D}^L$ the
Dirac
operator with respect to this connection.  Then using the connection
$\nabla^B$ on $\Lambda$ and $\Lambda^{p,0}T^*X$, the Dirac operator
${\Cal D}^L$ extends as an operator
$$
   {\Cal D}^L : \Gamma(X,\fs^+\otimes\Lambda\otimes\Lambda^{p,0}T^*X\otimes
      {\Cal E}) \to \Gamma(X,\fs^-\otimes\Lambda\otimes\Lambda^{p,0}
      T^*X\otimes{\Cal E})
$$
and one has the formula
$$
   \sqrt{2}(\nabla'' + \nabla^{\prime\prime*}) = {\Cal D}^L - {1\over 4}c(B)
      = {\Cal D}^L + {1\over2} \sum_{i=1}^n c(S(e_i)e_i)
$$
where $c(B)$ denotes Clifford multiplication by the 3-form $B$ and
$S=\nabla^B-\nabla^L$ is a 1-form on
$X$ with values in skew-Hermitian endomorphisms of $TX$.
We now work in a local
normal coordinate ball, where we trivialize the bundles using parallel
transport along geodesics.
Scale the metric on $X$ by $r^{-1}$ and let $I_r$ denote the operator
$2\square_{p,*} = (\sqrt{2}(\nabla'' + \nabla^{\prime\prime*}))^2$ in this
scaled metric.  In local normal coordinates, one has the following expression
for $I_r$ (cf. \cite{B})
$$
\align
   I_r &= -rg^{ij}\left(\partial_i + {1\over 4}\Gamma_{iab} c(e_a \wedge
e_b) + A_i
       + {1\over 2\sqrt{r}} c(S_{il\alpha}(e_l)e(f_\alpha)) + {1\over 4r}
      S_{i\beta\gamma} e(f_\beta\wedge f_\gamma)\right) \\
   &\quad \times\left(\partial_j + {1\over 4}\Gamma_{jab} c(e_a\wedge e_b) +
      A_j + {1\over 2\sqrt{r}} S_{jl\alpha} c(e_l)e(f_\alpha)
      \phantom{\bigl(} + {1\over 4r} S_{j\beta\gamma}
      e(f_\beta\wedge f_\gamma)\right)  \\
   & \quad + {1\over 4}rk - {1\over 2}rc(e_i\wedge e_j) L_{ij} -
      {1\over 2}e(f_\alpha\wedge f_\beta) L_{\alpha\beta}
      - \sqrt{r} c(e_i) e(f_\alpha\wedge L_{i\alpha})\\
   & \quad + rg^{ij}\Gamma_{ij}^k \left( \partial_k + {1\over 4}\Gamma_{kab}
      c(e_a\wedge e_b) + A_k +{1\over 2\sqrt{r}} S_{kl\alpha} c(e_l)
      e(f_\alpha) + {1\over 4r} S_{k\beta\gamma} e(f_\beta \wedge
f_\gamma)\right)
\endalign
$$
where $k$ denotes the scalar curvature of $X$.

Consider the heat equation on sections of
$\fs\otimes\Lambda\otimes\Lambda^{p,0}T^*X\otimes{\Cal E}$,
$$
\align
   (\partial_t + I_r) g(x,t) &= 0 \\
   g(x,0) &= g(x).
\endalign
$$
By parabolic theory, there is a fundamental solution $e^{-tI_r}(x,y)$ which
is smooth for $t>0$.  We will consider the case when $t=1,\ e^{-I_r}(x,y)$
and prove the existence of an asymptotic expansion on the diagonal, as
$r\to 0$.  A difficulty arises because of the singularities arising in the
coefficients of $I_r$, as $r\to 0$.

\proclaim{6.5. Proposition}
For some positive integer $p\ge n$, one has the following asymptotic expansion
as $r\to 0$,
$$
   e^{-I_r}(x,x) \sim r^{-p}\sum_{i=0}^\infty r^i E_i(x,x)
$$
where $E_i$ are endomorphisms of $\fs\otimes \Lambda\otimes
\Lambda^{p,0}T^*X\otimes{\Cal E}$.
\endproclaim

\demo{Proof}
Consider the operator
$$
\align
   J_r &= -rg^{ij}(\delta_i + {1\over 4}\Gamma_{iab} c(e_a\wedge e_b) +
      A_i) \times (\partial_j + {1\over 4}\Gamma_{jab}c(e_a\wedge e_b) + A_j) \\
   &\quad + rg^{ij}\Gamma_{ij}^k (\partial_k + {1\over 4} \Gamma_{kab}
      c(e_a\wedge e_b) + A_k) + {1\over 4}rk - {1\over 2}c(e_i\wedge e_j)
      L_{ij}.
\endalign
$$
Since $J_r$ has no singular terms as $r\to 0$, it has a well known
asymptotic expansion, as $r\to 0$ with $p=n$.

We can construct $\exp(-I_r)$ as a perturbation of $\exp(-J_r)$, using
Duhamel's principle.  More precisely,
$$
   \exp(-I_r) = \exp(-J_r) + \sum_{k=1}^\infty
      \underbrace{e^{-J_r}(J_r-I_r)e^{-J_r}\dots e^{-J_r}}_{k \text{terms}}
$$
Each coefficient in the difference $J_r-I_r$ contains at least one term
which is exterior multiplication by $f_\alpha$.  Therefore the infinite
series on the right hand side collapses to a finite number of terms.  The
proposition then follows from the asymptotic expansion for
$\exp(-J_r)(x,x)$.
\enddemo\hfill$\square$

\def\Tr{\operatorname{Tr}}
\def\ch{\operatorname{ch}}
Let $R^B$ denote the curvature of the holomorphic Hermitian connection and
$R^L$ denote the curvature of the Levi-Civita connection.  Let $\hat{A}$
denote $\hat{A}$-invariant polynomial and $\ch$ the Chern character
invariant polynomial.  Then
$$
   \hat{A}(R^{-B})\ch(\Tr(R^L) \ch (\Lambda^{p,0}R^L) \in \Lambda^*T^*X.
$$
The goal is to prove the following decoupling result in the adiabatic limit. It
resembles the local index theorem for almost K\"ahler manifolds by Bismut
\cite{Bi}
(he calls them non-K\"ahler manifolds). However, we use instead the
techniques of
the proofs in \cite{BGV}, \cite{Ge} and \cite{D} of the local index theorem
for families. In particular, we borrow a local conjugation trick due to
Donnelly \cite{D},
which is adjusted to our situation.

\proclaim{6.6. Theorem (Adiabatic decoupling)} Let $(X,g)$ be an almost
K\"ahler manifold.
In the notation above, one has the following
decoupling result in the adiabatic limit
$$
  \lim_{r\to 0} \Tr_\tau^s(Z_u e^{-I_r}(x,x)) = \Tr_\tau(Z_u)(x) \
      [ \hat{A}(R^{-B})\ch(\Tr(R^L)) \ch (\Lambda^{p,0}R^L)
      ]^{\max}_x\in\Lambda^{2n}T^*_xX $$
for all $x\in X$.
\endproclaim

\demo{Proof}
We first consider the corresponding problem on $\Bbb{R}^{2n}$, using the
exponential map.  Let $\bar{I}_r$ denote the operator on $\Bbb{R}^{2n}$, whose
expressions agrees with the local coordinate expression for $I_r$ near $p$,
where $p$ is identified with the origin in $\Bbb{R}^{2n}$.

Consider the heat equation on $\Bbb{R}^{2n}$,
$$
\gather
   (\partial_t + \bar{I}_r) g(x,t) = 0 \\
   g(x,0) = g(x).
\endgather
$$
Then one has
\enddemo\hfill$\square$

\proclaim{6.7. Proposition}
There is a unique fundamental solution $e^{-t\bar{I}_r}(x,y)$ which
satisfies the decay estimate
$$
   \bigl| e^{-t\bar{I}_r}(x,y)\bigr| \le c_1t^{-n} e^{-{c_2|x-y|^2\over t}}
$$
as $t\to 0$, with similar estimates for the derivatives in $x,y,t$.
\endproclaim

\demo{Proof} The proof is standard, as in 3.10.
\enddemo\hfill$\square$

By Duhamel's principle applied in a small enough normal coordinate neighborhood,
there is a positive constant $c$ such that
$$
   e^{-\bar{I}_r}(0,0) = e^{-I_r}(x,x) + O(e^{-c/r})\quad \text{as }r\to 0.
$$
Therefore $$\lim_{r\to0} \Tr_\tau^s(Z_u e^{-I_r}(x,x)) = \lim_{r\to0}
\Tr_\tau^s(Z_u e^{-\bar{I}_r}(0,0))\tag29$$
and it suffices to compute the right hand side of (29).  This is done using
Getzler's scaling idea \cite{Ge},  $x\to\epsilon x,\ t\to\epsilon^2t\
e_i\to \epsilon^{-1}e_i$.  Then Clifford
multiplication scales as $c_\epsilon(\cdot) = e(\cdot)  + \epsilon^2
i(\cdot)$, where
$e(\cdot)$ denotes exterior multiplication by the covector $\cdot$ and
$i(\cdot)$ denotes contraction by the dual vector.
$$
\align
   \bar{I}_\epsilon &= -rg^{ij}(\epsilon x) \Big(\partial_i +
      {\epsilon^{-1} \over 4}\Gamma_{iab}(\epsilon x)
      c_\epsilon(e_a\wedge e_b) + \epsilon A_i(\epsilon x) +
{\epsilon^{-1}\over 2\sqrt{r}}
      c_\epsilon(S_{il\alpha}(\epsilon x) e_i) e(f_\alpha)\\
   &\quad  + {\epsilon^{-1}\over 4r} S_{i\beta\gamma}{(\epsilon x)}
      c(f_\beta\wedge f_\beta)\Big) \times\Big(\partial_j +
      {\epsilon^{-1}\over 4} \Gamma_{jab}
   (\epsilon    x)  c_\epsilon(e_a\wedge e_b) + \epsilon A_j(\epsilon x)\\
   &\quad\phantom{\bigl(} + {\epsilon^{-1}\over 2\sqrt{r}} c_\epsilon
      (S_{jl\alpha} (\epsilon x) e_l) e(f_\alpha) + {\epsilon^{-1}\over 4r}
      S_{j\beta\gamma}(\epsilon x) e(f_\beta\wedge f_\gamma)\Big)
      + rg^{ij}(\epsilon x) \Gamma_{ij}^k(\epsilon x) \Big( \epsilon
      \partial_k\\
   &\quad  + {1\over4}\Gamma_{kab}(\epsilon x) c_\epsilon(e_a\wedge
      e_b) +\epsilon^2A_k(\epsilon x) + {1\over 2\sqrt{r}} S_{kl\alpha}
      (\epsilon x) c_\epsilon(e_i) e(f_\alpha)
      + {1\over 4r} S_{k\beta\gamma}(\epsilon x) e(f_\beta\wedge
      f_\gamma)\Big)  \\
   &\quad + {\epsilon^2\over 4} rk(\epsilon x)
   -{r\over 2} c_\epsilon(e_i\wedge e_j) L_{ij}(\epsilon x) -
      {1\over 2}f_\alpha\wedge f_\beta \wedge L_{\alpha\beta}(\epsilon x)
      -\sqrt{r} c_\epsilon(e_i) f_\alpha L_{i\alpha}(\epsilon x).
\endalign
$$
The asymptotic expansion in $r$ as in Propositions 6.5 and 6.7, for
$e^{-\bar{I}_r}(0,0)$ yields an asymptotic expansion in $\epsilon$ for
$e^{-\bar{I}_\epsilon}(0,0)$ and one
has
$$
   \lim_{r\to 0} \Tr_\tau^s\bigl(Z_u e^{-\bar{I}_r}(0,0)\bigr)
      = \lim_{\epsilon\to
      0}\Tr_\tau^s\bigl(Z_ue^{-\bar{I}_\epsilon}(0,0)\bigr)\tag30
$$
That is, if either limit exists, then both exist and are equal.

However, in the limit as $\epsilon\to 0$, there are singularities in the
coefficients of $S$ tensor in the expression for $\bar{I}_\epsilon$ and
one cannot immediately apply Getzler's theorem. Therefore one first makes
the following local conjugation trick, as in Donnelly \cite{D}.

Define the expression
$$
   h(x,\epsilon,r) = \exp\Bigl( {\epsilon^{-1}\over 2\sqrt{r}}
      S_{il\alpha}(0)x_ie_l\wedge f_\alpha + {\epsilon^{-1}\over 4r}
      S_{i\beta\gamma}(0) x_i f_\beta\wedge f_\gamma\Bigr).
$$
Note that $h(x,\epsilon,r)$ has polynomial growth in $x$, since its expression
contains exterior multiplication. We claim that if the operator
$\bar{I}_\epsilon$ is conjugated by $h$, then
the resulting operator is {\it not} singular as $\epsilon\to 0$.  More
precisely,
$$
\align
   J_\epsilon &= h\bar{I}_\epsilon h^{-1} \\
   &= rg^{ij}(\epsilon x) \Bigl(\partial_i + {\epsilon^{-1}\over4}\Gamma_{iab}
      (\epsilon x) e_a\wedge e_b + {\epsilon^{-1}\over 2\sqrt{r}}
      \bigl(S_{il\alpha} (\epsilon x) - S_{il\alpha}(0)\bigr) e_l\wedge
      f_\alpha+ {\epsilon^{-1}\over 4r} \bigl(S_{i\beta\gamma}(\epsilon x)\\
   &\quad   -
      S_{i\beta\gamma}(0)\bigr) f_\beta\wedge f_\gamma -{1\over 4r}
S_{il\alpha}(0) S_{kl\beta}(0)
      x_k f_\alpha\wedge f_\beta\Bigr) \times\Bigl( \partial_j +
{\epsilon^{-1}\over 4}\Gamma_{jab}(\epsilon x)
      e_a\wedge e_b\\
   &\quad + {\epsilon^{-1}\over 2\sqrt{r}}\bigl(S_{jl\alpha}
      (\epsilon x) - S_{i\beta\gamma}(0)\bigr) e_l\wedge f_\alpha +
{\epsilon^{-1}\over 4r}
      \bigl(S_{j\beta\gamma}(\epsilon x) - S_{j\beta\gamma}(0)\bigr)
      f_\beta\wedge f_\gamma -{1\over 4r} S_{jl\alpha}(0) S_{kl\beta}(0) x_k\\
   &\quad  f_\alpha\wedge f_\beta\Bigr)
   - {1\over 2} r e_i\wedge e_j L_{ij}(\epsilon x) -{1\over 2}
f_\alpha\wedge f_\beta L_{\alpha\beta}(\epsilon x) -
      \sqrt{r}e_i \wedge f_\alpha L_{i\alpha} (\epsilon x) +
R(x,\epsilon).\tag31
\endalign
$$
Here $R(x,\epsilon)$ denotes the terms which vanish as $\epsilon\to 0$, and
which therefore do not contribute to the limit. Clearly there are no singular
terms in $J_\epsilon$ as $\epsilon\to 0$.

A fundamental solution for the heat equation for $J_\epsilon$ can be
obtained by conjugating the one for $\bar{I}_\epsilon$, that is
$$
e^{-tJ_\epsilon}(x,y) = h(x, \epsilon, r) e^{-t\bar{I}_\epsilon}(x,y)
h^{-1}(y, \epsilon, r).
$$
The right hand side satisfies the heat equation
$(\partial_t+J_\epsilon)g(x,t) = 0,\ \ g(x,0) = \delta_x$.  Since $h(0)=1$,
one has $\forall\epsilon>0$,
$$
   \Tr_\tau^s\bigl(Z_ue^{-\bar{I}_\epsilon}(0,0)\bigr) =
      \Tr_\tau^s\bigl(Z_u e^{-J_\epsilon}(0,0)\bigr).\tag32
$$
Therefore it suffices to compute the limit as $\epsilon\to 0$ of the right
hand side of (31).

Using the following Taylor expansions in a normal coordinate neighborhood,
$$
\align
   \Gamma_{iab}(\epsilon x) &= -{1\over 2} R_{ijab}(0)\epsilon x_j +
      R(x,\epsilon^2) \\
   S_{il\alpha}(\epsilon x) &= S_{il\alpha}(0) + S_{il\alpha,j}(0)
      \epsilon x_j + R(x,\epsilon^2) \\
   S_{i\beta\gamma}(\epsilon x) &= S_{i\beta\gamma}(0) +
      S_{i\beta\gamma,j}(0) \epsilon x_j + R(x,\epsilon^2)
\endalign
$$
one sees that
$$
   J_0 = \lim_{\epsilon\to 0} J_\epsilon = -r\sum_i(\partial_i -
      {1\over4}B_{ij}x_j)^2 + r{\Cal L}
$$
where
$$
\align
   B_{ij} &= {1\over 2} R_{ijab}(0) e_a\wedge e_b -{2\over\sqrt{r}}
      S_{il\alpha,j}(0) e_l\wedge f_\alpha \\
   &\quad -{1\over r}\bigl(S_{i\beta\gamma,j}(0) - S_{il\beta}(0)
      S_{jl\gamma}(0)\bigr) f_\beta \wedge f_\gamma
\endalign
$$
and
$$
   {\Cal L} = {1\over2} L_{ij}(0) e_i\wedge e_j + {1\over\sqrt{r}}
      L_{i\alpha}(0) e_i\wedge f_\alpha + {1\over2r} L_{\alpha\beta}(0)
      f_\alpha\wedge f_\beta.
$$
Using Mehler's formula (cf. \cite{Ge}), one can obtain an explicit fundamental
solution $e^{-sJ_0}(x,y)$.  First decompose $B$ into its symmetric and skew
symmetric parts, that is
$B=C+D$ where $C={1\over2}(B+B^t)$ and $D={1\over2}(B-B^t)$, where $B,C,D$
are matrices of 2-forms.  Then

$$
   e^{-J_0}(x,0) = (4\pi r)^{-n/2}\hat{A}(r D) e^{{x^tCx\over 8}}
      \times\exp\Bigl(r{\Cal L} - {1\over 4r} x^t
      \Bigl({r D/2 \over \tanh(r D/2)}\Bigr) x\Bigr)
$$
Now $\lim_{\epsilon\to 0} e^{-J_\epsilon}(0,0) = e^{-J_0}(0,0)$.  Therefore
$$
   \lim_{\epsilon\to 0} \Tr_\tau^s(Z_u e^{-J_\epsilon}(0,0)) = \bigl(
      {2\over i}\bigr)^{n/2} (4\pi r)^{-n/2} \Tr_\tau(Z_u)(0) \bigl[
      \hat{A}(rD)\ch (r{\Cal L}) \bigr]^{\max}\tag33
$$
Here $D=R^{-B}(0)$ and ${\Cal L} = \Tr(R^L(0)) + \Lambda^{p,0} R^L(0)$.
Using (29), (30), (32) and (33), one completes the proof of Theorem 6.6.
$\square$
\vskip .2in

\proclaim{6.8. Theorem}
Let $\Cal{E}$ be a flat Hilbertian bundle of $D$-class,
over an almost K\"ahler manifold $(X,g)$ and let $h$, $h'$
be Hermitian metrics on $\Cal{E}$ such that $h=h'$ in
$\det(\Cal{E})$.  Then
$$
\rho_\Cal{E}^p(g,h)=\rho_\Cal{E}^p(g,h')\in\det(H^{p,*}(X,\Cal{E}))
$$
\endproclaim

\demo{Proof}
Since $h=h'$ in $\det(\Cal{E})$, there is a positive, self-adjoint
bundle map $A:\Cal{E}\rightarrow\Cal{E}$ satisfying
$$
\gather
h(A v,w)=h'(v,w) \qquad \forall v,w\in\Cal{E} \\
\text{and} \qquad \Det_\tau(A)=1.
\endgather
$$
By Corollary 6.2, there is a smooth 1-parameter family of
positive, self-adjoint
bundle maps $u\to A_u:\Cal{E}\rightarrow\Cal{E}$ joining
$A$ to the identity and satisfying
$$
\Det_\tau(A_u)=1. \tag34
$$
for all $u\in (-\epsilon,1+\epsilon)$. Here $A_0 = I$ and $A_1 = A$.
Let $u\rightarrow h_u$ be a smooth family of Hermitian metrics on
$\Cal{E}$ defined by
$$
h(A_u v,w)=h_u(v,w) \qquad \forall v,w\in\Cal{E}.
$$
Then $h_0=h$, $h_1=h'$ in $\E$ and $h=h_u$ in
$\det(\Cal{E})$ for all $u\in (-\epsilon,1+\epsilon)$
by (72).
Note that by differentiating (34), one has
$$
0 = {\partial \over\partial u}\Det_\tau(A_u) = \Tr_\tau(Z_u)
\tag35
$$
where $Z_u=A_u^{-1}\dot{A_u}$.

We wish to compute ${\partial \over\partial u}\rho_\Cal{E}^p (g,h_u)$.
By Theorem 4.5, one has
$$
{\partial \over\partial u}\rho_\Cal{E}^p(g,h_u)=c_\E^p(g,h_u)
   \rho_\Cal{E}^p(g,h_u)
$$
By Theorem 6.6 and (35), one sees that
$$
\lim_{t\rightarrow 0} \Tr_\tau^s(Z_u e^{-t\square(u)})=0.\tag36
$$
By the small time asymptotic expansion of the heat kernel,
one has
$$
\align
\lim_{t\rightarrow 0} \Tr_\tau^s(Z_u e^{-t\square(u)})&=
   \sum_{q=0}^n(-1)^q m_{n,p,q}(u)\\
&= c_\E^p (g,h_u)\tag37
\endalign
$$
Therefore by (36) and (37), one has $c_\E^p(g,h_u)=0$, that is,
$$
{\partial \over\partial u}\rho_\Cal{E}^p(g,h_u)=0.
$$
\enddemo\hfill$\square$

\subheading{6.9. Remarks} Theorem 6.8 says that on an almost
K\"ahler manifold $(X,g)$, the
holomorphic $L^2$ torsion $\rho_\Cal{E}^p(g,h)$ depends only on the
equivalence class of the Hermitian metric $h$ in $\det(\Cal{E})$.
We however do not believe that the almost K\"ahler
hypothesis in Theorem 6.8 is necessary. However, we use the techniques
of the proof of the local index theorem, and the
situation to date is that the local index theorem for the operator
$\overline\partial +\overline\partial^*$
has not yet been established for a general Hermitian manifold.

\subheading{6.10} Let $\Cal{E}$ and $\F$ be two flat Hilbertian bundles
of $D$-class over over an almost K\"ahler manifold $(X,g)$
and $\varphi:\det(\Cal{E})\rightarrow \det(\F)$ be an isomorphism of the
determinant line bundles.  Then using the theorem above, we will
construct a canonical isomorphism between determinant lines
$$
\gather
\widehat{\varphi}^p:\det H^{p,*}(X,\Cal{E})\rightarrow \det H^{p,*}(X,\F)\\
\widehat{\varphi}^p(\lambda \rho_\Cal{E}^p(g,h))=\lambda \rho_{\F}^p(g,h'),
\quad \lambda\in\Bbb{R}
\endgather
$$
where $h$ and $h'$ are Hermitian metrics on $\Cal{E}$ and $\F$
respectively, such that $\varphi(h)=h'$ in $\det(\F)$.
Then $\widehat{\varphi}$ is called a {\it correspondence} between determinant
line bundles. It is well defined by Theorem 6.8 and Remarks 6.9.
We next state some obvious properties of correspondences.

\proclaim{6.11. Proposition}
Let $\Cal{E}$  be a flat Hilbertian bundle of $D$-class over
over an almost K\"ahler manifold $(X,g)$ and
$\varphi:\det(\Cal{E})\rightarrow \det(\E)$ be the identity
map. Then $$\widehat{\varphi}^p= identity$$.

Let $\Cal{E}, \F$ and $\G$ be flat Hilbertian bundles of $D$-class over
over an almost K\"ahler manifold $(X,g)$ and
$\varphi:\det(\Cal{E})\rightarrow \det(\F)$,
$\psi:\det(\Cal{F})\rightarrow \det(\G)$
be isomorphisms of the determinant line bundles. Then the composition
satisfies
$$
\widehat{\varphi {\small\text o}\psi}^p = \widehat{\varphi}^p {\small\text o}
\widehat{\psi}^p.$$
\endproclaim

We next prove one of the main results in the paper.

\proclaim{6.12. Theorem}
Let $\Cal{E}$ and $\F$ be two flat Hilbertian bundles of $D$-class over
over an almost K\"ahler manifold $(X,g)$ and
$\varphi:\det(\Cal{E})\rightarrow \det(\F)$ be an isomorphism of the
corresponding determinant line bundles.  Consider smooth 1-parameter
families of almost K\"ahler metrics $g_u$ on $X$ and Hermitian metrics
$h_{1,u}$ on $\Cal{E}$, where $u$ varies in an internal $(-\epsilon,\epsilon)$.
Choose a smooth family of Hermitian metrics $h_{2,u}$ on $\F$ in such
a way that $\varphi(h_{1,u})=h_{2,u}$ in $\det(\F)$.  Then the relative
holomorphic torsion
$$
\rho_\varphi^p(u)=\rho_\Cal{E}^p(g_u,h_{1,u})\otimes\rho_{\F}^p(g_u,h_{2,u})
^{-1}
\in\det H^{p,*}(X,\Cal{E})\otimes\det H^{p,*}(X,\F)^{-1}
$$
is a smooth function of $u$ and satisfies ${\partial \over\partial u}
\rho_\varphi(u)=0$.  That is, the relative holomorphic
$L^2$ torsion $\rho_\varphi^p$
is independant of the choices of metrics on $X$, $\Cal{E}$
and $\F$ which are needed to define it.
\endproclaim

\demo{Proof}From the data in the theorem,
 one can define a correspondence as in 6.10,
$$
\hat{\varphi}^p:\det(H^{p,*}(X,\Cal{E}))\rightarrow \det(H^{p,*}(X,\F))
$$
which is an isomorphism of determinant lines.  It is defined as
$$
\hat{\varphi}^p(\lambda \rho_\Cal{E}^p(g_u,h_{1,u}))=\lambda\rho_{\F}^p
(g_u,h_{2,u})
\tag{38}
$$
for $\lambda\in\Bbb{R}$ and $u\in(-\epsilon,\epsilon)$.  Therefore
using theorem 4.5 and (38) above, one has
$$
\align
{\partial \over\partial u}\hat{\varphi}^p(\rho_\Cal{E}^p(g_u,h_{1,u}))&=
   \hat{\varphi}^p({\partial \over\partial u}\rho_\Cal{E}^p(g_u,h_{1,u}))\\
&=c_\Cal{E}(g_u,h_{1,u})\hat{\varphi}^p(\rho_\Cal{E}^p(g_u,h_{1,u}))\\
&=c_\Cal{E}(g_u,h_{1,u})\rho_{\F}^p(g_u,h_{2,u}).
\tag{39}
\endalign
$$
But by differentiating equation (38) above, one has
$$
\align
{\partial \over\partial u}\hat{\varphi}(\rho_\Cal{E}^p(g_u,h_{1,u}))
   &={\partial \over\partial u}\rho_{\F}^p(g_u,h_{2,u})\\
&=c_{\F}(g_u,h_{2,u})\rho_{\F}^p(g_u,h_{2,u})
\tag{40}
\endalign
$$
Equating (39) and (40), one has
$$
c_\Cal{E}(g_u,h_{1,u})=c_{\F}(g_u,h_{2,u})\tag41
$$
Then the relative holomorphic $L^2$ torsion
$$
\gather
\rho_\varphi^p\in\det H^{p,*}(X,\Cal{E})\otimes(\det H^{p,*}(X,\F))^{-1}\\
\rho_\varphi^p(u)=\rho_\Cal{E}^p(g_u,h_{1,u})\otimes(\rho_{\F}^p(g_u,h_{2,u}
))^{-1}
\endgather
$$
satisfies
$$
\align
{\partial \over\partial u}\rho_\varphi^p(u) =& \left({\partial \over\partial u}
   \rho_\Cal{E}^p(g_u,h_{1,u})\right)\otimes\left(\rho_{\F}^p
   (g_u,h_{2,u})\right)^{-1}\\
&-\rho_\Cal{E}^p(g_u,h_{1,u})\otimes{\partial \over\partial u}
   \rho_{\F}^p(g_u,h_{2,u})\cdot \rho_{\F}^p(g_u,h_{2,u})^{-2}\\
=& c_\Cal{E}(g_u,h_{1,u})\rho_\varphi^p(u)
-c_{\F}(g_u,h_{2,u}) \rho_\varphi^p(u)\\
=& 0
\endalign
$$
using Theorem 4.5 and (41) above.  This proves the theorem.
\enddemo\hfill$\square$

\heading{\bf \S 7. Calculations} \endheading

In this section, we  calculate the holomorphic $L^2$ torsion for K\"ahler
locally symmetric
spaces. We will restrict ourselves to the special case of the Hilbert
$({\Cal U}(\Gamma)-\Gamma)$-bimodule $\ell^2(\Gamma)$,
where $\Gamma$ is a countable discrete group. Let $\E\to X$ denote the
associated flat Hilbert ${\Cal U}(\Gamma)$-bundle over the compact complex
manifold
$X$. Then it is well known that
the Hilbert ${\Cal U}(\Gamma)$-complexes
$\left(\Omega_{(2)}^{\bullet,\bullet}({X},\E),\nabla''\right)$ and
$\left(\Omega_{(2)}^{\bullet,\bullet}(\widetilde{X}),
\bar{\partial}\right)$ are canonically
isomorphic, where $\Gamma\to{\widetilde X}\to X$ denotes the universal
covering space
of $X$ with structure group $\Gamma$. We will denote the
$\bar{\partial}$-Laplacian
acting on $\Omega_{(2)}^{p,q}(\widetilde{X})$ by $\square_{{p,q}}$.

Firstly, we will discuss the $D$-class condition in this case. Let $X$ be a
K\"ahler hyperbolic manifold. Recall that this means that $X$ is a K\"ahler
manifold with K\"ahler form $\omega$, which has the property that $p^*(\omega)
= d\eta$, where $\Gamma\to{\widetilde X}\to X$ denotes the universal cover
of $X$
and $\eta$ is a bounded 1-form on ${\widetilde X}$. Any Riemannian manifold of
negative sectional curvature,
which also supports a K\"ahler metric, is a K\"ahler hyperbolic manifold.
Note that
the K\"ahler metric is not assumed to be compatible with the Riemannian
metric of
negative sectional curvature. Then Gromov \cite{G} proved that on the universal
cover of a K\"ahler hyperbolic manifold, the Laplacian $\square_{{p,q}}$
has a spectral gap at zero
on all $L^2$ differential forms. Therefore it follows that the associated
flat bundle
$\E\to X$ is of $D$-class. By a vanishing theorem of Gromov \cite{G}
for the $L^2$ Dolbeault cohomology of the universal cover, one has
$$
H_{(2)}^{p,q}(\widetilde X) = 0
$$
unless $p+q = n$, where $n$ denotes the complex dimension of $X$.

In particular, let
$G$ be a connected semisimple Lie group, and $K$ be a maximal compact
subgroup such that
$G/K$ carries an invariant complex structure, and let $\Gamma$ be a
torsion-free uniform lattice
in $G$. Then it is known that $\Gamma\backslash G/K$ is a K\"ahler
hyperbolic manifold
(cf. \cite{BW})
and therefore the canonical flat Hilbert bundle $\E\to X$ is of $D$-class.
In this K\"ahler
metric, the Laplacian $\square_{{p,q}}$ is $G$-invariant,
so it follows that the theta function
$$
\theta_{p,q}(t) = C_{p,q}(t) vol(\Gamma\backslash G/K)
$$
is proportional to the volume of $\Gamma\backslash G/K$. Here $C_{p,q}(t)$
depends only
on $t$ and on $G$ and $K$, but {\it not} on $\Gamma$. It follows that the
zeta function
$\zeta_{p,q}(s, \lambda, \E)$ is also proportional to the volume of
$\Gamma\backslash G/K$. Therefore
the holomorphic $L^2$ torsion is given by
$$
\rho_\Cal{E}^p  = e^{C_p vol(\Gamma\backslash G/K)} \rho^{\prime p} \in
\det \left( H_{(2)}^{p,n-p}(G/K)\right)^{(-1)^{n-p}}
$$
where we have used the vanishing theorem of Gromov. Here $C_{p}$ is a
constant that depends only
on $G$ and $K$, but {\it not} on $\Gamma$. Using representation theory, as
for instance in \cite{M},
\cite{L} and \cite{Fr}, it is possible to determine $C_p$ explicitly. This
will be done elsewhere.
Using the proportionality principle again, one sees that the Euler
characteristic of
$\Gamma\backslash G/K$ is proportional to its volume.
By  a theorem of Gromov \cite{G}, the Euler characteristic of
$\Gamma\backslash G/K$ is non-zero. Therefore we can also express the
holomorphic $L^2$ torsion as
$$
\rho_\Cal{E}^p  = e^{C'_p \chi(\Gamma\backslash G/K)} \rho^{\prime p} \in
\det \left( H_{(2)}^{p,n-p}(G/K)\right)^{(-1)^{n-p}}
$$
where $\chi(\Gamma\backslash G/K)$ denotes the Euler characteristic of
$\Gamma\backslash G/K$,
and $C'_{p}$ is a constant that depends only
on $G$ and $K$, but {\it not} on $\Gamma$. This discussion is summarized in
the following
proposition.

\proclaim{7.1. Proposition} In the notation above, the holomorphic $L^2$
torsion of the
semisimple locally symmetric space $\Gamma\backslash G/K$, which is assumed to
carry an invariant complex structure, is given by
$$
\rho_\Cal{E}^p  = e^{C_p vol(\Gamma\backslash G/K)} \rho^{\prime p} \in
\det \left( H_{(2)}^{p,n-p}(G/K)\right)^{(-1)^{n-p}}
$$
Here $C_{p}$ is a constant that depends only
on $G$ and $K$, but {\it not} on $\Gamma$. Equivalently, the holomorphic
$L^2$ torsion of
$\Gamma\backslash G/K$ is given as
$$
\rho_\Cal{E}^p = e^{C'_p \chi(\Gamma\backslash G/K)} \rho^{\prime p} \in
\det \left( H_{(2)}^{p,n-p}(G/K)\right)^{(-1)^{n-p}}
$$
where $\chi(\Gamma\backslash G/K)$ denotes the Euler characteristic of
$\Gamma\backslash G/K$,
and $C'_{p}$ is a constant that depends only
on $G$ and $K$, but {\it not} on $\Gamma$.
\endproclaim

We will now compute the
holomorphic $L^2$ torsion for a Riemann surface, which is a special case of
the theorem
above, and we will show that the constants $C_p$ and $C'_p$ are not zero.

Let $X$ be a closed Riemann surface of genus $g$, which is greater than 1,
which can be realised as
a compact quotient complex hyperbolic space $\Bbb H$ of complex dimension 1,
by the torsion-free discrete group $\Gamma$.
Recall that
$$
{\square}_{0,1} = \frac{1}{2} \Delta_1
$$
acting on the subspace $\Omega_{(2)}^{0,1}(\Bbb H)$.
Also,
$$
{\square}_{1,0} = \frac{1}{2} \Delta_1
$$
acting on the subspace $\Omega_{(2)}^{1,0}(\Bbb H)$.
The $\star$ operator intertwines the operators ${\square}_{0,1}$ and
${\square}_{1,0}$, showing in particular that they are isospectral.
So using
$$
\Omega_{(2)}^{1}(\Bbb H) =
\Omega_{(2)}^{1,0}(\Bbb H) \oplus \Omega_{(2)}^{0,1}(\Bbb H),
$$
we see that in order to calculate the von Neumann determinant of the operator
${\square}_{1,0}$, it suffices to
first scale the metric $g\to 2g$ and then calculate the square root
von Neumann determinant of the Laplacian $\Delta_1$. However, this is easily
seen to be equal to the von Neumann determinant of the Laplacian $\Delta_0$
acting on $L^2$ functions on the hyperbolic disk. Recall that the von Neumann
determinant of the operator $A$ is by definition $e^{-\zeta'_A(0)}$, where
$\zeta'_A(s)$ denotes the zeta function of the operator $A$.

Using the work of Randol \cite{R}, one obtains the following expression for the
meromorphic continuation of the zeta function of $\Delta_0$ to the
half-plane $\Re(s)<1$
$$
\zeta_0(s,0,\Cal{E}) = (g-1) \frac{\pi}{(s-1)} \int_0^\infty \left( {1\over
4} +r^2\right)^{1-s}
sech^2(\pi r) dr.
$$
It follows that
$$
\align
\zeta'_0(0,0,\Cal{E}) & = \lim_{s\to 0} \left(\zeta_0(s,0,\Cal{E}) -
\zeta_0(0,0,\Cal{E}) \right)\Gamma(s)\\
& =  (g-1){\pi} \int_0^\infty \left( {1\over 4} +r^2\right)
sech^2(\pi r) \left( -1 + \log\left( {1\over 4} +r^2\right)\right)dr.
\endalign
$$
A numerical approximation for the last integral shows that
$\zeta'(0,0,\Cal{E}) \sim -0.677 (g-1) $. We can summarize
the discussion in the following proposition.

\proclaim{7.2. Proposition} In the notation above, the holomorphic $L^2$
torsion of a
compact Riemann surface $X = \Gamma \backslash \Bbb H$ of genus $g$, is given by
$$
\rho_\Cal{E}^0 = e^{C (g-1)} \rho^{\prime 0} \in
\det \left( H_{(2)}^{0,1}(\Bbb H)\right)^{(-1)}\tag42
$$
Here $C = {\pi\over 2} \int_0^\infty \left( {1\over 4} +r^2\right)
sech^2(\pi r) \left( -1 + \log\left( {1\over 4} +r^2\right)\right)dr$
is a constant that depends only on $\Bbb H$, but {\it not} on $\Gamma$.
$C$ is approximately $-0.338$, and in particular, it is not equal to zero.
Also,
$$
\rho_\Cal{E}^1 = e^{- C (g-1)} \rho^{\prime 1} \in
\det \left( H_{(2)}^{1,0}(\Bbb H)\right)^{(-1)}
$$
where the constant $C$ is as in (42).

\endproclaim

\subheading{Acknowledgement} We thank John Phillips for his cleaner proof
of Lemma 6.9.

\enddocument


--============_-1366095342==_============--

--============_-1364082820==_============--

From acarey Mon Nov 18 09:38:54 1996
Received: (from acarey@localhost) by spam.maths.adelaide.edu.au (8.8.0/8.8.0/UofA-1.5) id JAA09732 for vmathai@spam; Mon, 18 Nov 1996 09:38:54 +1030 (CST)
From: Alan Carey <acarey@maths.adelaide.edu.au>
Message-Id: <199611172308.JAA09732@spam.maths.adelaide.edu.au>
Subject: no subject (file transmission)
To: vmathai@spam.maths.adelaide.edu.au
Date: Mon, 18 Nov 1996 09:38:53 +1030 (CST)
X-Mailer: ELM [version 2.4 PL22]
MIME-Version: 1.0
Content-Type: text/plain; charset=US-ASCII
Content-Transfer-Encoding: 7bit
Status: RO


%
%
%
\input amstex
\documentstyle{amsppt}
\NoBlackBoxes
\TagsOnRight
\voffset= 1cm
\hoffset=.5cm
\hsize=6.1in
\vsize=8.2in
\strut
\vskip5truemm
\font\small=cmr8

\def\det{\operatorname{det}}
\def\herm{\operatorname{Herm}}
\define\C{{\Bbb C}}
\define\R{{\Bbb R}}

\define\M{{\Cal E}}

\define\G{{\Cal G}}

\define\fs{{\Cal S}}

\define\ch{\operatorname{ch}}
\define\Tr{\operatorname{Tr}}

\define\Det{\operatorname{Det}}

\define\End{\operatorname{End}}

\define\GL{\operatorname{GL}}

\define\A{{\Cal A}}
\redefine\B{{\Cal B}}

\define\F{{\Cal F}}

\define\im{\operatorname{im}}

\define\cl{\operatorname{cl}}
\redefine\H{\Cal H}

\define\E{\Cal E}

\def\<{\langle}
\def\>{\rangle}

\define\pd#1#2{\dfrac{\partial#1}{\partial#2}}
\documentstyle{amsppt}

\topmatter

\title{Correspondences, von Neumann algebras}\\ {and holomorphic L$^2$
torsion} \endtitle

\author ALAN L. CAREY, MICHAEL S. FARBER and VARGHESE MATHAI\endauthor

\address Department of Pure Mathematics, University of Adelaide, Adelaide 5005,
Australia. \endaddress

\email acarey$\@$maths.adelaide.au.edu\endemail

\address School of Mathematical Sciences, Tel-Aviv University, Tel-Aviv 69978,
Israel. \endaddress

\email farber$\@$math.tau.ac.il\endemail

\address Department of Pure Mathematics, University of Adelaide, Adelaide 5005,
Australia. \endaddress

\email vmathai$\@$maths.adelaide.au.edu\endemail



\abstract
Given a holomorphic Hilbertian bundle on a compact complex manifold,
we introduce the notion of holomorphic $L^2$ torsion, which lies in the
determinant line of the twisted $L^2$ Dolbeault cohomology and represents a
volume element there. Here we utilise the theory of determinant lines of
Hilbertian modules over finite von Neumann algebras as developed in \cite{CFM}.
This specialises to the Ray-Singer-Quillen holomorphic torsion in the
finite dimensional case.
We compute a metric variation formula for the holomorphic $L^2$ torsion,
which shows that it is {\it not} in general independent of the choice of
Hermitian metrics on the complex manifold and on the holomorphic Hilbertian
bundle, which are needed to define it. We therefore initiate the theory of
correspondences of determinant lines, that enables us to define a relative
holomorphic $L^2$ torsion for a pair of flat Hilbertian bundles,
which we prove is independent of the choice of
Hermitian metrics on the complex manifold and on the flat Hilbertian
bundles.
\endabstract




\endtopmatter

\document
\heading{\S 0. Introduction}\endheading

Ray and Singer (cf. \cite{RS}) introduced the notion of holomorphic torsion of a
holomorphic bundle over a compact complex manifold. In  \cite{Q}, Quillen
viewed the
holomorphic torsion as an element in the real determinant line of the twisted
Dolbeault cohomology, or equivalently, as a metric in the dual of the
determinant line of the twisted Dolbeault cohomology. Since then there have
been many generalisations in the finite dimensional case, particularly by
Bismut, Freed, Gillet and Soule, \cite{BF, BGS}.

In this paper, we investigate generalisations
of aspects of this previous work to the case of infinite dimensional
representations of
the fundamental group.
Our approach is to introduce the concepts of holomorphic Hilbertian bundles
and of connections compatible with the holomorphic structure.
These bundles have fibres which are von Neumann algebra modules.
 We are able to define the
{\it determinant line bundle} of a holomorphic Hilbertian
bundle over a compact complex manifold, generalising the
construction of the determinant line of a finitely
generated Hilbertian module
that was developed in our earlier paper \cite{CFM}.
A nonzero element of the
determinant line bundle can be naturally
viewed as a volume form on the Hilbertian bundle.
This enables us to make sense of the notions of volume form and determinant
line bundle in this infinite dimensional and non-commutative situation.
Given an isomorphism of the determinant line bundles of holomorphic
Hilbertian bundles, we introduce the concept of a {\it correspondence}
between the  determinant lines of the twisted $L^2$ Dolbeault cohomologies.
This was previously studied in the finite dimensional situation in \cite{F}.

Restricting our attention to the class of manifolds studied
in [BFKM] (the so-called determinant or $D$-class
examples)
we then define the holomorphic $L^2$ torsion of a holomorphic Hilbertian
bundle;
it reduces to the classical constructions in the finite dimensional situation.
This new torsion invariant lives in the determinant line
of the twisted $L^2$ Dolbeault cohomology. Some key results in our
paper are a metric variation formula for the holomorphic $L^2$ torsion,
and the definition of a correspondence
between the determinant lines of the twisted $L^2$ Dolbeault cohomologies
for a pair of flat holomorphic Hilbertian bundles, and finally the definition
of a metric independent relative holomorphic $L^2$ torsion associated
to a correspondence between determinant line bundles of flat 
Hilbertian bundles.
To prove that a correspondence between determinant
line bundles of flat Hilbertian bundles is well defined, we need to prove
a generalised local index theorem for almost K\"ahler manifolds, and as a
consequence,
we give an alternate proof of Bismut's local index theorem for almost
K\"ahler manifolds
\cite{B}, where we use instead the methods of Donnelly \cite{D} and
Getzler \cite{Ge}.

The paper  is organized as follows. In the first section, we recall
some preliminary material
on Hilbertian modules over finite von Neumann algebras,
the canonical trace on the commutant of a finitely generated Hilbertian
module, the Fuglede-Kadison determinant on Hilbertian modules
and the construction of determinant lines for finitely generated
Hilbertian modules. Details of the material in this section can be
found in \cite{CFM}. In section \S 2, we define Hilbertian bundles
and connections on these. The definition of a connection is tricky in the
infinite dimensional context, and we use some fundamental theorems in von
Neumann
algebras to make sense of our definition. Then we define holomorphic
Hilbertian bundles
and connections compatible with the holomorphic structure
as well as Cauchy-Riemann operators on these.
In section \S 3, we study the properties of the zeta function associated to
holomorphic
Hilbertian bundles of $D$-class. In section \S 4, we define
the holomorphic $L^2$ torsion as an element in the determinant line
of reduced $L^2$ Dolbeaut cohomology. Here we also prove metric variation
formulae and we deduce that holomorphic $L^2$ torsion {\it does} depend
on the choices of Hermitian metrics on the compact complex manifold and
on the holomorphic Hilbertian bundle. However, in sections \S 5 and \S 6,
we give situations when a relative version of the holomorphic $L^2$ torsion
is indeed independent of the choice of metric.
In section \S 5, we are able to deduce
the following theorem (Theorem 5.5 in the text) from the variation formula: let
${\E}$ and ${\F}$ be two flat Hilbert bundles of $D$-class over
a compact Hermitian manifold $X$. Then one can define a relative holomorphic
$L^2$ torsion
$$
\rho_{\E,\F}^p\in \det(H^{p,\ast}(X,\E))\otimes \det(H^{p,\ast}(X,\F))^{-1}
$$
which is independent of the choice of Hermitian metric on $X$.
In section \S 6, we
define the notion of the determinant line bundle of a Hilbertian bundle and
also of
correspondences between determinant lines. The proof that a correspondence
is well defined, uses techniques of Bismut \cite{B}, Donnelly \cite{D} and
Getzler \cite{Ge} in their proof of the local index theorem in
different situations. Using the notion of a
correspondence of determinant line bundles,
we prove one of the main theorems in our paper (Theorem 6.12 in the text),
which can be briefly  stated as follows: let
${\E}$ and ${\F}$ be two flat Hilbertian bundles of $D$-class over
a compact almost K\"ahler manifold $X$
and $\varphi:\det(\E)\to\det(\F)$ be an isomorphism of the corresponding
determinant line bundles. Then one can define a relative holomorphic
$L^2$ torsion
$$
\rho_\varphi^p\in \det(H^{p,\ast}(X,\E))\otimes \det(H^{p,\ast}(X,\F))^{-1}.
$$
Using the correspondence defined by the isomorphism $\varphi$, we show
that the relative holomorphic
$L^2$ torsion $\rho_\varphi^p$ is independent of the choices of
Hermitian metrics on $\E$ and $\F$
and the choice of almost K\"ahler metric on $X$ which are needed to define it.
Recall that an almost K\"ahler manifold is a Hermitian manifold whose
"K\"ahler" 2-form $\omega$ is not necessarily closed, but satisfies
the weaker condition $\overline\partial\partial\omega = 0$. A result
of Gauduchon (cf. \cite{Gau}) asserts that every compact complex surface is
almost K\"ahler,
whereas there are many examples of complex surfaces which are not K\"ahler.
In  section \S 7, we give some examples of calculation of the
holomorphic $L^2$ torsion for locally symmetric spaces and Riemann surfaces.

\heading{\bf \S 1. Preliminaries}\endheading

This section contains some preliminary material from \cite{CFM}.

\subheading{1.0. Hilbertian modules over von Neumann algebras}

 Throughout the paper ${\A}$ will denote
 a finite  von Neumann algebra
with a fixed finite,
normal, and faithful trace $\tau:{\Cal A}\rightarrow \C$.  The involution
in ${\A}$ will be denoted $*$ while $\ell^{2}({\Cal A})$ denotes the completion
of  ${\A}$ in the norm derived from the inner product
$\tau(a^*b),\  a,b\in \A$.
A {\it Hilbert module} over ${\A}$
is a Hilbert space $M$
together with a continuous left ${\Cal A}$-module structure such that there
exists an isometric ${\Cal A}$-linear embedding of $M$ into
$\ell^{2}({\Cal A})\otimes H$, for some Hilbert space $H$.
(Note that this embedding is not
part of the structure.)  A Hilbert module $M$ is {\it finitely generated} if it
admits an imbedding as above with
finite dimensional $H$.
To introduce the notion of determinant line requires us
to  forget
the scalar product on $H$ but keep the topology and the ${\Cal A}$-action.

\subheading{1.1. Definition} A {\it  Hilbertian module}
is a topological vector space $M$
with continuous left ${\Cal A}$-action such that there exists a scalar product
$\langle\;,\;\rangle$ on $M$ which generates the topology of $M$ and such
that $M$ together
with $\langle\;,\;\rangle$ and with the ${\Cal A}$-action is a Hilbert
module. Any scalar product $\langle\;,\;\rangle$ on
$M$ with the above properties will be called {\it admissible}.

\subheading{1.2. Remarks and further definitions} The
choice of any other
 admissible scalar product $\langle\;,\;\rangle_1$ gives an isomorphic
 Hilbert module.  In fact there exists an operator
$
A:M\rightarrow M
$
such that
$$
\langle v,w\rangle_1=\langle Av,w\rangle\tag1
$$
for any $v,w\in M$.  The operator $A$ must be a self-adjoint, positive
linear homeomorphism (since the scalar products $\<\ ,\ \>$ and
$\<\ ,\ \>_1$ define the same topology),
which commutes with the $\A$-action. A
{\it finitely generated Hilbertian
module} is one for which the corresponding Hilbert module
 is finitely
generated. Finally, a
 {\it morphism} of Hilbertian modules is a continuous linear map
$f:M\to N$, commuting with the $\A$-action.
Note that the kernel of any morphism $f$ is again a Hilbertian module as is
the closure of the image $\cl(\im(f))$.

\subheading{1.3. The canonical trace on the commutant}

Any choice of an admissible scalar product $\langle\;,\;\rangle$ on $M$,
defines obviously a $*$-operator on $\B$ (by assigning to an operator
its adjoint) and turns $\B$ into a von Neumann algebra.
If we choose another admissible scalar product $\<\ ,\ \>_1$
on $M$ then the new involution will be given by
$$
f\mapsto \ A^{-1}f^*A\qquad \text{for}\quad f\in \B,\tag2
$$
where $A\in\B$ satisfies
$\<v,w\>_1\ =\ \<Av,w\>$ for $v,w\in M$.
The trace on the commutant may now be defined as in \cite{Di}
 and here will be denoted
$\Tr_{\tau}$. It
is finite, normal, and faithful. If $M$ and $N$ are two finitely
generated modules over $\A$, then the canonical traces $\Tr_\tau$ on
$\B(M)$, $\B(N)$ and on $\B(M\oplus N)$ are compatible
in the following sense:
$$\Tr_\tau \left(\matrix
A&B\\
C&D
\endmatrix \right)
= \Tr_\tau (A) \ +\ \Tr_\tau(D),\tag3$$
for all $A\in \B(M),\quad D\in\B(N)$ and any
morphisms $B:M\to N$, and $C:N\to M$.
Note that the {\it von Neumann dimension}
of a Hilbertian submodule $N$ of $M$ is defined
as $\dim_{\tau}(M)=\Tr_{\tau}(P_N)$
where $P_N$ is the orthogonal projection onto $N$.

\subheading{1.4. Fuglede-Kadison determinant for
Hilbertian modules}

Let $\GL(M)$ denote the group of all invertible elements
of the algebra $\B(M)$ equipped with
the norm topology. With this topology it is a Banach Lie group whose
 Lie algebra may be identified with the commutant $\B(M)$. The canonical
trace $\Tr_{\tau}$ on the commutant $\B(M)$
 is a homomorphism of the Lie algebra $\B(M)$ into
 $\C$ and by standard theorems, it defines a group homomorphism
of the universal covering group of $\GL(M)$ into $\C$. This approach
leads to following construction of the Fuglede-Kadison determinant,
compare \cite{HS}.

\proclaim{1.5. Theorem} There exists a function
$\Det_{\tau}: \GL(M)\to \R^{>0}$
(called the Fuglede-Kadison determinant) whose key properties are:
\roster
\item $\Det_{\tau}$ is a group homomorphism and
is continuous if
$\GL(M)$ is supplied with the norm topology;
\item If $A_t$ for $t\in [0,1]$
is a continuous piecewise smooth path in $\GL(M)$ then
$$
\log \lbrack\frac{\Det_{\tau}(A_1)}{\Det_{\tau}(A_0)}\rbrack\ =\
\int_0^1\Re\Tr_{\tau}\lbrack A_t^{-1}A_t^\prime\rbrack dt.\tag4
$$
Here $\Re$ denotes the real part and $A_t^\prime$ denotes the
derivative of $A_t$ with respect to $t$.
\item Let $M$ and $N$ be two finitely generated modules over $\A$, and
$A\in\GL(M)$ and $B\in\GL(N)$ two automorphisms, and
$\gamma:N\to M$ be a homomorphism. Then the map given by the matrix
$$\left(\matrix
A&\gamma\\
0&B
\endmatrix\right)$$
belongs to $\GL(M\oplus N)$ and
$$
\Det_\tau\left(\matrix
A&\gamma\\
0&B
\endmatrix\right)\ =\
\Det_\tau(A)\cdot\Det_\tau(B)\tag5
$$
\endroster
\endproclaim

Given an operator $A\in\GL(M)$, there is a continuous piecewise
smooth path
$A_t\in\GL(M)$ with $t\in [0,1]$ such that $A_0=I$ and $A_1=A$ (it is well
known that the group $\GL(M)$ is pathwise connected, cf. \cite{Di}).
Then from (4) we have the formula:
$$\log \Det_{\tau}(A)\ =\
\int_0^1\Re\Tr_{\tau}\lbrack A_t^{-1}A_t^\prime\rbrack dt.\tag6$$
This integral does not depend on the choice of the path.
 As an example  consider the following situation.
Suppose that a self-adjoint operator $A\in \GL(M)$
has spectral resolution
$$A\ =\ \int_0^\infty \lambda dE_\lambda\tag7$$
where $dE_\lambda$ is the spectral measure.
Then we can choose the path
$$A_t\ =\ t(A-I)\ +\ I,\quad t\in [0,1]$$
joining $A$ with $I$ inside $\GL(M)$. Applying (6) we obtain
$$
\log \Det_\tau(A)\ = \int_0^\infty \ln\lambda d\phi_\lambda\tag8
$$
where $\phi_\lambda = \Tr_\tau E_\lambda$ is the spectral density function.

\subheading{1.6. Operators of determinant class}.

Following \cite{BFKM} and \cite{CFM} we extend the
previous ideas to a wider class of operators. An operator $A$ as in (7)
is said to be $D-class$ ($D$ for determinant) if
$$\int_0^\infty \ln\lambda d\phi_\lambda > -\infty\tag9$$
A scalar product
$\langle v,w\rangle=\langle Av,w\rangle_1$
is said to be $D-admissible$ if $A$ is $D$-class and $\langle\ ,\ \rangle_1$
is any admissible scalar product.
The Fuglede-Kadison determinant extends to such operators
via the formula:
$$\Det_\tau(A)=\exp[\int_0^\infty \ln\lambda d\phi_\lambda].\tag10$$

\subheading{1.7. Determinant line of a Hilbertian module}

 For a Hilbertian module $M$ we defined in \cite{CFM} the
determinant line $\det(M)$
as a real vector
space generated by symbols
$\<\ ,\ \>$, one for any admissible scalar product on $M$, subject to the
following relations: for any pair  $\<\ ,\ \>_1$ and $\<\ ,\ \>_2$
of admissible scalar products on $M$ we require
$$\<\ ,\ \>_2\ =\ \sqrt {\Det_\tau(A)}^{\ -1} \cdot \<\ ,\ \>_1,\tag11$$
where $A\in \GL(M)\cap\B(M)$ is such that
$\<v,w\>_2\ =\ \<Av,w\>_1$
for all $v,w\in M$.
It is not difficult to see  that $\det(M)$ {\it is one-dimensional
generated by the
symbol $\<\ ,\ \>$ of any admissible scalar product on $M$}.
Note also, that the real line has {\it the canonical orientation}, since the
transition coefficient $\sqrt {\Det_\tau(A)}$ is always positive. Thus we
may speak of {\it positive and negative} elements of $\det(M)$.
We think of elements of $\det(M)$ as ``volume forms" on $M$.
If $M$ is trivial module, $M=0$, then we set $\det(M)=\R$, by definition.

 Given two finitely generated Hilbertian modules
$M$ and $N$ over $\A$,
with admissible scalar products $\<\ ,\ \>_M$ and $\<\ ,\ \>_N$
respectively, we may obviously define the
scalar product $\<\ ,\ \>_M \oplus \<\ ,\ \>_N$ on the direct sum. This
defines the isomorphism
$$\det(M)\otimes\det(N)\to\det(M\oplus N).\tag12$$
By property (5) of the Fuglede-Kadison determinant
it is easy to show that
 this homomorphism does not depend on the
choice of the metrics
$\<\ ,\ \>_M$ and $\<\ ,\ \>_N$ and
 preserves the orientations.
Note that, any isomorphism $f:M\to N$ between
finitely generated
Hilbertian modules induces canonically an orientation
preserving isomorphism of the determinant
lines
$f^\ast:\det(M)\to\det(N).$
Indeed, if $\<\ ,\ \>_M$ is an admissible scalar product on
$M$ then set
$$f^\ast(\<\ ,\ \>_M)= \<\ ,\ \>_N,\tag13$$
where $\<\ ,\ \>_N$ is the scalar product on $N$ given by
$\<v,w\>_N=\<f^{-1}(v),f^{-1}(w)\>_M$ for $v,w\in N$.
 This definition does not depend on the choice
of the scalar product $\<\ ,\ \>_M$ on $M$: if we have a different admissible
scalar product $\<\ ,\ \>_M^\prime$ on $M$, where $\<v,w\>_M^\prime =
\<A(v),w\>_M$ with $A\in\GL(M)$ then the induced scalar product on $N$
will be
$$\<v,w\>_N^\prime\ =\ \<(f^{-1}Af)v,w\>_N$$
and our statement follows from property (5) of the Fuglede-Kadison
determinant. Finally we note the
{\it functorial} property:
if $f:M\to N$ and $g:N\to L$ are two isomorphisms between finitely generated
Hilbertian modules then
$(g\circ f)^\ast\ =\ g^\ast\circ f^\ast.$

\proclaim{1.8. Proposition} If $f:M\to M$ is an automorphism of a
finitely generated Hilbertian module $M$, $f\in\GL(M)$, then the induced
homomorphism $f^\ast:\det(M)\to \det(M)$ coincides with the multiplication
by $\Det_\tau(f)\in \R^{>0}$. Furthermore any exact sequence
$$0\to M^\prime@>{\alpha}>>M@>{\beta}>>M^{\prime\prime}\to 0$$
of finitely generated Hilbertian modules determines canonically an isomorphism
$$\det(M^\prime)\otimes\det(M^{\prime\prime})\ \to \det(M),$$
which preserves the orientation of the determinant lines.
\endproclaim

\subheading{1.9. Extension to $D$-admissible scalar products}

Any $D$-admissible scalar product
 determines a
 non-zero element of the determinant
line det($M$) namely
$\Det_\tau(A)^{-1/2} \langle\ ,\ \rangle_1.$
A $D-admissible$ isomorphism $f:M\to N$ is one
for which the inner product $\<v,w\>_M=\<f(v),f(w)\>_N$
on $M$ is $D$-admissible for some and hence any admissible
inner product on $N$. Proposition 1.9 extends to
$D$-admissible isomorphisms and to the obvious notion of
$D$-admissible exact sequence.

\heading{\bf \S 2. Holomorphic Hilbertian $\A$-bundles Bundles and
$\A$-linear Connections}\endheading

In this section, we define Hilbertian $\A$-bundles and
$\A$-linear connections on these. The definition of ($\A$-linear) connection
is tricky in the infinite dimensional case, if one wants to be able to
horizontally lift curves. We use some fundamental theorems in von Neumann
algebras to make sense of our definition.
We also define holomorphic Hilbertian $\A$-bundles bundles and
holomorphic $\A$-linear connections on these.

\subheading{2.1. Hilbertian $\A$-bundles}
A Hilbertian $\A$-bundle with fibre $M$ over $X$ is given by the
following data.
\roster
\item $p:\E\to X$ a smooth bundle of topological vector spaces,
possibly infinite dimensional, such that each fibre $p^{-1}(x),\ x\in X$ is
a separable Hilbertian space (cf.\cite{Lang}).
\item There is a smooth fibrewise action
$\A\times \E\to\E$ which endows each fibre
$p^{-1}(x),\ x\in X$ with a Hilbertian $\A$-module structure, such
that for all $x\in X, p^{-1}(x)$ is isomorphic to $M$ as Hilbertian
$\A$-modules.
\item There is a local trivializing cover of $p : \E\to X$ which
intertwines the $\A$-actions.  More precisely, there is an open
cover $\{U_\alpha\}$ of $X$ such that for each $\alpha$, there is a smooth
isomorphism
$$
   \tau_\alpha : p^{-1}(U_\alpha) \to U_\alpha\times M
$$
which intertwines the $\A$-actions on $p^{-1}(U_\alpha) \subset
\E$ and on $U_\alpha\times M$, and such that ${\operatorname{pr}_1}\circ
\tau_\alpha = p$,
where $\operatorname{pr}_1 : U_\alpha\times M\to U_\alpha$ denotes the
projection onto the first factor.  The restriction of $\tau_\alpha$
$$
   \tau_\alpha : p^{-1}(x) \to \{x\}\times M
$$
is the isomorphism of Hilbertian $\A$-modules $\forall x\in
U_\alpha$, as given in (2).
\endroster

\subheading{2.2. Remarks}
If $\{U_\alpha\}$ is a trivializing open cover of $p : \E\to X$,
then the isomorphisms
$$
   \tau_\beta \circ \tau_\alpha^{-1} : (U_\alpha \cap U_\beta) \times M
      \to (U_\alpha\cap U_\beta)\times M
$$
are of the form $\tau_\beta\circ\tau_\alpha^{-1} = (\operatorname{id},
g_{\alpha\beta})$ where $g_{\alpha\beta} : U_\alpha\cap U_\beta \to
\operatorname{GL}(M)$ are smooth maps and are called the
transition functions of $p : \E\to X$, and they satisfy the
cocycle identity
$$
   g_{\alpha\beta}g_{\beta\gamma}g_{\gamma\alpha}=1\quad
      \forall \alpha,\beta,\gamma.
$$
Now suppose that $\{U_\alpha\}_\alpha$ is an open cover
of $X$, and on each
intersection $U_\alpha\cap U_\beta$, we are given smooth maps
$$
   g_{\alpha\beta} : U_\alpha \cap U_\beta \to
      \operatorname{GL}(M)
$$
satisfying $g_{\alpha\beta} g_{\beta\gamma} g_{\gamma\alpha}=1$ on
$U_\alpha\cap U_\beta \cap U_\gamma$ and $g_{\alpha\alpha}=1$ on
$U_\alpha$, then one can construct a Hilbertian $\A$-bundle $p :
\E\to X$ via the clutching construction viz, consider the disjoint
union $\tilde{\E}=\bigcup_\alpha(U_\alpha\times M)$ with the
product topology, and define the equivalence relation $\sim$ on
$\tilde{\E}$ by $(x,v)\sim (y,w)$ for $(x,v)\in U_\alpha\times M$
and $(y,w)\in U_\beta\times M$ if and only if $x=y$ and
$w=g_{\alpha\beta}(x)v$.  Then the quotient $\tilde{\E}/\sim =
\E\to X$ is easily checked to be a Hilbertian $\A$-bundle
over $X$

\subheading{2.3. Remarks}
This definition generalizes and is compatible with Breuer's definition of
Hilbert $\A$-bundles (cf.\cite{B}, \cite{BFKM}) and also with Lang's
definition \cite{Lang}, where the action of the von Neumann algebra is not
considered.
Actually Breuer \cite{B} considers von Neumann algebras $\A$ which are not
necessarily finite.

\subheading{2.4. Examples}
(a). It follows from Breuer's work (\cite{B}) that there are many examples of
Hilbertian $\A$-bundles, even in the case of simply connected
manifolds.  For example, on the 2-sphere $S^2$, the isomorphism classes of
Hilbertian $\A$-bundles with fibre $\ell^2(A)$, are in 1-1
correspondence with homotopy classes of maps from $S^1$ to
$\operatorname{GL}(\ell^2(A))$.  If $\A$ is a type $II_1$ factor, then by a
result of Araki, Smith and Smith \cite{ASS}, it follows that the
isomorphism classes
of Hilbertian $\A$-bundle
 over $S^2$ is isomorphic to
$\R$ (considered as a discrete group).

(b) Let $\E\to X$ be a Hilbertian $\A$-bundle over $X$. Then
$\Lambda^jT^*_\C X\otimes\E$ is also a Hilbertian $\A$-bundle over $X$, where
$\Lambda^jT^*_\C X$ denotes the jth exterior power of the complexified
cotangent bundle of $X$.
This can be seen as follows.
Let
$$
   g_{\alpha\beta} : U_\alpha \cap U_\beta \to GL(M)
$$
denote the transition functions of the Hilbertian $\A$-bundle $\E$ with fibre
$M$, and
$$
   g'_{\alpha\beta} : U_\alpha \cap U_\beta \to GL(r, \C)
$$
denote the transition functions of the $\C$ bundle $\Lambda^j T^*_\C X\to X$.
Then
$$
   g''_{\alpha\beta} : U_\alpha \cap U_\beta \to GL(\C^r\otimes M)
$$
denotes the transition functions of the Hilbertian $\A$-bundle
$\Lambda^jT^*_\C X\otimes\E$ with fibre $\C^r\otimes M$.

\subheading{2.5. Sections of Hilbertian $\A$-bundles}
A section of a Hilbertian $\A$-bundle $p : \E\to X$ is a
smooth map ${s} : X\to\E$ such that $p\circ {s}$ is
the identity map on $X$.  Let $\{ U_\alpha \}_\alpha$ be a local trivialization
of $p : \E \to X$.  Then a smooth section ${s}$ is given
on $U_\alpha$ by a smooth map $s_\alpha : U_\alpha\to M$.  On $U_\alpha\cap
U_\beta$ one has the relation ${s}_\alpha = g_{\alpha\beta}
{s}_\beta$.

\subheading{2.6. $\A$-linear connections on Hilbertian $\A$-bundles}
An $\A$-linear connection on a Hilbertian $\A$-bundle $p : \E\to X$ is
an $\A$-morphism
$$
   \nabla : \Omega^j (X,\E) \to \Omega^{j+1}(X,\E)
$$
such that for any $A\in\Omega^0
(X,\End_{\A}(\E))$ and
$w\in\Omega^j(X,\E)$, there is
$\nabla A\in\Omega^1(X,\End_{\A}(\E))$
such that
$$
   \nabla(Aw) - A(\nabla w) = (\nabla A)w.
$$
Here $\Omega^j (X,\E)$ denotes the space of smooth sections of the
Hilbertian $\A$-bundle $\Lambda^jT^*_\C X\otimes\E$, and
$\Omega^1(X,\End_{\A}(\E))$
denotes the space of smooth sections of the
Hilbertian $\A$-bundle $T^*_\C X\otimes\End_{\A}(\E)$

\subheading{2.7. Remarks}
Let $V$ be a vector field on $X$.  Then
$$
   \nabla_V A \in \Omega^0(X, \End_{\A}(\E))
$$

\proclaim{2.8. Proposition}
Let $\nabla, \nabla'$ be two connections on the Hilbertian
$\A$-bundle $p : \E\to X$ with fibre $M$.  Then
$$
   \nabla - \nabla' \in \Omega^1 (X,\End_{\A}\E)
$$
\endproclaim

\demo{Proof}
Let $V$ be a vector field on $X$.  Then $\delta_V = \nabla_V - \nabla'_V$
in $C^\infty(X)$ linear, and hence by (\cite{Lang}) is defined pointwise.
$(\delta_V)_x$ is a derivation on the von Neumann algebra
$\End_{\A}(\E_x)$.  Since $(\delta_V)_x$ is everywhere defined, by Lemma 3,
part III, chapter 9 of  (\cite{Dix}), $(\delta_V)_x$ is
bounded.  By Theorem 1, part III, chapter 9 of (\cite{Dix}), there is an element
$B_x(V)\in \End_{\A}(\E_x)$ such that $(\delta_V)_x =
ad B_x(V)$.  That is, $x\to ad B_x(V)$ is
smooth.  The remainder of the proof establishes that there is a
smooth choice $x\to \Tilde{\Tilde{B}}_x(V)$ such that $ad
\Tilde{\Tilde{B}}_x(V) = (\delta_V)_x$.  We first discuss the local
problem.

Let $U$ be an open subset of $X$ and $M$ be a Hilbertian $\A$-module.  Consider
the trivial bundle $U\times M\to U$ over $U$.  By Dixmier's result cited
 above, there
is a map
$$
   x\to ad B_x(V)\quad x\in U
$$
where $B_x(V)\in\End_{\A}(M)$ for all $x\in U$, such that
$$
   ad B_x(V) = (\nabla_V - \nabla'_V)_x
$$
since $\nabla, \nabla'$ are connections and $V$ is smooth, we deduce that
$x\to ad B_x(V)$ is smooth.  However, it isn't {\it a priori} clear
that one can choose $x\to B_x(V)$ to be smooth, as $B_x(V)$ is only defined
modulo the centre of the von Neumann algebra $\End_{\A}(M) = \B(M)$.
To complete the proof we need the next result.

\proclaim{2.9. Lemma}
Let $\A$ be a von Neumann algebra with centre $Z$.  Then there is a smooth
section $s : \A/Z \to \A$ to the natural projection $p : \A\to \A/Z$.
\endproclaim

\demo{Proof}
Let $Z\subset \A\subset B(\ell^2(\A))$, then since $Z$ is a type I von
Neumann algebra
and hence injective, there exists a
projection of norm 1, $ P:B(H) \to Z$ (\cite{HT}). Then $\A \cap ker P$ is a
complementary subspace
to $Z$ and one defines a section to the projection  $p : \A\to \A/Z$.
$$
   s : \A/Z \to \A\quad\text{as } s([v]) = (1-P) v.
$$
Then $s$ is smooth since it is linear.

More explicitly, given a subgroup $G$ of the unitaries in the commutant of
$Z$, $U(Z')$,
which is amenable and whose span is ultra weakly
dense in $Z'$, one can use the invariant mean on $G$ to average over the closure
of the orbit $\{uxu^*: u \in G\}$
and thus obtain a map $P$ so that $P(x)$ is this average for each $x$ and hence
commutes with every $u \in Z'$. That is, $P(x)$ is in $Z''=Z$. Such projections
are called Schwartz projections, according to Kadison. (cf. \cite{Ph}).
\enddemo\hfill$\square$

Returning now to the proof of proposition
2.8, we define the smooth map $\tilde{B}(V)$ by
$$
   \tilde{B}(V) = s \circ ad B(V),
$$
where $s: End_\A(M)/Z \to End_\A(M)$ is the section as in Lemma 2.9
(with $End_\A(M)$ replacing $\A$) . Then
clearly
$$
   \nabla_V = \nabla'_V = ad\tilde{B}(V),
$$
where $x\to\tilde{B}_x(V)$ is smooth.  This solves the problem locally.

Let $\E\to X$ be a Hilbertian bundle with fibre $M$, and
$\{U_\alpha\}$ be a trivialization of $\E\to X$.  We have seen
that on $U_\alpha$, there is a smooth section
$$
   x\to\tilde{B}_{\alpha, x}(V)\quad\text{for } x\in U_{\alpha}
$$
on $\E\big|_{U_\alpha\cap U_\beta}$, we can compare the 2 sections
obtained, $x\to\tilde{B}_{\alpha,x}(V) - \tilde{B}_{\beta, x}(V)\in Z$,
since $ad\tilde{B}_{\alpha, x}(V) = ad\tilde{B}_{\beta, x}(V)$.
Therefore we can define $\lambda_{\alpha\beta}(x) = \tilde{B}_{\alpha,x}(V)
- \tilde{B}_{\beta, x}(V)$ i.e.\ $\lambda_{\alpha\beta} : U_\alpha \cap
U_\beta \to Z$ is a Cech 1-cocycle with values in the sheaf of smooth $Z$
valued functions.  As $Z$ is contractable,
lemme 22 of \cite{DD} applies and so the 1st cohomology
with values in the sheaf of smooth $Z$ valued functions is trivial.
Therefore $\lambda_{\alpha\beta}$ is a coboundary i.e.\ there are
smooth maps
$$
   \varphi_\alpha : U_\alpha\to Z
$$
such that $\lambda_{\alpha\beta} = \varphi_\beta - \varphi_\alpha$.
Then $\{ x\to \tilde{B}_{\alpha, x}(V) + \varphi_{\alpha,x}(V)\}_\alpha$ is
a global section, since on $U_\alpha\cap U_\beta$, one has
$$
   \tilde{B}_{\alpha, x}(V) + \varphi_{\alpha, x}(V) =
      \tilde{B}_{\beta,x}(V) + \varphi_{\beta, x}(V),
$$
i.e.\ one gets a smooth section
$$
   X \to \End_{\A}(\E), \qquad x\to\Tilde{\Tilde{B}}_x(V)
$$
where $\Tilde{\Tilde{B}}_x(V) = \tilde{B}_{\alpha,x}(V) + \varphi_{\alpha,x}(V)$
for $x\in U_\alpha$.  It follows that $\Tilde{\Tilde{B}}\in \Omega^1
(X,\End_{\A}(\E)$.
\enddemo\hfill$\square$

Let $\nabla$ be a connection on $p : \E\to X$ and let $\{ U_\alpha
\}_\alpha$ be a trivialization of $p : \E\to X$.  Since
$\E\big|_{U_\alpha} \cong U_{\alpha}\times M$, one sees that the
 differential $d$ is a connection on $p : \E
\big|_{U_\alpha}\to U_\alpha$.  By proposition 2.8, $\nabla - d\in
\Omega^1(U_\alpha, \End_{\A}M)$ i.e.\ $\nabla=d+B_\alpha$ where
$B_\alpha\in\Omega^1(U_\alpha, \End_{\A}M)$.

On $U_\alpha \cap U_\beta$, one easily derives the relation
$$
B_\beta = g_{\alpha\beta}^{-1} B_\alpha g_{\alpha\beta} +
      g_{\alpha\beta}^{-1} dg_{\alpha\beta}.\tag14
$$
So a connection can also be thought of as a collection $\{d +
B_\alpha\}_\alpha$ where $B_\alpha\in\Omega^1(U_\alpha,
\End_{\A} M)$ and satisfying the relation (14) on the intersection.

\subheading{2.10. Parallel sections and horizontal lifts of curves}

Let $\nabla$ be a connection on $p : \E\to X$.  Let $p :
\E\to X$ be a Hilbertian $\A$-bundle and $I=[0,1]$ be the unit interval.
Let $\gamma : I\to X$ be a curve.  Let $\xi : I\to\E$ be a curve such that
$p_0\xi =
\gamma$.  Then $\xi$ is called a lift of $\gamma$.  $\xi$ is said to be a
{\it horizontal lift} of $\gamma$ if it is parallel along $\gamma$, that
is, if it satisfies the following equation,
$$
   \nabla_{\dot{\gamma}(t)} \xi(t)=0\quad \forall t\in I
$$
where dot denotes the derivative with respect to $t$.
In a local trivialization $U_\alpha$, the equation looks as,
$$
   \dot{\xi}(t) + B_\alpha(\dot{\gamma}(t))\xi(t) = 0\quad
      \forall t\in I\tag15
$$
where $\nabla = d+B_\alpha$ on $U_\alpha$ as before.  Since
$B_\alpha(\dot{\gamma}(t))$ is {\it bounded}, we use a theorem of
ordinary differential equations for Banach space valued functions (see prop
1.1, chapter IV in \cite{Lang}) to see that there is a unique solution to
equation (15) with initial
condition $\xi(0) = v\in M$.  It follows that a connection enables one to
lift curves horizontally.  This enables one to define a ``horizontal''
subbundle $\H$ of $T\E$, which is a complement to the
``vertical'' subbundle $p^*\E\subset T\E$.  This is how
[Lang] discusses connections on infinite dimensional vector bundles.
Conversely, given a choice of ``horizontal'' subbundle $\H$ of
$T\E$, one can define a ``covariant derivative'' (that is,
a connection) as follows.  By
hypothesis $T\E=\H\oplus p^*\E$.  Let
$\operatorname{pr}_2 : T\E\to p^*\E$ denote projection to
the 2nd factor and $\kappa : T\E \to \E$ be the
composition $p\circ \operatorname{pr}_2$ where $p :
p^*\E\to\E$.  Let $V$ be a vector field on $X$.  Define
$\nabla_V s = \kappa(D s(V))$ where $s :
X\to\E$ is a smooth section, and $D{s}$ is its
differential.  Then $\nabla$ locally has the form $\{ d + B_\alpha \}$ on a
trivialization $\{ U_\alpha \}$ of $p : \E \to X$, where $B_\alpha
\in \Omega^1(U_\alpha, \End_{\A}M)$ (see [Lang, Chapter IV, Section
3]) and it satisfies relation (14).  Therefore $\nabla$ defines
a connection on $p : \E \to X$ in the sense of 2.6.

\subheading{2.11. Holomorphic Hilbertian $\A$-bundles}

A Hilbertian $\A$-bundle $p : \E\to X$ with fibre $M$, is
said to be a holomorphic Hilbertian $\A$-bundle if the transition
functions of $p : \E\to X$,
$$
   g_{\alpha\beta} : U_\alpha \cap U_\beta \to
      \operatorname{GL}(M)
$$
are holomorphic maps.  We call $\{ U_\alpha \}_\alpha$ a holomorphic
trivialization of $p : \E\to X$.

\subheading{2.12 Remarks}
$\operatorname{GL}(M)$ is an open subset of a Banach space,
and so it is a complex manifold (of infinite dimension).

\subheading{2.13 Examples of holomorphic Hilbertian $\A$-bundles}
(a) By using the clutching construction again, we see that holomorphic
Hilbertian $\A$-bundles over $S^2$ correspond to holomorphic maps
$$
   g : A_\epsilon \to \operatorname{GL}(M)
$$
where $A_\epsilon = \{ z\in\C : 1 - \epsilon < |z| < 1+\epsilon \}$ is
an annulus, for some small $\epsilon>0$.  Therefore by 2.4, there
are many examples of holomorphic Hilbertian $\A$-bundles over
$S^2$.

(b) Let $p : \E\to X$ be a flat Hilbertian $\A$-bundle
over $X$, i.e.\ $M$ is a finitely generated $(\pi-\A)$ bimodule,
where $\varphi : \pi \to \operatorname{GL}(M)$ is the left
action of $\pi$ on $M$.  Then
$$
   \E = (M\times \tilde{X})/\sim \to X
$$
where $(v,x)\sim (\varphi(g)v, g.x)$ for $g\in\pi,\ v\in M$ and
$x\in\tilde{X}$. Let
$$
   g_{\alpha\beta} : U_\alpha \cap U_\beta \to \pi
$$
denote the transition functions of the universal cover $\tilde{X}$, which
is a principal $\pi$ bundle over $X$.  Here $\{ U_\alpha\}_\alpha$ forms an
open cover of $X$.  Since $\pi$ is a discrete group and $g_{\alpha\beta}$
is smooth, it follows that $g_{\alpha\beta}$ is locally constant, and therefore
holomorphic.  The transition functions of $\E$ are
$\varphi(g_{\alpha\beta})$, which again are locally constant, and therefore
holomorphic.

(c) Let $E\to X$ be a holomorphic $\C$-vector bundle over $X$ and
$\E\to X$ a flat Hilbertian $\A$-bundle over $X$.  Let
$$
   g_{\alpha\beta} : U_\alpha \cap U_\beta \to \operatorname{GL}(r,\C)
$$
denote the holomorphic transition functions where $\{ U_\alpha \}_\alpha$
form an open cover of $X$. Let
$$
   g'_{\alpha\beta} : U_\alpha \cap U_\beta \to
      \operatorname{GL}(M)
$$
denote the transition functions of the flat Hilbertian $\A$-bundle
$\E\to X$.  Since $\E\to X$ is flat, $g'_{\alpha\beta}$
are locally constant and this holomorphic (by the previous example).
Consider the new bundle whose transition functions are given by
$$
   g''_{\alpha\beta} \equiv g_{\alpha\beta} \otimes g'_{\alpha\beta}
      : U_\alpha \cap U_\beta \to \operatorname{GL}(\C^r\otimes M).
$$
Since the $g''_{\alpha\beta}$ are holomorphic, so is the new bundle which is
the tensor product bundle, and which is denoted by
$$
   E\otimes_{\C} \E \to X
$$
We have shown that it is a holomorphic Hilbertian $\A$-bundle over $X$,
with fibre
$\C^r \otimes M$.

\subheading{2.14. Holomorphic sections of holomorphic Hilbertian $\A$-bundles}

Let $p : \E\to X$ be a holomorphic Hilbertian
$\A$-bundle.  A section ${a} : X\to\E$ is said to
be a {\it holomorphic section} if in a holomorphic
local trivialization, $\{ U_\alpha \}_\alpha$,
the expression for ${s}$ in $U_\alpha$,
$$
   {s}_\alpha : U_\alpha \to M
$$
is a holomorphic map.  Note that $M$ is a Banach space, and therefore
a complex manifold.  On $U_\alpha\cap U_\beta$, one has the relation
$$
   {s}_\alpha = g_{\alpha\beta}{s}_\beta
$$
which is holomorphic, since $g_{\alpha\beta}$ is holomorphic.

\subheading{2.15. $\A$-linear Cauchy-Riemann operators}

Let $p : \E\to X$ be a holomorphic Hilbertian $\A$-bundle
over $X$.
With respect to the decomposition
$$
T^*_\C X = T^*X\otimes_\R \C = \big(T^{1,0}X\big)^*\oplus
\big(T^{0,1}X\big)^*, \tag16
$$
the space of smooth differential $j$-forms on $X$
with values in $\E$ decomposes as a direct sum of
spaces of smooth differential $(p,q)$-forms on $X$
with values in $\E$, where $p+q =j$.
This space, which is an $\A$ module, will be denoted by $\Omega^{p,q}(X,\E)$.

Then there is a unique operator
$$
   \bar{\partial} : \Omega^{p,q}(X,\E) \to \Omega^{p,q+1}
      (X,\E)
$$
which in any holomorphic trivialization of $p : \E\to X$, is equal
to
$$
   \bar{\partial} = \sum_{i=1}^n e(d\bar{z}^i)
      \frac{\partial}{\partial \bar{z}^i}
$$
where $e(d\bar{z}^i)$ denotes exterior multiplication by the $1$-form
$d\bar{z}^i$ and $n=\dim_{\C} X$. Note that $\bar{\partial}^2 = 0$.

\subheading{2.16. Holomorphic $\A$-linear connections}

Let $\nabla : \Omega^p(X,\E) \to \Omega^{p+1} (X,\E)$ be
an $\A$-linear connection on a holomorphic Hilbertian
$\A$-bundle $p : \E \to X$.  Then with respect to (16),
there is a decomposition
$$
   \nabla = \nabla' + \nabla''.
$$
Here
$$
   \nabla' : \Omega^{p,q} (X,\E) \to \Omega^{p+1,q}
       (X,\E)
$$
is an $\A$-morphism such that for $A \in
\Omega^0(X,\End_{\A} \E)$ and $w \in
\Omega^j (X,\E)$,
$$
   \nabla'(Aw) - A(\nabla' w) = (\nabla' A) w
$$
where $\nabla' A \in \Omega^{1,0} (X, \End_{\A}
\E)$ is the $(1,0)$ component of $\nabla A$, while
$$
   \nabla'' : \Omega^{p,q} (X, \E) \to \Omega^{p,q+1}(X,\E)
$$
is an $\A$-morphism such that
$$
   \nabla'' (Aw) - A(\nabla'' w) = (\nabla'' A)\, w
$$
where $\nabla'' A \in \Omega^{0,1} (X,\End_{\A}
\E)$ is the $(0,1)$ component of $\nabla A$.

An $\A$-linear connection $\nabla$ on a  holomorphic
Hilbertian $\A$-bundle $p : \E \to X$ is said to be a
{\it holomorphic $\A$-linear connection} if $\nabla'' =
\bar{\partial}$.  In this case, $(\nabla'')^2 = 0$.

Since every holomorphic Hilbertian $\A$-bundle has a
$\A$-linear Cauchy-Riemann operator, it follows that it also
has a holomorphic
$\A$-linear connection.

\subheading{2.17. Examples of holomorphic $\A$-linear connections}

(a) Let $\E\to X$ be a flat Hilbertian $\A$-bundle. Then $\E$
has a {\it canonical flat
$\A$-linear connection} $\nabla$ given by the de Rham exterior derivative,
where we identify the space of smooth differential $j$-forms on $X$
with values in $\E$,  denoted $\Omega^{j}(X,\E)$,
as $\pi$-invariant differential forms in
$M\otimes_\C \Omega^{j}(\widetilde X)$. Here
$M\otimes_\C \Omega^{j}(\widetilde X)$ has the diagonal action.
 (See \cite{CFM} for more details).
Since the de Rham differential $d = \bar{\partial} + {\partial}$,
it is a {\it canonical flat holomorphic $\A$-linear connection}.

(b) Let $E \to X$ be a holomorphic $\C$-vector bundle over $X$, and
$\E \to X$ a flat Hilbertian $\A$-bundle over $X$.  Then
we have seen that $E \otimes_{\C} \E \to X$ is a
holomorphic Hilbertian $\A$-bundle over $X$, with fibre
$\C^r \otimes M$.  Let $\tilde{\nabla}$ be a holomorphic connection
on $E \to X$, and let $\Tilde{\Tilde{\nabla}}$ be the canonical flat
$\A$-linear connection on $\E \to X$.  Then $\nabla =
\tilde{\nabla} \otimes 1  + 1 \otimes \Tilde{\Tilde{\nabla}}$ is easily
checked to yield a holomorphic $\A$-linear connection on the
holomorphic Hilbertian $\A$-bundle $E \otimes_{\C}
\E \to X$.

\heading{\bf \S 3. Zeta functions and $D$-class bundles}\endheading

We now have most of the notation and preliminary results we need
to generalize the classical construction of the holomorphic torsion
of D.B.Ray and I.M.Singer \cite{RS}
to the infinite dimensional case.
This section generalizes \cite{BFKM} and \cite{CFM} for the notion of a
$D$-class
holomorphic Hilbertian bundle and the definition of zeta-functions
for complexes of such bundles.

\subheading{3.1 Hermitian metrics, Hilbert $\A$ bundles, L$^2$ scalar
products and the canonical
holomorphic (Hermitian) $\A$-linear connection}
{\it A Hermitian metric $h$} on a Hilbertian $\A$-bundle $p:\E\to X$
is a smooth family of admissible scalar products on the fibers. Any
Hermitian metric on $p:\E\to X$ defines a wedge product
$$
\wedge : \Omega^{p,q}(X,\E)\otimes \Omega^{r,s}(X,\E)\rightarrow
\Omega^{p+r,q+s}(X)
$$
similar to the finite dimensional case.

Let $p:\E\to X$ be a {\it holomorphic} Hilbertian $\A$-bundle
and $h$ be a Hermitian metric on $\E$. The Hermitian metric on $p:\E\to X$
determines a {\it canonical holomorphic $\A$-linear connection} on $\E$ as
follows.
Let $\nabla$ be a holomorphic $\A$-linear connection on $\E$ which
preserves the Hermitian metric $\E$, that is,
$$
d h(\xi,\eta) = h(\nabla\xi,\eta) + h(\xi,\nabla\eta)
$$
where $\xi$ and $\eta$ are smooth sections of $\E$. Equating forms of the
same type, one has
$$
\partial h(\xi,\eta) = h(\nabla'\xi,\eta) + h(\xi,\nabla''\eta)
$$
and
$$
\bar{\partial} h(\xi,\eta) = h(\nabla''\xi,\eta) + h(\xi,\nabla'\eta).
$$
Since $\nabla'' = \bar{\partial}$, we see that a choice of Hermitian
metric determines a holomorphic $\A$-linear connection, which is called the
{\it canonical holomorphic $\A$-linear connection}.

The Hermitian metric on $p:\E\to X$ together with a Hermitian metric
on $X$ determines a scalar product on $\Omega^{p,q}(X,\E)$ in the
standard way; namely, using the Hodge star operator
$$\ast:\Omega^{p,q}(X,\E)\to \Omega^{n-q, n-p}(X,\E)$$
one sets
$$(\omega,\omega^\prime)\ =\ \int_X \omega\wedge\ast\overline\omega^\prime $$

With this scalar product $\Omega^i(X,\E)$
becomes a pre-Hilbert space. Define the space of
$L^2$ differential ${p,q}$-forms on $X$ with coefficients in
$\E$, denoted $\Omega_{(2)}^{p,q}(X,\E)$, to be the Hilbert
space completion of $\Omega^{p,q}(X,\E)$. We will tend to ignore the
scalar product on $\Omega_{(2)}^{p,q}(X,\E)$ and view it as an
infinite Hilbertian $\A$ module.

\subheading{3.2 Reduced L$^2$ Dolbeault cohomology}
Given a holomorphic Hilbertian $\A$ bundle $p:\E\to X$ together with a
Hermitian metric on $\E$, one
defines the {\it reduced $L^2$ Dolbeault cohomology with coefficients
in $\E$} as the quotient
$$
H^{{p,q}}(X,\E)=\frac{\ker
{\nabla''}/\Omega^{{p,q}}_{(2)}(X,\E)}{\cl(\im\;{\nabla''}/
\Omega^{p,q-1}(X,\E))},
$$
where the Cauchy-Riemann operator ${\nabla''}$ is associated to the canonical
$\A$-linear connection $\nabla$ on $\E$. ${\nabla''}$ on $\E$ extends
to an unbounded, densely defined operator $\Omega_{(2)}^{p,q}(X,\M)\to
\Omega_{(2)}^{p,q+1}(X,\M)$.
Then $H^{p,q}(X,\E)$ is naturally defined as a Hilbertian module over $\A$.
It can also be considered as the cohomology of $X$ with coefficients in a
locally constant sheaf, determined by $\E$.

\subheading{3.3 Hodge decomposition}
The Laplacian $\square_{p,q}$ acting
on $L^2$ $\M$-valued $(p,q)$-forms on $X$ is defined to be
$$
\square_{p,q}={\nabla''} {\nabla''}^{*} + {\nabla''}^*
{\nabla''}:\Omega_{(2)}^{p,q}(X,\M)
\rightarrow\Omega_{(2)}^{p,q}(X,\M)
$$
where ${\nabla''}^{*}$ denotes the formal adjoint of ${\nabla''}$ with
respect to
the
$L^{2}$ scalar product on $\Omega_{(2)}^{p,q}(X,\M)$.
Note that by definition, the
Laplacian is a formally self-adjoint operator which is densely defined. We
also denote by $\square_{p,q}$ the self adjoint extension of the Laplacian.

Let ${\H}^{{p,q}}(X,\M)$ denote the closed subspace of $L^{2}$ harmonic
${p,q}$-forms with coefficients in $\M$, that is, the kernel of
$\square_{p,q}$. Note that ${\H}^{p,q}(X,\M)$ is a Hilbertian ${\A}$-module.
By elliptic regularity (cf. section 2, \cite{BFKM}), one sees that
${\H}^{p,q}(X,\M) \subset \Omega^{p,q}(X,\E)$, that is, every $L^2$ harmonic
$(p,q)$-form with coefficients in $\M$ is smooth.
Standard arguments then show that one has the following Hodge decomposition
(cf. \cite{D}; section 4, \cite{BFKM} and also section 3, \cite{GS})
$$
\Omega_{(2)}^{p,q}(X,\M) = {\H}^{p,q}(X,\M) \oplus \cl(\im\;{\nabla''}/
\Omega^{p,q-1}(X,\E)) \oplus \cl(\im\;{\nabla''}^*/
\Omega^{p,q+1}(X,\E)).
$$
Therefore it follows that the natural map
$
{\H}^{{p,q}}(X,\M)\rightarrow H^{{p,q}}(X,\M)
$
is an isomorphism Hilbertian ${\A}$-modules.
The corresponding $L^2$ Betti numbers are denoted by
$$
b^{{p,q}}(X,\M) = \dim_\tau \left( H^{{p,q}}(X,\M)\right).
$$

\subheading{3.4 Definition }
Let $\square_{p,q} = \int_0^\infty \lambda dE_{p,q}(\lambda)$ denote the
spectral
decomposition of the Laplacian. The {\it spectral density function} is
defined to be
$N_{p,q}(\lambda) = \Tr_\tau(E_{p,q}(\lambda))$ and the {\it theta function}
is defined to be
$\theta_{p,q}(t) = \int_0^\infty e^{-t\lambda} dN_{p,q}(\lambda)=
\Tr_\tau(e^{-t\square_{p,q}}) -
b^{{p,q}}(X,\M)$. Here we use the well known fact that the projection
$E_{p,q}(\lambda)$ and the heat operator $e^{-t\square_{p,q}}$ have smooth
Schwartz kernels which are smooth sections of a bundle over $X\times X$ with
fiber the commutant of $M$, cf. \cite{BFKM}, \cite{GS}, \cite{Luk}.
The symbol $\Tr_\tau$ denotes
application of the canonical trace
on the commutant to the restriction of the kernels to the
diagonal followed by integration over the manifold $X$. This is a trace;
it vanishes
on commutators of smoothing operators.
See also \cite{M},
\cite{L}  and \cite{GS} for the case of the flat bundle defined by the
regular representation of the fundamental group.

\subheading{3.5 Definition} A holomorphic Hilbertian $\A$-bundle
$\M\rightarrow X$
together with a choice of Hermitian metric $h$ on $\E$, is said to be  {\it
D-class} if
$$
   \int_0^1 \log (\lambda) dN_{p,q}(\lambda) >\,-\infty
$$
or equivalently
$$
\int_1^\infty t^{-1} \theta_{p,q}(t) dt < \infty
$$
for all ${p,q}=0,....,n$.
Note that the
$D$-class property of a holomorphic Hilbertian $\A$ bundle does
not depend on the choice of metrics $g$ on $X$ and $h$ on $\E$.

For the most of the paper, we make the assumption that the holomorphic
Hilbertian $\A$-bundle
$\M\rightarrow X$ is D-class. Under this assumption, we
will next define and study the zeta function of the Laplacian $\square_{p,q}$
acting on $\E$ valued $L^2$
differential forms on $X$.

\subheading{3.6 Definition} For $\lambda>0$
the {\it zeta function of the Laplacian} $\square_{p,q}$ is defined on the
half-plane
$\Re(s)>n$ as
$$
\zeta_{{p,q}}(s,\lambda, \E)  = \frac{1}{\Gamma(s)}
\int_{0}^{\infty}t^{s-1}e^{-\lambda t} \theta_{{p,q}}(t)dt.
\tag17
$$

\proclaim{3.7 Lemma }
$\zeta_{{p,q}}(s,\lambda,\E)$ is a holomorphic function in the half-plane $
\Re(s)>n$ (where $n=\dim_\C X$) and has a meromorphic continuation to
$\C$ with no pole at $s=0$. If we assume that the holomorphic
Hilbertian $\A$-bundle $\M\rightarrow X$ is D-class then
$\lim_{\lambda\rightarrow 0} \zeta'_{{p,q}}(0,\lambda, \E)$
exists (where the prime denotes differentiation with respect to $s$)
\endproclaim

\demo{Proof} There is an asymptotic expansion as
$t\rightarrow 0^{+}$ of the trace of the heat kernel
$\Tr_\tau(e^{-t\square_{p,q}})$ (cf. \cite{BFKM} and chapter 13 \cite{R}),
$$
\Tr_\tau(e^{-t\square_{p,q}})\sim t^{-n}\sum_{i=0}^{\infty}t^{i}c_{i,{p,q}}
\tag18
$$
In particular, $\Tr_\tau(e^{-t\square_{p,q}})\leq C  t^{-n}$ for $
0<t\leq 1$. From this we deduce that $\zeta_{{p,q}}(s,\lambda,\E)$ is well
defined on the
half-plane $\Re(s)>n$ and it is holomorphic there.
The meromorphic continuation of $\zeta_{{p,q}}(s,\lambda,\E)$
to the half-plane $
\Re(s)>n-N$ is obtained by considering the first $N$ terms of the
small time asymptotic expansion (18) of $\Tr_\tau(e^{-t\square_{p,q}})$,
$$
\align
\zeta_{{p,q}}(s,\lambda,\E)\ & =
-\;\sum_j\frac{b^{{p,q}}(X,\E)(-\lambda)^j}{(s+j)j!}+\frac{1}{\Gamma(s)}
\left[
\sum_{0\leq i+j\leq N}  \frac{(-\lambda)^jc_{i,{p,q}}}{(s+i+j-{n})j!}
+ R_{N}(s,\lambda) \right]\\
& +\frac{1}{\Gamma(s)}\int_1^\infty t^{s-1}
\theta_{p,q}(t)e^{-t\lambda} dt
\tag19
\endalign
$$
where $R_{N}(s,\lambda)$ is holomorphic in the half plane $\Re(s)>n-N$
with a meromorphic extension to a neighbourhood of $s=0$.
Since the Gamma function has a simple pole at $s=0$, we observe
that the meromorphic continuation of $\zeta_{{p,q}}
(s,\lambda,\E)$ has no pole at $s=0$. The last part of the lemma
now follows
cf  \cite{BFKM}.

\enddemo\hfill$\square$

Let $ \zeta'_{{p,q}}(0,0, \E)=
\lim_{\lambda\rightarrow 0} \zeta'_{{p,q}}(0,\lambda, \E)$.
The following corollary is clear from (19).

\proclaim{3.8 Corollary }
One has
$$
\zeta_{{p,q}}(0,0,\E) = - b^{{p,q}}(X,\M)+ c_{n,{p,q}}
$$
where $c_{n,{p,q}}$ is the ${n}$-th coefficient in the small time
asymptotic expansion of the theta function, cf. (18).
\endproclaim

\heading{\bf \S 4. Holomorphic $L^2$-torsion}\endheading

In this section, we define and study the generalization of Ray-Singer
holomorphic torsion
to the case of holomorphic Hilbertian $\A$-bundles.
{\it For the rest of the section, we make the assumption that the holomorphic
Hilbertian $\A$-bundle $\M\rightarrow X$ is D-class}.
Given a Hermitian manifold $X$ and a metric on a
holomorphic Hilbertian $\A$-bundle $\M$ over $X$ with fibre a Hilbertian ${\A}$
module $M$, the holomorphic $L^2$ torsion $\rho_\E^p$ defined in this section is
a {\it positive element of the determinant line}
$$
\det(H^{p,*}(X,\M)).
$$
We also prove a variational formula for the holomorphic $L^2$ torsion.

\subheading{4.1. The construction of holomorphic $L^2$ torsion} Let $(X,g)$ be
a compact, connected Hermitian manifold of complex dimension $n$ with
$\pi=\pi_{1}(X)$. Let $\E\to X$ be a holomorphic Hilbertian $\A$-bundle over $X$
with fibre $M$ and let $h$ be a Hermitian metric on $\M$. We assume that
$\E$ is of $D$-class.

As before,
let $H^{{p,q}}(X,\M)$
denote the $L^{2}$ cohomology groups of $X$ with coefficients in $\M$.
Then we know that $H^{{p,q}}(X,\M)$ is a Hilbertian ${\Cal A}$-module.
If ${\Cal H}^{{p,q}}(X,\M)$ denotes the space of $L^{2}$ harmonic
${p,q}$-forms with coefficients in $\M$, then it is a Hilbert
${\Cal A}$-module with the admissible scalar product induced from
$\Omega_{(2)}^{p,q}(X,\M)$.  By the Hodge theorem, the natural map
$$
{\Cal H}^{{p,q}}(X,\M)\rightarrow H^{{p,q}}(X,\M)
$$
is an isomorphism of Hilbertian ${\A}$-modules. Thus, we may identify
these modules via this isomorphism, or equivalently,
we may say that this isomorphism defines an admissible scalar product
on the reduced $L^2$ cohomology  $H^{{p,q}}(X,\M)$.
  These admissible scalar products on $H^{{p,q}}(X,\M)$ for all ${p,q}$,
determine elements of the determinant lines $\det(H^{{p,q}}(X,\M))$
 and thus, their product in
$$\det(H^{p,\ast}(X,\M))
=\prod_{q=0}^n \det(H^{{p,q}}(X,\M))^{(-1)^q}$$
is defined. This last element we will denote $\rho^{\prime p}(g, h)$; the
notation emphasizing the dependence on the metrics $g$ and $h$.

Using the results of the previous section, we introduce the graded  zeta
function
$$
{\zeta^p}(s,\lambda,\M)=\sum_{q=0}^n (-1)^{q}q\zeta_{{p,q}}(s,\lambda,\M).
$$
It is a meromorphic function with no pole at $s=0$.
Note also that this
zeta-function depends on the choice of the trace $\tau$ and on the metrics
$g$ and $h$.

\subheading{4.2. Definition} Define the {\it holomorphic} $L^{2}$ {\it torsion}
to be the element of the determinant line
$$
\rho_{\E}^p(g,h)\in\det(H^{p,\ast}(X,\M)),\qquad
\rho_{\E}^p(g, h)=e^{\frac{1}{2} {\zeta^p}'(0,0,\M)}
\cdot
\rho^{\prime p}(g, h).
$$
where ${\zeta^p}'$ denotes the derivative with respect to $s$.
Thus, the holomorphic $L^2$ torsion is a volume form on the reduced
$L^2$ Dolbeault cohomology.

\subheading{4.3. Remarks}  1.  In the case when ${\Cal A}=\C$, we
arrive at the classical definition of the Ray-Singer-Quillen
metric on the determinant of the Dolbeault cohomology.

2. We will prove later in this section a metric
variation formula for the holomorphic $L^{2}$  torsion
as defined in 4.2. Using this, we prove that a relative version of the
holomorphic $L^2$ torsion is independent of the choice of Hermitian metric.

3. Assuming that the reduced $L^{2}$ Dolbeault cohomology $H^{p,*}(X,\M)$
vanishes, we can identify canonically the determinant line $\det({H}^{p,*} (X,
\M))$ with $\R$, and so the torsion $\rho_{\E}^p$ in this case is just
a number.

\subheading{4.4 Metric Variation Formulae}
Suppose that a holomorphic Hilbertian $\A$-bundle $\M\to X$ of $D$-class
is given. This property
does not depend on the choice of the metrics.
Consider a smooth 1-parameter family
of metrics $g_{u}$ on $X$ and $h_{u}$ on $\M$, where $u$ varies in an interval
$(-\epsilon,\epsilon)$.
Let $(,)_u$ denote the $L^2$
scalar product on $\Omega_{(2)}^{p,*}(X, \M)$ determined by $g_u$ and $h_u$.
This family determines an invertible, positive, self-adjoint bundle map
$A_u:\M\to \M$ which is uniquely determined by the relation
$$
(\omega,\omega')_u = (A_u \omega,\omega')_0
$$
for $\omega,\omega'\in \Omega_{(2)}^{p,*}(X, \M)$; it depends smoothly on $u$.

Let $\nabla$ be the canonical $\A$-linear connection on $\E$.
Define the operator
$$
D_{u}={\nabla''}+{\nabla''}^{*}_{u}:
\Omega_{(2)}^{p,*}(X,\M)\rightarrow\Omega_{(2)}^{p,*}(X,\M)
$$
where ${\nabla''}^{*}_{u}$ denotes the formal adjoint of ${\nabla''}$ with
respect to the
$L^{2}$ scalar product $(,)_u$ on $\Omega_{(2)}^{p,*}(X,\M)$.  Then
${\nabla''}^{*}_{u}=
A_u^{-1} {\nabla''}^*_0 A_u$ acting on $\Omega^{p,*}_{(2)}(X,\M)$.
Denote $Z_u = A_u^{-1}{\dot A_u}$, where the dot means the derivative with
respect to $u$.

As in 4.1, let $\zeta_u^p(s,\lambda,\M)$ denote the graded zeta function with
respect to the
metrics $g_u,\,h_u$. The scalar product $(,)_u$ induces a scalar
product on the space of harmonic
forms ${\H}^{p,*}_u(X,\M)$, and via the canonical isomorphism
${\H}^{p,*}_u(X,\M)\to
H^{p,*}(X,\M)$, it induces an admissible scalar product on the
reduced $L^2$ cohomology
$H^{p,*}(X,\M)$. Let $\rho'(u)$ denote the class
in $\det(H^{p,*}(X,\M))$ of this scalar product.  Then the
holomorphic $L^2$ torsion
with respect to the metrics $g_u,\,h_u$ is given, as in 4.2, by
$$
\rho_{\E}^p(u)=e^{{1\over 2}{\zeta^p}'_{u}(0,0,\M)}\rho^{\prime p}(u)\in
\det(H^{p,*}(X,\M)), $$
where ${\zeta^p}'$ means the derivative with respect to $s$.

\proclaim{4.5. Theorem} Let $\E\to X$ be a holomorphic Hilbert bundle of
$D$-class.
Then in the notation above, $u\mapsto \rho_{\E}^p(u)$ is a
smooth map and one has
$$
\frac{\partial}{\partial u}\rho_{\E}^p(u) = c_\E^p(u)  \rho_{\E}^p(u),
$$
where $c_\E^p(u)\in \R$ (cf. (24)) is a local
term.
\endproclaim

The proof of this theorem
will follow from two propositions which we will prove in this section.

Let $P_p(u)$ denote the orthogonal
projection from $\Omega_{(2)}^{p,*}(X, \M)$ onto $\ker D^{2}_{u}$ and
$\Tr^s_\tau(.)$ denote the graded trace,
that is the alternating sum of the von Neumann traces $\Tr_\tau$
on operators on $\Omega_{(2)}^{p,*}(X, \M)$ having smooth Schwartz kernels.

\proclaim{4.6. Proposition} Let $\E\to X$ be a holomorphic Hilbert bundle
of $D$-class.
Then in the notation above, one has
$$
\frac{\partial}{\partial u}\;{\zeta^p}'_{u}(0,0,\M)=\;\Tr^s_\tau
(Z_u P_p(u))-2 c_\E^p(u)
$$
where $c_\E^p(u)\in \R$ (cf. (24)) is a local term.
\endproclaim

\demo{Proof}  We consider the function
$$
F(u,\lambda,s) =\sum_{q=0}^{n}(-1)^{q}q\int^{\infty}_{0}t^{s-1}e^{-t\lambda}
\Tr_\tau(e^{-t
\square_{p,q}(u)}-P_{p,q}(u))dt
$$
which is defined on the half-plane $\Re(s)>n$ and is holomorphic there.
As in (18), one has for each $u$, the small time asymptotic expansion
of the heat kernel,
$$
\Tr_\tau(e^{-t\square_{{p,q}}(u)})\sim\sum_{k=0}^{\infty}c_{k,{p,q}}(u)t^{-n
+k}\quad.
\tag20
$$
Using (20), we see that $F(u,\lambda,s)$ has a meromorphic continuation to
$\C$ with
no pole at $s=0$.
This assertion is analogous to that in Lemma 2.8, and is proved by
an easy modification of that proof.

If we know that $u\to F(u,\lambda,s)$ is a smooth function then
$$
\frac{\partial}{\partial
u}\;{\zeta^p}'_{u}(0,0,\M)
= \left.\lim_{\lambda\rightarrow 0}\frac{\partial}{\partial
s}\Big(\frac{1}{\Gamma(s)}\frac{\partial}{\partial
u} F(u,\lambda,s)\Big)\right|_{s=0}
$$
by the $D$-class assumption.
Hence:
$$\frac{\partial}{\partial
u}\;{\zeta^p}'_{u}(0,0,\M)
= \left.\lim_{\lambda\rightarrow 0}\frac{\partial}{\partial u}\;F(u,\lambda,s)
\right|_{s=0}.
$$
Observing that $\Tr_\tau(P_{{p,q}}(u))=b^{{p,q}}(X,\M)$ is independent of $u$
we see that $u\to F(u,s)$ is smooth provided we can show
that $u\to \Tr_\tau(e^{-t\square_{{p,q}}(u)})$ is a smooth function. By
an application of Duhamel's principle, one has
$$
\frac{1}{u'-u}\Big(\Tr_\tau((e^{-{t\over 2}\square_{{p,q}}(u')}
-  e^{-{t\over 2}\square_{{p,q}}(u)})e^{-{t\over 2}\square_{{p,q}}(u)}) \Big)
$$
$$
= -\int_0^{t\over 2}
\Tr_\tau(e^{-s\square_{p,q}(u')} \frac{1}{u'-u}(\square_{p,q}(u') -
\square_{p,q}(u))
e^{-{t\over 2}\square_{{p,q}}(u)}e^{-({t\over 2}-s)\square_{{p,q}}(u)})ds.
\tag21
$$
Since $||\frac{1}{u'-u}(\square_{p,q}(u') - \square_{p,q}(u)) - {\dot
\square_{p,q}(u)})
e^{-{t\over 2}\square_{{p,q}}(u)}||$ is $O(u'-u)$ as $u'\to u$, one sees
that the limit as $u'\to u$ of (21) exists and
$$
\align
\Tr_\tau\Big(\Big(\frac{\partial}{\partial u}e^{-{t\over
2}\square_{{p,q}}(u)}\Big)
e^{-{t\over 2}\square_{{p,q}}(u)}\Big)
& =-\int_0^{t\over 2}
\Tr_\tau(e^{-s\square_{p,q}(u)} {\dot \square_{p,q}(u)}
e^{-{t\over 2}\square_{{p,q}}(u)}e^{-({t\over 2}-s)\square_{{p,q}}(u)})ds\\
& = -{t\over 2} \Tr_\tau( {\dot \square_{p,q}(u)} e^{-t\square_{{p,q}}(u)}).
\endalign
$$
Therefore $u\to \Tr_\tau(e^{-t\square_{{p,q}}(u)})$ is a smooth function
(and hence so is $u\to F(u,s)$) and by the semigroup property of the
heat kernel, one has
$$
\align
\frac{\partial}{\partial u}\;\Tr_\tau(e^{-t\square_{{p,q}}(u)}-P_{{p,q}}(u))& =
\frac{\partial}{\partial u}\;\Tr_\tau(e^{-t\square_{{p,q}}(u)}) \\
&= 2 \Tr_\tau\Big(\Big(\frac{\partial}{\partial u}e^{-{t\over
2}\square_{{p,q}}(u)}\Big)
e^{-{t\over 2}\square_{{p,q}}(u)}\Big)\\
&= -t {\Tr_\tau}(\dot{\square}_{{p,q}}(u)e^{-t\square_{{p,q}}(u)}).
\endalign
$$
A calculation similar to \cite{RS}, page 152 yields
$$
\dot{\square}_{{p,q}}(u)=-Z_u {\nabla''}^{*}_{u}{\nabla''}+{\nabla''}^{*}_{u}Z_u
{\nabla''}
-\;{\nabla''} Z_u {\nabla''}^{*}_{u}+{\nabla''}\;{\nabla''}^{*}_{u} Z_u.
$$
Since ${\nabla''}\;\square_{{p,q}}(u)=\square_{p,q+1}(u){\nabla''}$ and
${\nabla''}^{*}_{u}\square_{{p,q}}(u)=
\square_{p,q-1}(u){\nabla''}^{*}_{u}$ and using the fact that $\Tr_\tau$ is a
trace,
one has
$$
\align
\Tr_\tau(\dot{\square}_{{p,q}}(u)e^{-t\square_{{p,q}}(u)}) \ & =
\Tr_\tau(Z_u {\nabla''}\;{\nabla''}_{u}^ {*}e^{-t\square_{{p,q}}(u)}) -\;
\Tr_\tau(Z_u {\nabla''}^{*}_{u}{\nabla''}\;e^{-t\square_{p,q-1}(u)}) \\
 & + \;\Tr_\tau(Z_u {\nabla''}\;{\nabla''}^{*}_{u}e^{-t\square_{p,q+1}(u)})
-\;\Tr_\tau(Z_u {\nabla''}^{*}_{u}{\nabla''}_{u}e^{-t\square_{{p,q}}(u)}).
\endalign
$$
So one sees that
$$
\align
\frac{\partial}{\partial
u}\;\sum_{q=0}^{n}(-1)^{q}q\Tr_\tau(e^{-t\square_{{p,q}}(u)}-P_{{p,q}}(u))
&= -t
\sum_{q=0}^{n}(-1)^{q}q\Tr_\tau(\dot{\square}_{{p,q}}(u)e^{-t\square_{{p,q}}
(u)})
\\
& = -t \sum_{q=0}^{n}(-1)^{q}q\Tr_\tau(Z_u
{\square}_{{p,q}}(u)e^{-t\square_{{p,q}}(u)})\\
& = t \frac{\partial}{\partial t}\sum_{q=0}^{n}(-1)^{q}\Tr_\tau(Z_u e^{-t
\square_{{p,q}}(u)}).
\endalign
$$
Using this, one sees that for $\Re(s)>n$,
$$
\frac{\partial}{\partial u}\;F(u,\lambda,s) =\sum_{q=0}^{n}(-1)^{q}
\int^{\infty}_{0}t^{s}e^{-t\lambda}
\;\frac{\partial}{\partial t}\;\Tr_\tau(Z_u
(e^{-t\square_{{p,q}}(u)}
-P_{{p,q}}(u))dt
\tag22$$
Since $Z_u$ is a bounded endomorphism, by a straightforward generalization
of lemma 1.7.7 in \cite{Gi},
there is a small time asymptotic expansion
$$
\Tr_\tau(Z_u e^{-t\square_{{p,q}}(u)})\sim\sum_{k=0}^{\infty}m_{k,{p,q}}
(u)t^{-n+k}.\tag23
$$
In particular, one has
$$
|\Tr_\tau(Z_u e^{-t\square_{{p,q}}(u)})|\leq ct^{-n}
$$
for all $0< t\leq 1$.  If $\Re(s)>n$, we can integrate
 the right-hand side of
(22) by parts to get
$$
\;\sum_{q=0}^{n}(-1)^{q+1}
\int^{\infty}_{0}(st^{s-1}-\lambda t^s)e^{-t\lambda}\Tr_\tau(Z_u
(e^{-t\square_{{p,q}}(u)}-P_{{p,q}}(u)))
dt
$$
By splitting the integral into two parts, one from 0 to 1
and the other from 1 to $\infty$ and using
(23) on the first integral together with  the observations above, one gets
the following
explicit expression
for the meromorphic continuation of
$\displaystyle\frac{\partial}{\partial u}\;F(u,s)$ to the half-plane
$\Re(s)>n-N$
$$
\align
\frac{\partial}{\partial u}\;F(u,s)\ & =
\sum_{q=0}^{n}(-1)^{q}\Tr_\tau(Z_u P_{{p,q}}(u))\frac{1}{\lambda^s}
\int_0^\lambda(st^{s-1}-t^s)e^{-t}dt\\
&+
\sum_{q=0}^{n}(-1)^{q+1}\sum_{0\leq k+r\leq N}
\frac{(-\lambda)^rm_{k,{p,q}}(u)}{r!}(\frac{s} {s-n+k+r}-
\frac{\lambda}{s-n+k+r+1}) + R_{N}(u,\lambda,s)
\endalign
$$
where $R_N(u,\lambda,s)$ is holomorphic
in a neighbourhood of zero. At $s=0$
we have
$$
R_N(u,\lambda,0)= \int^{\infty}_{1}\Tr_\tau(Z_u
(e^{-t\square_{{p,q}}(u)}-P_{{p,q}}(u)))
e^{-t\lambda}dt.
$$
Thus we have
$$
\align
\frac{\partial}{\partial u}\;{\zeta'}_{u}^p(0,\lambda,\M)\ & =
\frac{\partial}{\partial u}\;F(u,0) \\
&= \sum_{q=0}^{n}(-1)^{q+1}\Big(\sum_{k+r=n}\frac{(-\lambda)^r}{r!}(1-\lambda)
m_{k,{p,q}}(u) -\Tr_\tau(Z_u P_{{p,q}}(u))\Big)
+R_N(u,\lambda,0)
\endalign
$$
Hence
$${\zeta^p}'_{u}(0,0,\M)
=  \sum_{q=0}^{n}(-1)^{q}\Tr_\tau(Z_u P_{{p,q}}(u))- 2c_\E^p(u)
$$
where
$$
c_\E^p(u)   =
\;\frac{1}{2}\;\sum_{q=0}^{n}(-1)^{q}m_{n,{p,q}}(u) \tag24
$$
This completes the proof of the proposition.
\enddemo\hfill$\square$

The 1-parameter family of scalar products on $\Omega^{p,*}_{(2)}({X}, \M)$
which are induced
by the 1-parameter family of metrics on $X$ and $\M\to X$, defines an
inclusion isomorphism of Hilbertian modules
$$
I_u : {\Cal H}^{p,\ast}_{u}(X,\M)\rightarrow H^{p,\ast}(X,\M).
$$
Here ${\Cal H}^{p,q}_{u}(X,\M)$ denotes the kernel of $\square_{p,q}(u)$.
There is an induced isomorphism of determinant lines cf. (13) and the discussion
in the paragraph above it.
$$
I_u^* : \det(H^{p,\ast}(X,\M))\rightarrow \det({\Cal
H}^{p,\ast}_{u}(X,\M)).
$$
We first identify $H^{p,\ast}(X,\M)$ with ${\Cal H}^{p,\ast}_{0}(X,\M)$.
Then $I_u$ defines a 1-parameter family of admissible scalar products on
$H^{p,\ast}(X,\M)$, which we can write explicitly as follows:
$$
\langle \eta, \eta'\rangle_u = ( P(u)\eta,
P(u)\eta')_u = (A_u  P(u)\eta, P(u)\eta')_0
$$
where $\eta,\eta'$ are harmonic forms in
${\Cal H}^{p,\ast}_{0}(X,\M)$. The relation between these scalar products
in the determinant line $\det(H^{p,\ast}(X,\M))$ is given as in 1.8  and
(11), by
$$
\langle\ ,\  \rangle_u = \prod_{q=0}^n {{\Det}_{\tau'}
(P_{p,q}(u)^\dagger A_u P_{p,q}(u))}^{{(-1)^{q+1}\over 2}}
\langle\ ,\  \rangle_0.\tag25
$$
where $P_{p,q}(u)^\dagger$ denotes the adjoint of $P_{p,q}(u)$ with respect to
the
fixed admissible scalar product $\langle\ \ ,\  \rangle_0$
and $\Tr_{\tau'}(\cdot)$ is the trace on
$H^{{p,q}}({X}, \M)$.
Using the fact that $H^{p,q}( X, \M)$ is isomorphic
to a submodule of a free Hilbertian module as is
${\Cal H}^{{p,q}}_{u}(X, \M)$, it follows that
$\Tr_{\tau'}(.)$
is equal to $\Tr_\tau(P_{p,q}(u)\cdot P_{p,q}(u))$.  We begin with the following

\proclaim{4.7. Proposition} Let $\E\to X$ be a holomorphic Hilbert bundle
of $D$-class.
Then the function $u\to P_{p,q}(u)$ is smooth and
in the notation of 4.4 and 4.5, one has
$$
\frac{\partial}{\partial u} \rho^{\prime p}(u) =
-\frac{1}{2}\Tr^s_\tau( Z_u P_p(u)) \rho^{\prime p}(u).
$$
\endproclaim

\demo{Proof} We will first prove that $u\to P_{p,q}(u)$ is a smooth function.
First consider the Hodge decomposition in the $u$-metric in the context,
$$
\Omega^{{p,q}}_{(2)}({X}, \M) = {\H}^{p,q}_u(X,\E) \oplus \cl(\im {\nabla''})
\oplus
\cl(\im {\nabla''}^*_u)
$$
and let $\pi$ denote the projection onto $\cl(\im {\nabla''})$, which
does not depend on the $u$-metric.
Let $h \in {\H}^{p,q}_0(X,\E)$ be harmonic in the $u=0$ metric. We will
arrive at a formula for $h_u \equiv P_{p,q}(u) h$, from which which
the differentiability
of $u\to P_{p,q}(u)$ will be clear. Now define $r_u$ by the equation
$$
h_u = h + r_u.
$$
Since $h_u$ is harmonic in the $u$-metric, one has
${\nabla''}^*_u(h_u) = 0$. By the formula for ${\nabla''}^*_u$ in 4.4,
one sees that ${\nabla''}^*_0 A_u(h+r_u) = 0$.
Since ${\nabla''}^*_0$ is injective on $\cl(\im {\nabla''})$, one has that
$\pi(A_u(h+r_u)) = 0$. Since $B_u\equiv \pi\,A_u \pi : \cl(\im {\nabla''})\to
\cl(\im {\nabla''})$ is an isomorphism, one sees that
$r_u = - B_u^{-1}\pi A_u(h)$ and therefore
$$
h_u = h -  B_u^{-1}\pi A_u(h).
$$
Since $u\to A_u$ is smooth, it follows that $u\to B_u$ is smooth
and by the formula above, one concludes that $u\to P_{p,q}(u)$ is
also smooth.

Observe that
$$
P_{p,q}(u)^2 = P_{p,q}(u).
$$
Differentiating with respect to $u$, one has
$$
\dot{P_{p,q}(u)} = P_{p,q}(u)\dot{P_{p,q}(u)} +
\dot{P_{p,q}(u)} P_{p,q}(u).
$$
Therefore
$$
P_{p,q}(u)\dot{P_{p,q}(u)}P_{p,q}(u) = 0.
$$
Therefore
$$
\align
{\Tr_\tau}(\dot{P_{p,q}(u)})
&= 2{\Tr_\tau}(P_{p,q}(u)\dot{P_{p,q}(u)})\\
&=2{\Tr_\tau}(P_{p,q}(u)\dot{P_{p,q}(u)}P_{p,q}(u)) \\
& =0.
\endalign
$$
A similar argument shows that the projection $P_{p,q}^\dagger(u)$
also satisfies
$$
{\Tr_\tau}(\dot{P_{p,q}^\dagger (u)})  = 0
$$

By definition,
$\rho^{\prime p}(u) = \langle\  ,\  \rangle_u \in \det\Big( H^{p,*}(X,\M)\Big)$,
and therefore by differentiating the relation (25), one has
$$
\frac{\partial}{\partial u} \rho^{\prime p}(u) =
-{1\over 2} \Tr^s_{\tau'}( C_u^{-1}\frac{\partial}{\partial u} C_u)
\rho^{\prime p}(u)
$$
where $C_u \equiv P_p(u)^\dagger A_u  P_p(u)$ and
$\Tr^s_{\tau'}(.)$ denotes the graded von Neumann trace on
$H^{p,\bullet}( X, \M)$.
Therefore one sees that
$$
\align
\frac{\partial}{\partial u} \rho^{\prime p}(u)
& =-{1\over 2} \Tr^s_{\tau'}( Z_u P_p(u) + P_p(u) \dot{P_p(u)} +  P_p^\dagger(u)
\dot{P_p^\dagger(u)})\rho^{\prime p}(u)\\
& = -{1\over 2} \Tr^s_{\tau'}( Z_u P_p(u))\rho^{\prime p}(u) =
-{1\over 2} \Tr^s_{\tau}( Z_u P_p(u))\rho^{\prime p}(u).
\endalign
$$

\enddemo\hfill$\square$

\noindent{\bf Proof of Theorem 4.5}.
By Proposition 4.7, one calculates
$$
\align
\frac{\partial}{\partial u} \rho_{\E}^p (u) \ & =
 \frac{1}{2}e^{\frac{1}{2} {\zeta^p}'_{u}(0,0,\M)}  \frac{\partial}{\partial u}
{\zeta^p}'_{u}(0,0,\M)
\rho^{\prime p}(u)
+ e^{\frac{1}{2} {\zeta^p}'_{u}(0,0,\M)} \frac{\partial}
{\partial u} \rho^{\prime p}(u) \\
& = \frac{1}{2} \Big[ \frac{\partial}{\partial u} {\zeta^p}'_{u}(0,0,\M)
- \Tr^s_\tau( Z_u P_p(u)) \Big]
 e^{\frac{1}{2} {\zeta^p}'_{u}(0,0,\M)}
\rho^{\prime p}(u)\\
& = \frac{1}{2} \Big[ \frac{\partial}{\partial u} {\zeta^p}'_{u}(0,0,\M)
-  \Tr^s_\tau( Z_u P_p(u)) \Big] \rho_{\E}^p (u).
\endalign
$$
Therefore by Proposition 4.6, one has
$$
\frac{\partial}{\partial u} \rho_{\E}^p (u) =  c_\E^p(u)  \rho_{\E}^p (u)
$$
where $c_\E^p(u)\in \R$ is as in (24).
This completes the proof of the theorem.
\hfill$\square$

\heading{\bf \S 5. Flat Hilbert $\A$-bundles and Relative Holomorphic $L^2$
Torsion}
\endheading

In this section, we define the relative holomorphic
$L^2$ torsion with respect to a pair of {\it flat Hilbert} (unitary)
$\A$-bundles
$\E$ and $\F$, and we
prove that it is independent of the choice of Hermitian metric on the
complex manifold. Thus it can be viewed as an invariant
volume form on the reduced $L^2$ cohomology  $H^{p,*}(X,\E)\oplus
H^{p,*}(X,\F)^\prime$.
In section \S 6, we will prove the relative holomorphic
$L^2$ torsion with respect to a pair of {\it flat Hilbertian}  $\A$-bundles
$\E$ and $\F$, is independent of the choice of almost K\"ahler metric on an
almost K\"ahler manifold and on the choice of Hermitian metrics on $\E$ and
$\F$.

\subheading{5.1 Relative holomorphic $L^2$ torsion}
It follows from Theorem 4.5 that the holomorphic $L^2$ torsion is {\it not}
independent of the choice of metrics on the complex manifold and on the
flat Hilbertian
bundle. Therefore in order to obtain an invariant, we now consider the
{\it relative} holomorphic $L^2$ torsion for a pair of {\it unitary} flat
Hilbertian bundles over a complex manifold. In the next section, we will
study the
{\it relative} holomorphic $L^2$ torsion for an arbitrary pair of flat
Hilbertian bundles over a complex manifold.

A {\it distance function} $r$ on a manifold $X$ is a map
$r:X\times X\to\Bbb{R}$ such that

\vskip 0.1in

(1) Its square $r^2(x,y)$ is smooth on $X\times X$.

(2) $r(x,x)=0$ and $r(x,y)>0$ if $x\neq y$.

(3) $\displaystyle{\partial^2 \over \partial x_i\partial x_j} r^2(x,y)
      \bigg|_{x=y} = g_{ij}(x)$

\vskip 0.1in

Condition (3) says essentially that $r(x,y)$ coincides with the geodesic
distance from $x$ to $y$, whenever $x$ and $y$ are close.  One can easily
construct such a function using local coordinates and a partition of unity.
Let
$$
   k(t,x,y) = c_1 t^{-n} e^{-c_2{r^2(x,y)\over t}},\quad t>0
$$
and $c_1,c_2$ are some positive constants.  Then one has the following basic
theorem about the fundamental solution of the heat equation,

\proclaim{5.2. Proposition}
The heat kernel $e^{-t\square^{\Cal E}_{p,q}}(x,y)$ is a smooth, symmetric
double form on $X$ and has the property
$$
   {\nabla''}_x e^{-t\square^{\Cal E}_{p,q}}(x,y) =
      {\nabla''}_y^* e^{-t\square^{\Cal E}_{p,q}}(x,y). \tag26
$$
It satisfies the bounds
$$
   \big| D e^{-t\square^{\Cal E}_{p,q}}(x,y) \big| \leq c_3 t^{-{1\over 2}}
      k(t,x,y)\tag27
$$
for $D={\nabla''}$ or ${\nabla''}^*$,\ $x,y$ close to each other
and $0<t\leq 1$.  Finally, there is a small time asymptotic expansion
$$
   e^{-t\square^{\Cal E}_{p,q}}(x,x) \sim \sum_{j=0}^\infty t^{-n+j}
      C_{j,p,q}(x)\tag28
$$
as $t\to 0$, where $C_{j,p,q}$ is a smooth double form on $X$, for all
$j$.
\endproclaim

\demo{Proof} The result is local, and in a local normal coordinate
neighborhood of a
point $x \in X$, where the bundle ${\Cal E}$ is also trivialized, one can
proceed
exactly as in \cite{RS1}, Proposition 5.3.
(cf. \cite{R}, \cite{BFKM})
\enddemo\hfill$\square$

\subheading{5.3} By Theorem 4.5, we see that the holomorphic $L^2$
torsion is not necessarily independent
of the choice of Hermitian metrics on $X$ and $\E\to X$. We will now study the
case when the flat Hilbertian $\A$-bundle $\E\to X$ with fiber $M$
is defined by a
{\it unitary representation} $\pi\rightarrow
{\Cal B}_{\Cal A}(M)$, that is, $M$ is a {\it unitary}
Hilbertian $(\A-\pi)$ bimodule. That is,
$$
\E\equiv (M\times\widetilde X)/\sim \to X
$$
where $(v,x)\sim(vg^{-1},gx)$ for all $g\in \pi$, $x\in \widetilde X$ and $v\in
M$.
The unitary representation defines a flat Hermitian
metric $h$ on $\E\to X$. We call such a bundle a {\it flat Hilbert bundle},
or sometimes a {\it unitary} flat Hilbertian bundle.
Then by definition (cf. 4.2), one has
$$
\rho_{\E}^p(g,h)\in\det(H^{p,*}(X,\M)).
$$

Let ${\Cal F}\to X$ be another flat Hilbert $\A$ bundle with fibre $N$,
such that
$\dim_\tau(M) = \dim_\tau(N)$. Let $\square^{\Cal E}_{p,q}(u)$ and
$\square^{\Cal F}_{p,q}(u)$
denote the Laplacians in the metric $g_u$, acting on
$\Omega^{p,q}_{(2)}(X,{\Cal E})$ and $\Omega^{p,q}_{(2)}(X,{\Cal F})$
respectively.   We first prove the following Proposition.

\proclaim{5.4. Proposition} Let $\E$ and $\F$ be a pair of flat Hilbert
bundles over $X$,
as above. Then there are positive constants $C_1, C$ such that
$$
\Big|\Tr_\tau(Z_u\exp(-t\square^\E_{p,q}(u)) ) -
\Tr_\tau(Z_u\exp(-t\square^\F_{p,q}(u)) )\Big| \leq C_1 e^{-{C\over t}}
$$
for all $0<t\le 1$.
\endproclaim

\demo{Proof}
Let $x\in X$ and assume that the ball $U_\delta = \{y\in X : r^2(x,y)<
\delta\}$ is simply connected, where $r$ is a distance function on $X$
which coincides with the geodesic distance on $U_\delta$.  Since the
Laplacian is a local operator, it follows that $\square_{p,q}^{\Cal E}$
acting on $\Omega_{(2)}^{p,q}(X,{\Cal E})$ over $U_\delta$ coincides
with $\square^{\Cal F}_{p,q}$ acting on $\Omega_{(2)}^{p,q}(X, {\Cal
F})$ over $U_\delta$.  By Duhamel's Principle and by
applying Green's theorem, one has for $x,y\in U_\delta$, one has
$$
\align
    e^{-t\square^{\Cal E}_{p,q}(u)}(x,y) -
      e^{-t\square^{\Cal F}_{p,q}(u)}(x,y)
   & = \int_0^t \int_{r^2(x,z)=\delta} \Bigl[
      e^{-(t-s)\square^{\Cal F}_{p,q}(u)}(z,y) \wedge*
      {\nabla''}^* e^{-s\square^{\Cal E}_{p,q}(u)}(x,z)  \\
   &  \phantom{= \int_0^t \int_{r^2(x,z)=\delta} \Bigl[}\
      -{\nabla''}^* e^{-s\square^{\Cal E}_{p,q}(u)}(x,z) \wedge*
      e^{-(t-s)\square^{\Cal F}_{p,q}(u)}(z,y) \\
   & \phantom{= \int_0^t \int_{r^2(x,z)=\delta} \Bigl[}\
      -e^{-s\square^{\Cal E}_{p,q}(u)}(x,z) \wedge*
      {\nabla''}^* e^{-(t-s)\square^{\Cal F}_{p,q}(u)}(z,y) \\
   &  \phantom{= \int_0^t \int_{r^2(x,z)=\delta} \Bigl[}\
      +{\nabla''}^* e^{-(t-s)\square^{\Cal F}_{p,q}(u)}(z,y) \wedge*
      e^{-s\square^{\Cal E}_{p,q}(u)}(x,z) \Bigr]
\endalign
$$
Using the basic estimate $(27)$ for heat kernels, one has
$$
   \Bigl| \operatorname{Tr}_\tau\bigl(Z_u e^{-t\square^{\Cal E}_{p,q}(u)}
      \bigr) - \operatorname{Tr}_\tau\bigl(Z_u e^{-t\square^{\Cal F}_{p,q}(u)}
      \bigr) \Bigr| \le c_1 t^{-{1\over 2}} e^{-{c_2\delta\over t}}\le C_1
e^{-{C\over t}}
$$
for all $0<t\le 1$.
\enddemo\hfill$\square$

\proclaim{5.5 Theorem} In the notation of 4.3, if
$\E$ and $\F$ are a pair of flat Hilbert bundles over $X$ which are of
$D$-class,
then the relative holomorphic $L^2$ torsion
$$
\rho_{\E,\F}^p = \rho_{\E}^p\otimes ({\rho_{\F}^p})^{-1} \in \det\Big(
H^{p,*}(X,
\E)\Big)\otimes \det\Big( H^{p,*}(X, \F)\Big)^{-1}
$$
is independent of the choice of Hermitian metric on $X$ which is needed to
define it.

\endproclaim

\demo{Proof}
Let $u\to g_u$ be a smooth family of Hermitian metrics on $X$ and
$\square^\E_{p,q}(u)$ and $\square^\F_{p,q}(u)$ denote the Laplacians
on $\E$ and $\F$ respectively, as before.

By Proposition 5.4, one has
$$
\Big|\Tr_\tau(Z_u\exp(-t\square^\E_{p,q}(u)) ) -
\Tr_\tau(Z_u\exp(-t\square^\F_{p,q}(u)) )\Big| \leq C_1 e^{-{C\over t}}
$$
as $t\to 0$. That is, $\Tr_\tau^s(Z_u\exp(-t\square^\E_{p,q}(u)) )$ and
$\Tr_\tau^s(Z_u\exp(-t\square^\F_{p,q}(u)) )$ have the same asymptotic expansion
as $t\to 0$. In particular, one has in the notation of Theorem 4.5,
$$
c_\E(u) = c_\F(u).
$$
Then the relative holomorphic $L^2$ torsion
$$
\gather
\rho_{\E,\F}^p\in\det H^{p,*}(X,\Cal{E})\otimes(\det H^{p,*}(X,\F))^{-1}\\
\rho_{\E,\F}^p(u)=\rho_\Cal{E}^p(u)\otimes(\rho_{\F}^p(u))^{-1}
\endgather
$$
satisfies
$$
\align
{\partial \over\partial u}\rho_{\E,\F}^p(u) =& \left({\partial \over\partial u}
   \rho_\Cal{E}^p(u)\right)\otimes\left(\rho_{\F}^p
   (u)\right)^{-1}-\rho_\Cal{E}^p(u)\otimes{\partial \over\partial u}
   \rho_{\F}^p(u)\cdot \rho_{\F}^p(u)^{-2}\\
=& (c_\Cal{E}(u)-c_{\F}(u)) \rho_{\E,\F}^p(u)\\
=& 0
\endalign
$$
using Theorem 4.5 and the discussion above.  This proves the theorem.
\enddemo\hfill$\square$

\heading{\bf \S 6. Determinant Line Bundles, Correspondences and Relative
Holomorphic $L^2$ Torsion} \endheading

In this section, we introduce the notion of determinant line bundles of
Hilbertian $\A$-bundles over compact manifolds. A main result in this
section is Theorem 6.8, which says that the holomorphic $L^2$ torsion
associated to a
flat Hilbertian bundle over a compact almost K\"ahler manifold, depends only on
the class of the Hermitian metric in the determinant line bundle of the flat
Hilbertian bundle. This enables us to show that a correspondence of determinant
line bundles is well defined on almost K\"ahler manifolds.
Finally, using such a correspondence of
determinant line bundles, we prove in Theorem 6.12 that the relative holomorphic
$L^2$ torsion is independent of the choices of almost K\"ahler metrics on the
complex manifold and Hermitian metrics on the pair of flat Hilbertian bundles
over the complex manifold.

\proclaim{6.1. Lemma}
The subgroup $SL(M) = \Det_\tau^{-1}(1)$ of $GL(M)$ is connected.
\endproclaim

\demo{Proof}
Let $U(M)$ denote the subgroup of all unitary elements in $GL(M)$. Recall
the standard retraction of $GL(M)$ onto $U(M)$, which is given by
$$
T_s : GL(M)\to GL(M)
$$
$$
A \to |A|^s \frac{A}{|A|}
$$
where $T_0 :GL(M)\to U(M)$ is onto and $T_1 = identity$. Clearly $U(M)
\subset SL(M)$
and the retraction $T_s$ above restricts to be a retraction of $SL(M)$ onto
$U(M)$. By the results of \cite{ASS}, it follows that $SL(M)$ is connected.
\enddemo\hfill$\square$

Let $\E\to X$ be a Hilbertian $\A$-bundle over $X$ and
$GL(\Cal{E})$ denote the space of complex $A$-linear automorphisms
of $\Cal{E}$ which induce the identity map on $X$, that is,
$GL(\Cal{E})$ is the gauge group of $\Cal{E}$.
The Fuglede-Kadison determinant, cf Theorem 1.32.
$$
\Det_\tau : GL(M)\rightarrow \Bbb{R}^+
$$
extends to a homomorphism
$$
\Det_\tau : GL(\Cal{E})\rightarrow C^\infty(X,\Bbb{R}^+)
$$
where $C^\infty(X,\Bbb{R}^+)$ denotes the space of smooth positive
functions on $X$.  This extension has all the properties listed in
theorem 1.32. Using the long exact sequence in homotopy and the Lemma
above, one has

\proclaim{6.2. Corollary} Let $\E\to X$ be a Hilbertian $\A$-bundle over $X$
(recall that $X$ is assumed to be connected). Then
the subgroup $SL(\E) = \Det_\tau^{-1}(1)$ of $GL(\E)$ is connected.
\endproclaim

\subheading{6.3. Determinant Line Bundles}
Let $\E\to X$ be a Hilbertian $\A$-bundle over $X$. Then we can define a
natural {\it determinant line bundle} of $\E$ as follows:

Let $\herm(\Cal{E})$ denote the space of all Hermitian metrics on
$\Cal{E}$.  Clearly $\herm(\Cal{E})$ is a convex set and $GL(\Cal{E})$
acts on $\herm(\Cal{E})$ by
$$
\gather
GL(\Cal{E})\times \herm(\Cal{E})\rightarrow \herm(\Cal{E})\\
(a,h)\rightarrow {\bar{a}}^t ha
\endgather
$$
That is, $(a.h)_x(v,w)=h_x(av,aw)$ for all $v, w\in\Cal{E}_x$.

The action of $GL(\Cal{E})$ on $\herm(\Cal{E})$ is transitive, that
is, one can identify $\herm(\Cal{E})$ with the quotient
$$
GL(\Cal{E})\big/ U(\Cal{E}, h_0)\qquad \text{where} \quad U(\Cal{E},h_0)
$$
is the subgroup of $GL(\Cal{E})$ which leaves $h_0\in
\herm(\Cal{E})$ invariant, that is, $U(\Cal{E},h_0)$ is the unitary
transformations with respect to $h_0$.

For a Hilbertian bundle $\Cal{E}$ over $X$, we define
$\det(\Cal{E})$ to be the real vector space generated by the symbols
$h$, one for each Hermitian metric on $\Cal{E}$, subject to the
following relations : for any pair $h_1$, $h_2$ of Hermitian metrics
on $\Cal{E}$, we write the following relation
$$
h_2=\sqrt{\Det_\tau(A)}^{\ -1} h_1
$$
where $A\in GL(\Cal{E})$ is positive, self-adjoint and satisfies
$$
h_2(v,w)=h_1(Av,w)
$$
for all $v, w\in \Cal{E}_x$.

Assume that we have three different Hermitian metrics $h_1$, $h_2$ and
$h_3$ on $\Cal{E}$.

Suppose that
$$
h_2(v,w)=h_1(Av,w) \text{ and } h_3(v,w)=h_2(Bv,w)
$$
for all $v,w\in \Cal{E}_x$ and $A,B\in GL(\Cal{E})$.  Then
$h_3(v,w)=h_1(ABv,w)$a and we have the following relations in
$\det(\Cal{E})$,
$$
\align
h_2 &= \sqrt{\Det_\tau(A)}^{\ -1} h_1\\
h_3 &= \sqrt{\Det_\tau(B)}^{\ -1} h_2\\
h_3 &= \sqrt{\Det_\tau(AB)}^{\ -1} h_1
\endalign
$$

The third relation follows from the first two, from which is follows
that $\det(\Cal{E})$ is a line bundle over $X$.

To summarize, $\det(\Cal{E})$ is a real line bundle over $X$, which
has nowhere zero sections $h$, where $h$ is any Hermitian metric on
$\Cal{E}$.  It has a canonical orientation, since the transition
functions $\sqrt{\Det_\tau(A)}^{\ -1}$ are always positive.

Non zero elements of $\det(\Cal{E})$ should be viewed as volume forms
on $\Cal{E}$.

For {\it flat} Hilbertian $\A$ bundles, the determinant line
bundle can be described in the following alternate way.

Then
$\Cal{E}=M\times_\rho\widetilde{X}$, where $\rho:\pi\rightarrow GL(M)$ is a
representation.  The associated {\it determinant line bundle} is
defined as
$$
\det{\Cal{E}}=\det(M)\times_{\Det_\tau(\rho)}\widetilde{X}.
$$
Here $\Det_\tau(\rho):\pi\rightarrow \Bbb{R}^+$ is a representation
which is defined as
$$
\Det_\tau(\rho) (\gamma) = \Det_\tau(\rho(\gamma))
$$
for $\gamma \in \pi$. Then $\det(\Cal{E})$ has the property that
$$
\det(\Cal{E})_x=\det(\Cal{E}_x) \qquad \forall x\in X.
$$
Clearly
$\det(\Cal{E})$ coincides with the construcion given in the beginning
of 6.3, and $\det(\Cal{E})$ is a {\it flat} real line bundle over $X$.

\subheading{6.4 Almost K\"ahler manifolds} A Hermitian manifold $(X,g)$
is said to be {\it almost K\"ahler}
if the K\"ahler 2-form $\omega$ is not necessarily closed, but instead
satisfies the weaker condition $\overline\partial\partial \omega = 0$.
Gauduchon (cf. \cite{Gau}) proved that every complex manifold of real
dimension less than or equal to 4, is almost K\"ahler.

Let $\nabla^B$ denote the holomorphic Hermitian connection on $TX$ with the
torsion tensor $T^B$ and curvature tensor $R^B$.  Define the smooth 3-form
$B$ by
$$
   B(U,V,W) = (T^B(U,V),W)
$$
for all $U,V,W\in TX$.  Let $\omega$ denote the K\"ahler 2-form on $X$.  Then
one has
$$
   B = i(\partial - \overline\partial) \omega.
$$
Since $X$ is almost K\"ahler, it follows that $B$ is closed and therefore the
following curvature identity holds
$$
   (R^B(U,V)W,Z) = (R^{-B}(Z,W)V,U)
$$
for all $U,V,W,Z\in TX$.  The Dolbeault operator $\sqrt{2}({\nabla''} +
\nabla^{\prime\prime*})$ is a Dirac type operator.  More precisely, Let
$\Lambda = (\det T^{\prime\prime0}X)^{1\over 2}$ and $\fs$ denote the
bundle of spinors on $X$, then as $\Bbb{Z}_2$ graded bundles on $X$, one
has
$$
   \Lambda^{p,*}T^*X\otimes {\Cal E} =
      \fs\otimes\Lambda\otimes\Lambda^{p,0} T^*X\otimes {\Cal E}.
$$
Let $\nabla^L$ denote the Levi-Civita connection on $X$ and ${\Cal D}^L$ the
Dirac
operator with respect to this connection.  Then using the connection
$\nabla^B$ on $\Lambda$ and $\Lambda^{p,0}T^*X$, the Dirac operator
${\Cal D}^L$ extends as an operator
$$
   {\Cal D}^L : \Gamma(X,\fs^+\otimes\Lambda\otimes\Lambda^{p,0}T^*X\otimes
      {\Cal E}) \to \Gamma(X,\fs^-\otimes\Lambda\otimes\Lambda^{p,0}
      T^*X\otimes{\Cal E})
$$
and one has the formula
$$
   \sqrt{2}(\nabla'' + \nabla^{\prime\prime*}) = {\Cal D}^L - {1\over 4}c(B)
      = {\Cal D}^L + {1\over2} \sum_{i=1}^n c(S(e_i)e_i)
$$
where $c(B)$ denotes Clifford multiplication by the 3-form $B$ and
$S=\nabla^B-\nabla^L$ is a 1-form on
$X$ with values in skew-Hermitian endomorphisms of $TX$.
We now work in a local
normal coordinate ball, where we trivialize the bundles using parallel
transport along geodesics.
Scale the metric on $X$ by $r^{-1}$ and let $I_r$ denote the operator
$2\square_{p,*} = (\sqrt{2}(\nabla'' + \nabla^{\prime\prime*}))^2$ in this
scaled metric.  In local normal coordinates, one has the following expression
for $I_r$ (cf. \cite{B})
$$
\align
   I_r &= -rg^{ij}\left(\partial_i + {1\over 4}\Gamma_{iab} c(e_a \wedge
e_b) + A_i
       + {1\over 2\sqrt{r}} c(S_{il\alpha}(e_l)e(f_\alpha)) + {1\over 4r}
      S_{i\beta\gamma} e(f_\beta\wedge f_\gamma)\right) \\
   &\quad \times\left(\partial_j + {1\over 4}\Gamma_{jab} c(e_a\wedge e_b) +
      A_j + {1\over 2\sqrt{r}} S_{jl\alpha} c(e_l)e(f_\alpha)
      \phantom{\bigl(} + {1\over 4r} S_{j\beta\gamma}
      e(f_\beta\wedge f_\gamma)\right)  \\
   & \quad + {1\over 4}rk - {1\over 2}rc(e_i\wedge e_j) L_{ij} -
      {1\over 2}e(f_\alpha\wedge f_\beta) L_{\alpha\beta}
      - \sqrt{r} c(e_i) e(f_\alpha\wedge L_{i\alpha})\\
   & \quad + rg^{ij}\Gamma_{ij}^k \left( \partial_k + {1\over 4}\Gamma_{kab}
      c(e_a\wedge e_b) + A_k +{1\over 2\sqrt{r}} S_{kl\alpha} c(e_l)
      e(f_\alpha) + {1\over 4r} S_{k\beta\gamma} e(f_\beta \wedge
f_\gamma)\right)
\endalign
$$
where $k$ denotes the scalar curvature of $X$.

Consider the heat equation on sections of
$\fs\otimes\Lambda\otimes\Lambda^{p,0}T^*X\otimes{\Cal E}$,
$$
\align
   (\partial_t + I_r) g(x,t) &= 0 \\
   g(x,0) &= g(x).
\endalign
$$
By parabolic theory, there is a fundamental solution $e^{-tI_r}(x,y)$ which
is smooth for $t>0$.  We will consider the case when $t=1,\ e^{-I_r}(x,y)$
and prove the existence of an asymptotic expansion on the diagonal, as
$r\to 0$.  A difficulty arises because of the singularities arising in the
coefficients of $I_r$, as $r\to 0$.

\proclaim{6.5. Proposition}
For some positive integer $p\ge n$, one has the following asymptotic expansion
as $r\to 0$,
$$
   e^{-I_r}(x,x) \sim r^{-p}\sum_{i=0}^\infty r^i E_i(x,x)
$$
where $E_i$ are endomorphisms of $\fs\otimes \Lambda\otimes
\Lambda^{p,0}T^*X\otimes{\Cal E}$.
\endproclaim

\demo{Proof}
Consider the operator
$$
\align
   J_r &= -rg^{ij}(\delta_i + {1\over 4}\Gamma_{iab} c(e_a\wedge e_b) +
      A_i) \times (\partial_j + {1\over 4}\Gamma_{jab}c(e_a\wedge e_b) + A_j) \\
   &\quad + rg^{ij}\Gamma_{ij}^k (\partial_k + {1\over 4} \Gamma_{kab}
      c(e_a\wedge e_b) + A_k) + {1\over 4}rk - {1\over 2}c(e_i\wedge e_j)
      L_{ij}.
\endalign
$$
Since $J_r$ has no singular terms as $r\to 0$, it has a well known
asymptotic expansion, as $r\to 0$ with $p=n$.

We can construct $\exp(-I_r)$ as a perturbation of $\exp(-J_r)$, using
Duhamel's principle.  More precisely,
$$
   \exp(-I_r) = \exp(-J_r) + \sum_{k=1}^\infty
      \underbrace{e^{-J_r}(J_r-I_r)e^{-J_r}\dots e^{-J_r}}_{k \text{terms}}
$$
Each coefficient in the difference $J_r-I_r$ contains at least one term
which is exterior multiplication by $f_\alpha$.  Therefore the infinite
series on the right hand side collapses to a finite number of terms.  The
proposition then follows from the asymptotic expansion for
$\exp(-J_r)(x,x)$.
\enddemo\hfill$\square$

\def\Tr{\operatorname{Tr}}
\def\ch{\operatorname{ch}}
Let $R^B$ denote the curvature of the holomorphic Hermitian connection and
$R^L$ denote the curvature of the Levi-Civita connection.  Let $\hat{A}$
denote $\hat{A}$-invariant polynomial and $\ch$ the Chern character
invariant polynomial.  Then
$$
   \hat{A}(R^{-B})\ch(\Tr(R^L) \ch (\Lambda^{p,0}R^L) \in \Lambda^*T^*X.
$$
The goal is to prove the following decoupling result in the adiabatic limit. It
resembles the local index theorem for almost K\"ahler manifolds by Bismut
\cite{Bi}
(he calls them non-K\"ahler manifolds). However, we use instead the
techniques of
the proofs in \cite{BGV}, \cite{Ge} and \cite{D} of the local index theorem
for families. In particular, we borrow a local conjugation trick due to
Donnelly \cite{D},
which is adjusted to our situation.

\proclaim{6.6. Theorem (Adiabatic decoupling)} Let $(X,g)$ be an almost
K\"ahler manifold.
In the notation above, one has the following
decoupling result in the adiabatic limit
$$
  \lim_{r\to 0} \Tr_\tau^s(Z_u e^{-I_r}(x,x)) = \Tr_\tau(Z_u)(x) \
      [ \hat{A}(R^{-B})\ch(\Tr(R^L)) \ch (\Lambda^{p,0}R^L)
      ]^{\max}_x\in\Lambda^{2n}T^*_xX $$
for all $x\in X$.
\endproclaim

\demo{Proof}
We first consider the corresponding problem on $\Bbb{R}^{2n}$, using the
exponential map.  Let $\bar{I}_r$ denote the operator on $\Bbb{R}^{2n}$, whose
expressions agrees with the local coordinate expression for $I_r$ near $p$,
where $p$ is identified with the origin in $\Bbb{R}^{2n}$.

Consider the heat equation on $\Bbb{R}^{2n}$,
$$
\gather
   (\partial_t + \bar{I}_r) g(x,t) = 0 \\
   g(x,0) = g(x).
\endgather
$$
Then one has
\enddemo\hfill$\square$

\proclaim{6.7. Proposition}
There is a unique fundamental solution $e^{-t\bar{I}_r}(x,y)$ which
satisfies the decay estimate
$$
   \bigl| e^{-t\bar{I}_r}(x,y)\bigr| \le c_1t^{-n} e^{-{c_2|x-y|^2\over t}}
$$
as $t\to 0$, with similar estimates for the derivatives in $x,y,t$.
\endproclaim

\demo{Proof} The proof is standard, as in 3.10.
\enddemo\hfill$\square$

By Duhamel's principle applied in a small enough normal coordinate neighborhood,
there is a positive constant $c$ such that
$$
   e^{-\bar{I}_r}(0,0) = e^{-I_r}(x,x) + O(e^{-c/r})\quad \text{as }r\to 0.
$$
Therefore $$\lim_{r\to0} \Tr_\tau^s(Z_u e^{-I_r}(x,x)) = \lim_{r\to0}
\Tr_\tau^s(Z_u e^{-\bar{I}_r}(0,0))\tag29$$
and it suffices to compute the right hand side of (29).  This is done using
Getzler's scaling idea \cite{Ge},  $x\to\epsilon x,\ t\to\epsilon^2t\
e_i\to \epsilon^{-1}e_i$.  Then Clifford
multiplication scales as $c_\epsilon(\cdot) = e(\cdot)  + \epsilon^2
i(\cdot)$, where
$e(\cdot)$ denotes exterior multiplication by the covector $\cdot$ and
$i(\cdot)$ denotes contraction by the dual vector.
$$
\align
   \bar{I}_\epsilon &= -rg^{ij}(\epsilon x) \Big(\partial_i +
      {\epsilon^{-1} \over 4}\Gamma_{iab}(\epsilon x)
      c_\epsilon(e_a\wedge e_b) + \epsilon A_i(\epsilon x) +
{\epsilon^{-1}\over 2\sqrt{r}}
      c_\epsilon(S_{il\alpha}(\epsilon x) e_i) e(f_\alpha)\\
   &\quad  + {\epsilon^{-1}\over 4r} S_{i\beta\gamma}{(\epsilon x)}
      c(f_\beta\wedge f_\beta)\Big) \times\Big(\partial_j +
      {\epsilon^{-1}\over 4} \Gamma_{jab}
   (\epsilon    x)  c_\epsilon(e_a\wedge e_b) + \epsilon A_j(\epsilon x)\\
   &\quad\phantom{\bigl(} + {\epsilon^{-1}\over 2\sqrt{r}} c_\epsilon
      (S_{jl\alpha} (\epsilon x) e_l) e(f_\alpha) + {\epsilon^{-1}\over 4r}
      S_{j\beta\gamma}(\epsilon x) e(f_\beta\wedge f_\gamma)\Big)
      + rg^{ij}(\epsilon x) \Gamma_{ij}^k(\epsilon x) \Big( \epsilon
      \partial_k\\
   &\quad  + {1\over4}\Gamma_{kab}(\epsilon x) c_\epsilon(e_a\wedge
      e_b) +\epsilon^2A_k(\epsilon x) + {1\over 2\sqrt{r}} S_{kl\alpha}
      (\epsilon x) c_\epsilon(e_i) e(f_\alpha)
      + {1\over 4r} S_{k\beta\gamma}(\epsilon x) e(f_\beta\wedge
      f_\gamma)\Big)  \\
   &\quad + {\epsilon^2\over 4} rk(\epsilon x)
   -{r\over 2} c_\epsilon(e_i\wedge e_j) L_{ij}(\epsilon x) -
      {1\over 2}f_\alpha\wedge f_\beta \wedge L_{\alpha\beta}(\epsilon x)
      -\sqrt{r} c_\epsilon(e_i) f_\alpha L_{i\alpha}(\epsilon x).
\endalign
$$
The asymptotic expansion in $r$ as in Propositions 6.5 and 6.7, for
$e^{-\bar{I}_r}(0,0)$ yields an asymptotic expansion in $\epsilon$ for
$e^{-\bar{I}_\epsilon}(0,0)$ and one
has
$$
   \lim_{r\to 0} \Tr_\tau^s\bigl(Z_u e^{-\bar{I}_r}(0,0)\bigr)
      = \lim_{\epsilon\to
      0}\Tr_\tau^s\bigl(Z_ue^{-\bar{I}_\epsilon}(0,0)\bigr)\tag30
$$
That is, if either limit exists, then both exist and are equal.

However, in the limit as $\epsilon\to 0$, there are singularities in the
coefficients of $S$ tensor in the expression for $\bar{I}_\epsilon$ and
one cannot immediately apply Getzler's theorem. Therefore one first makes
the following local conjugation trick, as in Donnelly \cite{D}.

Define the expression
$$
   h(x,\epsilon,r) = \exp\Bigl( {\epsilon^{-1}\over 2\sqrt{r}}
      S_{il\alpha}(0)x_ie_l\wedge f_\alpha + {\epsilon^{-1}\over 4r}
      S_{i\beta\gamma}(0) x_i f_\beta\wedge f_\gamma\Bigr).
$$
Note that $h(x,\epsilon,r)$ has polynomial growth in $x$, since its expression
contains exterior multiplication. We claim that if the operator
$\bar{I}_\epsilon$ is conjugated by $h$, then
the resulting operator is {\it not} singular as $\epsilon\to 0$.  More
precisely,
$$
\align
   J_\epsilon &= h\bar{I}_\epsilon h^{-1} \\
   &= rg^{ij}(\epsilon x) \Bigl(\partial_i + {\epsilon^{-1}\over4}\Gamma_{iab}
      (\epsilon x) e_a\wedge e_b + {\epsilon^{-1}\over 2\sqrt{r}}
      \bigl(S_{il\alpha} (\epsilon x) - S_{il\alpha}(0)\bigr) e_l\wedge
      f_\alpha+ {\epsilon^{-1}\over 4r} \bigl(S_{i\beta\gamma}(\epsilon x)\\
   &\quad   -
      S_{i\beta\gamma}(0)\bigr) f_\beta\wedge f_\gamma -{1\over 4r}
S_{il\alpha}(0) S_{kl\beta}(0)
      x_k f_\alpha\wedge f_\beta\Bigr) \times\Bigl( \partial_j +
{\epsilon^{-1}\over 4}\Gamma_{jab}(\epsilon x)
      e_a\wedge e_b\\
   &\quad + {\epsilon^{-1}\over 2\sqrt{r}}\bigl(S_{jl\alpha}
      (\epsilon x) - S_{i\beta\gamma}(0)\bigr) e_l\wedge f_\alpha +
{\epsilon^{-1}\over 4r}
      \bigl(S_{j\beta\gamma}(\epsilon x) - S_{j\beta\gamma}(0)\bigr)
      f_\beta\wedge f_\gamma -{1\over 4r} S_{jl\alpha}(0) S_{kl\beta}(0) x_k\\
   &\quad  f_\alpha\wedge f_\beta\Bigr)
   - {1\over 2} r e_i\wedge e_j L_{ij}(\epsilon x) -{1\over 2}
f_\alpha\wedge f_\beta L_{\alpha\beta}(\epsilon x) -
      \sqrt{r}e_i \wedge f_\alpha L_{i\alpha} (\epsilon x) +
R(x,\epsilon).\tag31
\endalign
$$
Here $R(x,\epsilon)$ denotes the terms which vanish as $\epsilon\to 0$, and
which therefore do not contribute to the limit. Clearly there are no singular
terms in $J_\epsilon$ as $\epsilon\to 0$.

A fundamental solution for the heat equation for $J_\epsilon$ can be
obtained by conjugating the one for $\bar{I}_\epsilon$, that is
$$
e^{-tJ_\epsilon}(x,y) = h(x, \epsilon, r) e^{-t\bar{I}_\epsilon}(x,y)
h^{-1}(y, \epsilon, r).
$$
The right hand side satisfies the heat equation
$(\partial_t+J_\epsilon)g(x,t) = 0,\ \ g(x,0) = \delta_x$.  Since $h(0)=1$,
one has $\forall\epsilon>0$,
$$
   \Tr_\tau^s\bigl(Z_ue^{-\bar{I}_\epsilon}(0,0)\bigr) =
      \Tr_\tau^s\bigl(Z_u e^{-J_\epsilon}(0,0)\bigr).\tag32
$$
Therefore it suffices to compute the limit as $\epsilon\to 0$ of the right
hand side of (31).

Using the following Taylor expansions in a normal coordinate neighborhood,
$$
\align
   \Gamma_{iab}(\epsilon x) &= -{1\over 2} R_{ijab}(0)\epsilon x_j +
      R(x,\epsilon^2) \\
   S_{il\alpha}(\epsilon x) &= S_{il\alpha}(0) + S_{il\alpha,j}(0)
      \epsilon x_j + R(x,\epsilon^2) \\
   S_{i\beta\gamma}(\epsilon x) &= S_{i\beta\gamma}(0) +
      S_{i\beta\gamma,j}(0) \epsilon x_j + R(x,\epsilon^2)
\endalign
$$
one sees that
$$
   J_0 = \lim_{\epsilon\to 0} J_\epsilon = -r\sum_i(\partial_i -
      {1\over4}B_{ij}x_j)^2 + r{\Cal L}
$$
where
$$
\align
   B_{ij} &= {1\over 2} R_{ijab}(0) e_a\wedge e_b -{2\over\sqrt{r}}
      S_{il\alpha,j}(0) e_l\wedge f_\alpha \\
   &\quad -{1\over r}\bigl(S_{i\beta\gamma,j}(0) - S_{il\beta}(0)
      S_{jl\gamma}(0)\bigr) f_\beta \wedge f_\gamma
\endalign
$$
and
$$
   {\Cal L} = {1\over2} L_{ij}(0) e_i\wedge e_j + {1\over\sqrt{r}}
      L_{i\alpha}(0) e_i\wedge f_\alpha + {1\over2r} L_{\alpha\beta}(0)
      f_\alpha\wedge f_\beta.
$$
Using Mehler's formula (cf. \cite{Ge}), one can obtain an explicit fundamental
solution $e^{-sJ_0}(x,y)$.  First decompose $B$ into its symmetric and skew
symmetric parts, that is
$B=C+D$ where $C={1\over2}(B+B^t)$ and $D={1\over2}(B-B^t)$, where $B,C,D$
are matrices of 2-forms.  Then

$$
   e^{-J_0}(x,0) = (4\pi r)^{-n/2}\hat{A}(r D) e^{{x^tCx\over 8}}
      \times\exp\Bigl(r{\Cal L} - {1\over 4r} x^t
      \Bigl({r D/2 \over \tanh(r D/2)}\Bigr) x\Bigr)
$$
Now $\lim_{\epsilon\to 0} e^{-J_\epsilon}(0,0) = e^{-J_0}(0,0)$.  Therefore
$$
   \lim_{\epsilon\to 0} \Tr_\tau^s(Z_u e^{-J_\epsilon}(0,0)) = \bigl(
      {2\over i}\bigr)^{n/2} (4\pi r)^{-n/2} \Tr_\tau(Z_u)(0) \bigl[
      \hat{A}(rD)\ch (r{\Cal L}) \bigr]^{\max}\tag33
$$
Here $D=R^{-B}(0)$ and ${\Cal L} = \Tr(R^L(0)) + \Lambda^{p,0} R^L(0)$.
Using (29), (30), (32) and (33), one completes the proof of Theorem 6.6.
$\square$
\vskip .2in

\proclaim{6.8. Theorem}
Let $\Cal{E}$ be a flat Hilbertian bundle of $D$-class,
over an almost K\"ahler manifold $(X,g)$ and let $h$, $h'$
be Hermitian metrics on $\Cal{E}$ such that $h=h'$ in
$\det(\Cal{E})$.  Then
$$
\rho_\Cal{E}^p(g,h)=\rho_\Cal{E}^p(g,h')\in\det(H^{p,*}(X,\Cal{E}))
$$
\endproclaim

\demo{Proof}
Since $h=h'$ in $\det(\Cal{E})$, there is a positive, self-adjoint
bundle map $A:\Cal{E}\rightarrow\Cal{E}$ satisfying
$$
\gather
h(A v,w)=h'(v,w) \qquad \forall v,w\in\Cal{E} \\
\text{and} \qquad \Det_\tau(A)=1.
\endgather
$$
By Corollary 6.2, there is a smooth 1-parameter family of
positive, self-adjoint
bundle maps $u\to A_u:\Cal{E}\rightarrow\Cal{E}$ joining
$A$ to the identity and satisfying
$$
\Det_\tau(A_u)=1. \tag34
$$
for all $u\in (-\epsilon,1+\epsilon)$. Here $A_0 = I$ and $A_1 = A$.
Let $u\rightarrow h_u$ be a smooth family of Hermitian metrics on
$\Cal{E}$ defined by
$$
h(A_u v,w)=h_u(v,w) \qquad \forall v,w\in\Cal{E}.
$$
Then $h_0=h$, $h_1=h'$ in $\E$ and $h=h_u$ in
$\det(\Cal{E})$ for all $u\in (-\epsilon,1+\epsilon)$
by (72).
Note that by differentiating (34), one has
$$
0 = {\partial \over\partial u}\Det_\tau(A_u) = \Tr_\tau(Z_u)
\tag35
$$
where $Z_u=A_u^{-1}\dot{A_u}$.

We wish to compute ${\partial \over\partial u}\rho_\Cal{E}^p (g,h_u)$.
By Theorem 4.5, one has
$$
{\partial \over\partial u}\rho_\Cal{E}^p(g,h_u)=c_\E^p(g,h_u)
   \rho_\Cal{E}^p(g,h_u)
$$
By Theorem 6.6 and (35), one sees that
$$
\lim_{t\rightarrow 0} \Tr_\tau^s(Z_u e^{-t\square(u)})=0.\tag36
$$
By the small time asymptotic expansion of the heat kernel,
one has
$$
\align
\lim_{t\rightarrow 0} \Tr_\tau^s(Z_u e^{-t\square(u)})&=
   \sum_{q=0}^n(-1)^q m_{n,p,q}(u)\\
&= c_\E^p (g,h_u)\tag37
\endalign
$$
Therefore by (36) and (37), one has $c_\E^p(g,h_u)=0$, that is,
$$
{\partial \over\partial u}\rho_\Cal{E}^p(g,h_u)=0.
$$
\enddemo\hfill$\square$

\subheading{6.9. Remarks} Theorem 6.8 says that on an almost
K\"ahler manifold $(X,g)$, the
holomorphic $L^2$ torsion $\rho_\Cal{E}^p(g,h)$ depends only on the
equivalence class of the Hermitian metric $h$ in $\det(\Cal{E})$.
We however do not believe that the almost K\"ahler
hypothesis in Theorem 6.8 is necessary. However, we use the techniques
of the proof of the local index theorem, and the
situation to date is that the local index theorem for the operator
$\overline\partial +\overline\partial^*$
has not yet been established for a general Hermitian manifold.

\subheading{6.10} Let $\Cal{E}$ and $\F$ be two flat Hilbertian bundles
of $D$-class over over an almost K\"ahler manifold $(X,g)$
and $\varphi:\det(\Cal{E})\rightarrow \det(\F)$ be an isomorphism of the
determinant line bundles.  Then using the theorem above, we will
construct a canonical isomorphism between determinant lines
$$
\gather
\widehat{\varphi}^p:\det H^{p,*}(X,\Cal{E})\rightarrow \det H^{p,*}(X,\F)\\
\widehat{\varphi}^p(\lambda \rho_\Cal{E}^p(g,h))=\lambda \rho_{\F}^p(g,h'),
\quad \lambda\in\Bbb{R}
\endgather
$$
where $h$ and $h'$ are Hermitian metrics on $\Cal{E}$ and $\F$
respectively, such that $\varphi(h)=h'$ in $\det(\F)$.
Then $\widehat{\varphi}$ is called a {\it correspondence} between determinant
line bundles. It is well defined by Theorem 6.8 and Remarks 6.9.
We next state some obvious properties of correspondences.

\proclaim{6.11. Proposition}
Let $\Cal{E}$  be a flat Hilbertian bundle of $D$-class over
over an almost K\"ahler manifold $(X,g)$ and
$\varphi:\det(\Cal{E})\rightarrow \det(\E)$ be the identity
map. Then $$\widehat{\varphi}^p= identity$$.

Let $\Cal{E}, \F$ and $\G$ be flat Hilbertian bundles of $D$-class over
over an almost K\"ahler manifold $(X,g)$ and
$\varphi:\det(\Cal{E})\rightarrow \det(\F)$,
$\psi:\det(\Cal{F})\rightarrow \det(\G)$
be isomorphisms of the determinant line bundles. Then the composition
satisfies
$$
\widehat{\varphi {\small\text o}\psi}^p = \widehat{\varphi}^p {\small\text o}
\widehat{\psi}^p.$$
\endproclaim

We next prove one of the main results in the paper.

\proclaim{6.12. Theorem}
Let $\Cal{E}$ and $\F$ be two flat Hilbertian bundles of $D$-class over
over an almost K\"ahler manifold $(X,g)$ and
$\varphi:\det(\Cal{E})\rightarrow \det(\F)$ be an isomorphism of the
corresponding determinant line bundles.  Consider smooth 1-parameter
families of almost K\"ahler metrics $g_u$ on $X$ and Hermitian metrics
$h_{1,u}$ on $\Cal{E}$, where $u$ varies in an internal $(-\epsilon,\epsilon)$.
Choose a smooth family of Hermitian metrics $h_{2,u}$ on $\F$ in such
a way that $\varphi(h_{1,u})=h_{2,u}$ in $\det(\F)$.  Then the relative
holomorphic torsion
$$
\rho_\varphi^p(u)=\rho_\Cal{E}^p(g_u,h_{1,u})\otimes\rho_{\F}^p(g_u,h_{2,u})
^{-1}
\in\det H^{p,*}(X,\Cal{E})\otimes\det H^{p,*}(X,\F)^{-1}
$$
is a smooth function of $u$ and satisfies ${\partial \over\partial u}
\rho_\varphi(u)=0$.  That is, the relative holomorphic
$L^2$ torsion $\rho_\varphi^p$
is independant of the choices of metrics on $X$, $\Cal{E}$
and $\F$ which are needed to define it.
\endproclaim

\demo{Proof}From the data in the theorem,
 one can define a correspondence as in 6.10,
$$
\hat{\varphi}^p:\det(H^{p,*}(X,\Cal{E}))\rightarrow \det(H^{p,*}(X,\F))
$$
which is an isomorphism of determinant lines.  It is defined as
$$
\hat{\varphi}^p(\lambda \rho_\Cal{E}^p(g_u,h_{1,u}))=\lambda\rho_{\F}^p
(g_u,h_{2,u})
\tag{38}
$$
for $\lambda\in\Bbb{R}$ and $u\in(-\epsilon,\epsilon)$.  Therefore
using theorem 4.5 and (38) above, one has
$$
\align
{\partial \over\partial u}\hat{\varphi}^p(\rho_\Cal{E}^p(g_u,h_{1,u}))&=
   \hat{\varphi}^p({\partial \over\partial u}\rho_\Cal{E}^p(g_u,h_{1,u}))\\
&=c_\Cal{E}(g_u,h_{1,u})\hat{\varphi}^p(\rho_\Cal{E}^p(g_u,h_{1,u}))\\
&=c_\Cal{E}(g_u,h_{1,u})\rho_{\F}^p(g_u,h_{2,u}).
\tag{39}
\endalign
$$
But by differentiating equation (38) above, one has
$$
\align
{\partial \over\partial u}\hat{\varphi}(\rho_\Cal{E}^p(g_u,h_{1,u}))
   &={\partial \over\partial u}\rho_{\F}^p(g_u,h_{2,u})\\
&=c_{\F}(g_u,h_{2,u})\rho_{\F}^p(g_u,h_{2,u})
\tag{40}
\endalign
$$
Equating (39) and (40), one has
$$
c_\Cal{E}(g_u,h_{1,u})=c_{\F}(g_u,h_{2,u})\tag41
$$
Then the relative holomorphic $L^2$ torsion
$$
\gather
\rho_\varphi^p\in\det H^{p,*}(X,\Cal{E})\otimes(\det H^{p,*}(X,\F))^{-1}\\
\rho_\varphi^p(u)=\rho_\Cal{E}^p(g_u,h_{1,u})\otimes(\rho_{\F}^p(g_u,h_{2,u}
))^{-1}
\endgather
$$
satisfies
$$
\align
{\partial \over\partial u}\rho_\varphi^p(u) =& \left({\partial \over\partial u}
   \rho_\Cal{E}^p(g_u,h_{1,u})\right)\otimes\left(\rho_{\F}^p
   (g_u,h_{2,u})\right)^{-1}\\
&-\rho_\Cal{E}^p(g_u,h_{1,u})\otimes{\partial \over\partial u}
   \rho_{\F}^p(g_u,h_{2,u})\cdot \rho_{\F}^p(g_u,h_{2,u})^{-2}\\
=& c_\Cal{E}(g_u,h_{1,u})\rho_\varphi^p(u)
-c_{\F}(g_u,h_{2,u}) \rho_\varphi^p(u)\\
=& 0
\endalign
$$
using Theorem 4.5 and (41) above.  This proves the theorem.
\enddemo\hfill$\square$

\heading{\bf \S 7. Calculations} \endheading

In this section, we  calculate the holomorphic $L^2$ torsion for K\"ahler
locally symmetric
spaces. We will restrict ourselves to the special case of the Hilbert
$({\Cal U}(\Gamma)-\Gamma)$-bimodule $\ell^2(\Gamma)$,
where $\Gamma$ is a countable discrete group. Let $\E\to X$ denote the
associated flat Hilbert ${\Cal U}(\Gamma)$-bundle over the compact complex
manifold
$X$. Then it is well known that
the Hilbert ${\Cal U}(\Gamma)$-complexes
$\left(\Omega_{(2)}^{\bullet,\bullet}({X},\E),\nabla''\right)$ and
$\left(\Omega_{(2)}^{\bullet,\bullet}(\widetilde{X}),
\bar{\partial}\right)$ are canonically
isomorphic, where $\Gamma\to{\widetilde X}\to X$ denotes the universal
covering space
of $X$ with structure group $\Gamma$. We will denote the
$\bar{\partial}$-Laplacian
acting on $\Omega_{(2)}^{p,q}(\widetilde{X})$ by $\square_{{p,q}}$.

Firstly, we will discuss the $D$-class condition in this case. Let $X$ be a
K\"ahler hyperbolic manifold. Recall that this means that $X$ is a K\"ahler
manifold with K\"ahler form $\omega$, which has the property that $p^*(\omega)
= d\eta$, where $\Gamma\to{\widetilde X}\to X$ denotes the universal cover
of $X$
and $\eta$ is a bounded 1-form on ${\widetilde X}$. Any Riemannian manifold of
negative sectional curvature,
which also supports a K\"ahler metric, is a K\"ahler hyperbolic manifold.
Note that
the K\"ahler metric is not assumed to be compatible with the Riemannian
metric of
negative sectional curvature. Then Gromov \cite{G} proved that on the universal
cover of a K\"ahler hyperbolic manifold, the Laplacian $\square_{{p,q}}$
has a spectral gap at zero
on all $L^2$ differential forms. Therefore it follows that the associated
flat bundle
$\E\to X$ is of $D$-class. By a vanishing theorem of Gromov \cite{G}
for the $L^2$ Dolbeault cohomology of the universal cover, one has
$$
H_{(2)}^{p,q}(\widetilde X) = 0
$$
unless $p+q = n$, where $n$ denotes the complex dimension of $X$.

In particular, let
$G$ be a connected semisimple Lie group, and $K$ be a maximal compact
subgroup such that
$G/K$ carries an invariant complex structure, and let $\Gamma$ be a
torsion-free uniform lattice
in $G$. Then it is known that $\Gamma\backslash G/K$ is a K\"ahler
hyperbolic manifold
(cf. \cite{BW})
and therefore the canonical flat Hilbert bundle $\E\to X$ is of $D$-class.
In this K\"ahler
metric, the Laplacian $\square_{{p,q}}$ is $G$-invariant,
so it follows that the theta function
$$
\theta_{p,q}(t) = C_{p,q}(t) vol(\Gamma\backslash G/K)
$$
is proportional to the volume of $\Gamma\backslash G/K$. Here $C_{p,q}(t)$
depends only
on $t$ and on $G$ and $K$, but {\it not} on $\Gamma$. It follows that the
zeta function
$\zeta_{p,q}(s, \lambda, \E)$ is also proportional to the volume of
$\Gamma\backslash G/K$. Therefore
the holomorphic $L^2$ torsion is given by
$$
\rho_\Cal{E}^p  = e^{C_p vol(\Gamma\backslash G/K)} \rho^{\prime p} \in
\det \left( H_{(2)}^{p,n-p}(G/K)\right)^{(-1)^{n-p}}
$$
where we have used the vanishing theorem of Gromov. Here $C_{p}$ is a
constant that depends only
on $G$ and $K$, but {\it not} on $\Gamma$. Using representation theory, as
for instance in \cite{M},
\cite{L} and \cite{Fr}, it is possible to determine $C_p$ explicitly. This
will be done elsewhere.
Using the proportionality principle again, one sees that the Euler
characteristic of
$\Gamma\backslash G/K$ is proportional to its volume.
By  a theorem of Gromov \cite{G}, the Euler characteristic of
$\Gamma\backslash G/K$ is non-zero. Therefore we can also express the
holomorphic $L^2$ torsion as
$$
\rho_\Cal{E}^p  = e^{C'_p \chi(\Gamma\backslash G/K)} \rho^{\prime p} \in
\det \left( H_{(2)}^{p,n-p}(G/K)\right)^{(-1)^{n-p}}
$$
where $\chi(\Gamma\backslash G/K)$ denotes the Euler characteristic of
$\Gamma\backslash G/K$,
and $C'_{p}$ is a constant that depends only
on $G$ and $K$, but {\it not} on $\Gamma$. This discussion is summarized in
the following
proposition.

\proclaim{7.1. Proposition} In the notation above, the holomorphic $L^2$
torsion of the
semisimple locally symmetric space $\Gamma\backslash G/K$, which is assumed to
carry an invariant complex structure, is given by
$$
\rho_\Cal{E}^p  = e^{C_p vol(\Gamma\backslash G/K)} \rho^{\prime p} \in
\det \left( H_{(2)}^{p,n-p}(G/K)\right)^{(-1)^{n-p}}
$$
Here $C_{p}$ is a constant that depends only
on $G$ and $K$, but {\it not} on $\Gamma$. Equivalently, the holomorphic
$L^2$ torsion of
$\Gamma\backslash G/K$ is given as
$$
\rho_\Cal{E}^p = e^{C'_p \chi(\Gamma\backslash G/K)} \rho^{\prime p} \in
\det \left( H_{(2)}^{p,n-p}(G/K)\right)^{(-1)^{n-p}}
$$
where $\chi(\Gamma\backslash G/K)$ denotes the Euler characteristic of
$\Gamma\backslash G/K$,
and $C'_{p}$ is a constant that depends only
on $G$ and $K$, but {\it not} on $\Gamma$.
\endproclaim

We will now compute the
holomorphic $L^2$ torsion for a Riemann surface, which is a special case of
the theorem
above, and we will show that the constants $C_p$ and $C'_p$ are not zero.

Let $X$ be a closed Riemann surface of genus $g$, which is greater than 1,
which can be realised as
a compact quotient complex hyperbolic space $\Bbb H$ of complex dimension 1,
by the torsion-free discrete group $\Gamma$.
Recall that
$$
{\square}_{0,1} = \frac{1}{2} \Delta_1
$$
acting on the subspace $\Omega_{(2)}^{0,1}(\Bbb H)$.
Also,
$$
{\square}_{1,0} = \frac{1}{2} \Delta_1
$$
acting on the subspace $\Omega_{(2)}^{1,0}(\Bbb H)$.
The $\star$ operator intertwines the operators ${\square}_{0,1}$ and
${\square}_{1,0}$, showing in particular that they are isospectral.
So using
$$
\Omega_{(2)}^{1}(\Bbb H) =
\Omega_{(2)}^{1,0}(\Bbb H) \oplus \Omega_{(2)}^{0,1}(\Bbb H),
$$
we see that in order to calculate the von Neumann determinant of the operator
${\square}_{1,0}$, it suffices to
first scale the metric $g\to 2g$ and then calculate the square root
von Neumann determinant of the Laplacian $\Delta_1$. However, this is easily
seen to be equal to the von Neumann determinant of the Laplacian $\Delta_0$
acting on $L^2$ functions on the hyperbolic disk. Recall that the von Neumann
determinant of the operator $A$ is by definition $e^{-\zeta'_A(0)}$, where
$\zeta'_A(s)$ denotes the zeta function of the operator $A$.

Using the work of Randol \cite{R}, one obtains the following expression for the
meromorphic continuation of the zeta function of $\Delta_0$ to the
half-plane $\Re(s)<1$
$$
\zeta_0(s,0,\Cal{E}) = (g-1) \frac{\pi}{(s-1)} \int_0^\infty \left( {1\over
4} +r^2\right)^{1-s}
sech^2(\pi r) dr.
$$
It follows that
$$
\align
\zeta'_0(0,0,\Cal{E}) & = \lim_{s\to 0} \left(\zeta_0(s,0,\Cal{E}) -
\zeta_0(0,0,\Cal{E}) \right)\Gamma(s)\\
& =  (g-1){\pi} \int_0^\infty \left( {1\over 4} +r^2\right)
sech^2(\pi r) \left( -1 + \log\left( {1\over 4} +r^2\right)\right)dr.
\endalign
$$
A numerical approximation for the last integral shows that
$\zeta'(0,0,\Cal{E}) \sim -0.677 (g-1) $. We can summarize
the discussion in the following proposition.

\proclaim{7.2. Proposition} In the notation above, the holomorphic $L^2$
torsion of a
compact Riemann surface $X = \Gamma \backslash \Bbb H$ of genus $g$, is given by
$$
\rho_\Cal{E}^0 = e^{C (g-1)} \rho^{\prime 0} \in
\det \left( H_{(2)}^{0,1}(\Bbb H)\right)^{(-1)}\tag42
$$
Here $C = {\pi\over 2} \int_0^\infty \left( {1\over 4} +r^2\right)
sech^2(\pi r) \left( -1 + \log\left( {1\over 4} +r^2\right)\right)dr$
is a constant that depends only on $\Bbb H$, but {\it not} on $\Gamma$.
$C$ is approximately $-0.338$, and in particular, it is not equal to zero.
Also,
$$
\rho_\Cal{E}^1 = e^{- C (g-1)} \rho^{\prime 1} \in
\det \left( H_{(2)}^{1,0}(\Bbb H)\right)^{(-1)}
$$
where the constant $C$ is as in (42).

\endproclaim

\subheading{Acknowledgement} We thank John Phillips for his cleaner proof
of Lemma 6.9.

\enddocument